\documentclass{elsart}

\usepackage{graphicx}
\usepackage{psfrag}

\usepackage{amssymb}
\usepackage{amsmath}
\usepackage{dsfont}
\usepackage{rotating}

\newcommand{\Path}{\boldsymbol{\Psi}}

\newcommand{\tens}[1]{#1}
\newcommand{\K}{\tens{K}}
\newcommand{\M}{\tens{M}}
\newcommand{\D}{\mathbf{D}}

\renewcommand{\S}{\mathbf{S}} 
\newcommand{\Q}{\mathbf{Q}} 
\newcommand{\q}{\mathbf{q}} 
\newcommand{\f}{\mathbf{f}} 
\newcommand{\g}{\mathbf{g}} 
\newcommand{\A}{\mathbf{A}} 
\newcommand{\B}{\mathbf{B}} 
\newcommand{\F}{\mathbf{F}} 
\newcommand{\G}{\mathbf{G}} 
\newcommand{\x}{\mathbf{x}} 
\newcommand{\xxi}{\boldsymbol{\xi}} 
\newcommand{\X}{\mathbf{X}} 
\renewcommand{\P}{\mathbf{P}} 
\newcommand{\w}{\mathbf{w}}

\newcommand{\halb}{\frac{1}{2}}

\newcommand{\ar}{\phi_1\rho_1}
\newcommand{\arr}{\phi_2\rho_2}

\newcommand{\ub}{\textbf{u}_\textbf{1}}

\begin{document}

\begin{frontmatter}



\title{High--Order Unstructured Lagrangian One--Step WENO Finite Volume Schemes for Non--conservative Hyperbolic Systems: Applications to Compressible Multi--Phase Flows}


\author[UNITN]{Michael Dumbser}
\ead{michael.dumbser@unitn.it}
\author[UNITN]{Walter Boscheri}
\ead{walter.boscheri@unitn.it}
\address[UNITN]{Laboratory of Applied Mathematics \\ Department of Civil, Environmental and Mechanical Engineering \\ University of Trento, Via Mesiano 77, I-38123 Trento, Italy}

\begin{abstract}
In this article we present the first better than second order accurate \textit{unstructured Lagrangian--type} 
one--step WENO finite volume scheme for the solution of hyperbolic partial differential equations with 
\textit{non--conservative} products. The method achieves high order of accuracy in space together with 
essentially non--oscillatory behaviour using a nonlinear WENO reconstruction operator on unstructured triangular meshes. 
High order accuracy in time is obtained via a local Lagrangian space--time Galerkin predictor method that evolves 
the spatial reconstruction polynomials in time within each element. 
The final one--step finite volume scheme is derived by integration over a moving space--time control volume, where
the non--conservative products are treated by a path--conservative approach that defines the jump terms on the element
boundaries. The entire method is formulated as an Arbitrary--Lagrangian--Eulerian (ALE) method, where the mesh velocity 
can be chosen independently of the fluid velocity. 

The new scheme is applied to the full seven--equation Baer--Nunziato model of compressible multi--phase flows in two 
space dimensions. The use of a Lagrangian approach allows an excellent resolution of the solid contact and the resolution 
of jumps in the volume fraction. The high order of accuracy of the scheme in space and time is confirmed via a numerical
convergence study. Finally, the proposed method is also applied to a reduced version of the compressible Baer--Nunziato model 
for the simulation of free surface water waves in moving domains. In particular, the phenomenon of sloshing is studied in 
a moving water tank and comparisons with experimental data are provided.  

\end{abstract}
\end{frontmatter}


\section{Introduction}
\label{sec.introduction}

Multi--phase flow problems, such as liquid-vapour and solid-gas flows are encountered in numerous natural processes,  
such as avalanches, meteorological flows with cloud formation, volcano explosions, sediment transport in rivers and 
on the coast, granular flows in landslides, etc., as well as in many industrial applications, e.g., in aerospace engineering,  
automotive industry, petroleum and chemical process engineering, nuclear reactor safety, paper and food manufacturing and 
renewable energy production. Most of the industrial applications are concerned with 
\textit{compressible} multi--phase flows as they appear for example in combustion processes of liquid and solid fuels in 
car, aircraft and rocket engines, but also in solid bio--mass combustion processes. Already the mathematical description of
such flows is quite complex and up to now there is no universally agreed model for such flows. One wide--spread model is 
the Baer--Nunziato model for compressible two-phase flow, which has been introduced by Baer and Nunziato in \cite{BaerNunziato1986} 
for detonation waves in solid--gas combustion processes and which has been further extended by Saurel and Abgrall to liquid--gas 
flows in \cite{SaurelAbgrall}. It is therefore also often called in literature the Saurel--Abgrall model. The main difference between 
the original Baer--Nunziato model and the Saurel--Abgrall model is the definition of the pressure and velocity at the interface. 
In the present paper, we will use the original choice of Baer--Nunziato, which has also been used in several papers about the exact 
solution of the Riemann--Problem of the Baer--Nunziato model, see \cite{AndrianovWarnecke,DeledicquePapalexandris,Schwendeman}.  
A reduced five--equation model has been proposed in \cite{MurroneGuillard} and approximate Riemann solvers of Baer--Nunziato--type 
models of compressible multi--phase flows can be found for example in \cite{USFORCE2,OsherNC,TokarevaToro,TianToro}. 

The full seven--equation Baer--Nunziato model with inter-phase drag and pressure relaxation is given by the following 
non--conservative system of partial differential equations: 
\begin{equation}\left.\begin{array}{l}
\label{ec.BN}


	\frac{\partial}{\partial t}\left(\ar\right)+\nabla\cdot\left(\ar\ub\right)=0,
	
	\\


\frac{\partial}{\partial t}\left(\ar\ub\right)
+\nabla\cdot\left(\phi_1\rho_1\textbf{u}_\textbf{1} \textbf{u}_\textbf{1}\right)+\nabla\phi_1p_1 
= p_I\nabla\phi_1 - \lambda\left(\textbf{u}_\textbf{1}-\textbf{u}_\textbf{2} \right),

\\


\frac{\partial}{\partial t}\left(\phi_1\rho_1E_1\right)
+\nabla\cdot\left(\left(\phi_1\rho_1E_1+\phi_1p_1\right)\textbf{u}_\textbf{1}\right) = 
-p_I\partial_t\phi_1 - \lambda \, \mathbf{u_I} \cdot \left(\mathbf{u_1}-\mathbf{u_2}\right),

\\


\frac{\partial}{\partial t}\left(\phi_2\rho_2\right)+\nabla\cdot\left(\phi_2\rho_2\textbf{u}_\textbf{2}\right)=0,

\\ 


\frac{\partial}{\partial t}\left(\phi_2\rho_2\textbf{u}_\textbf{2}\right)
+\nabla\cdot\left(\phi_2\rho_2\textbf{u}_\textbf{2} \textbf{u}_\textbf{2}\right)+\nabla\phi_2p_2=p_I\nabla\phi_2 
- \lambda \, \left(\textbf{u}_\textbf{2}-\textbf{u}_\textbf{1}\right),

\\


\frac{\partial}{\partial t}\left(\phi_2\rho_2E_2\right)
+\nabla\cdot\left(\left(\phi_2\rho_2E_2+\phi_2p_2\right)\textbf{u}_\textbf{2}\right)=p_I\partial_t\phi_1 
- \lambda \, \mathbf{u_I} \cdot \left(\mathbf{u_2}-\mathbf{u_1}\right),

\\


\frac{\partial}{\partial t}\phi_1+\textbf{u}_\textbf{I}\nabla\phi_1 = \nu (p_1 - p_2).
\end{array}\right\}
\end{equation}

In the entire article, the system is closed by the so-called stiffened gas equation 
of state (EOS) for each phase: 
\begin{equation}
\label{eqn.eos} 
   e_k = \frac{p_k + \gamma_k \pi_k}{\rho_k (\gamma_k -1 )}. 
\end{equation}
Here, $\phi_k$ denotes the volume fraction of phase $k$, $\rho_k$ is the 
density, $\mathbf{u_k}$ is the velocity vector, 
$E_k = e_k + \halb \mathbf{u_k}^2$ and $e_k$ are the phase specific total 
and internal energies, respectively, $\lambda$ is a parameter 
characterizing the friction between both phases and $\nu$ characterizes 
pressure relaxation. For consistency the sum of the volume fractions must
always be unity, i.e. $\phi_1 + \phi_2 = 1$.  
In the literature, one of the phases is often also called the \textit{solid} 
phase and the other one the \textit{gas} phase. Defining arbitrarily the 
first phase as the solid phase in the rest of the paper we will therefore use 
the subscripts $1$ and $s$ as well as $2$ and $g$ as synonyms. 
For the interface velocity and pressure $\mathbf{u_I}$ and $p_I$ we choose 
$\mathbf{u_I} = \mathbf{u_1}$ and $p_I = p_2$ respectively, according to
 \cite{BaerNunziato1986}, although other choices are possible, 
see e.g. the paper by Saurel and Abgrall \cite{SaurelAbgrall}. 
We can cast system \eqref{ec.BN} in the general non-conservative form \eqref{eqn.pde.nc} below 
\begin{equation}
\label{eqn.pde.nc} 
  \frac{\partial \Q}{\partial t} + \nabla \cdot \tens{\F}(\Q) + \tens{\B}(\Q) \cdot \nabla \Q = \S(\Q), \qquad \x \in \Omega \subset \mathds{R}^2, t \in \mathds{R}_0^+ 
\end{equation} 
where $\Q=\left( \ar, \, \ar \mathbf{u}_1, \, \ar E_1, \, \arr, \, \arr \mathbf{u}_2, \, \arr E_2, \, \phi_1 \right) \in \Omega_Q \subset \mathds{R}^\nu$ is the 
state vector,  $\tens{\F} = (\mathbf{f}, \mathbf{g})$ is the flux tensor, i.e. the purely 
conservative part of the PDE system, $\tens{\B}=(\B_1,\B_2)$ contains the purely non-conservative part of the system in block--matrix notation and $\S(\Q)$ is the vector of  algebraic source terms, which, in our case, contains the inter--phase drag and the pressure relaxation terms. \\ 
We furthermore introduce the abbreviation 
$\P = \P(\Q,\nabla \Q) = \B(\Q) \cdot \nabla \Q$. System \eqref{eqn.pde.nc} can also be written equivalently in the following quasi--linear form   
\begin{equation}
\label{eqn.pde.ql} 
  \frac{\partial \Q}{\partial t} + \tens{\A}(\Q) \cdot \nabla \Q = \S(\Q), 
\end{equation} 
with $\tens{\A}(\Q) = (\A_1, \A_2) = \partial \tens{\F}(\Q) / \partial \Q + \tens{\B}(\Q)$. 
The fluxes $\tens{\F}$, the system matrix $\tens{\A}(\Q)$ and the source term vector $\S(\Q)$ can be readily computed from \eqref{ec.BN}. Hyperbolicity of the Baer--Nunziato model and exact Riemann solvers 
have been studied in \cite{AndrianovWarnecke,DeledicquePapalexandris,Schwendeman}. A very important requirement for numerical methods used to solve the Baer--Nunziato model is the sharp resolution of 
material interfaces, i.e. jumps in the volume fractions $\phi_k$. Improved resolution of the material interfaces can be achieved using high order finite volume methods together with little diffusive 
Riemann solvers, such as the Osher scheme and the HLLC method, both of which have already been successfully applied to the Baer--Nunziato model in \cite{OsherNC} and \cite{TokarevaToro}, respectively.  
However, if the interface velocity is small compared to the sound speed of at least one of the two phases, significant numerical smearing of the material contact will occur even with high order
schemes and little diffusive Riemann solvers due to the small CFL number associated with the material wave. It may therefore be necessary to further improve the resolution of material interfaces 
by using a \textit{Lagrangian} method instead of an \textit{Eulerian} one, since in the Lagrangian case the mesh moves with the flow field and therefore allows a sharp tracking of material interfaces 
independent of the CFL number associated with the material contact wave. 

The significantly improved resolution of material interfaces led to intensive research on Lagrangian schemes in the past decades. The construction of Lagrangian schemes can start either directly from 
the conservative quantities such as mass, momentum and total energy \cite{Maire2007,Smith1999}, or from the nonconservative form of the governing equations, as proposed in 
\cite{Benson1992,Caramana1998,Neumann1950}. Furthermore existing Lagrangian schemes can be divided into \textit{staggered mesh} and \textit{cell--centered} approaches, or combinations of them \cite{StagLag}.  
Cell centered Godunov--type finite volume schemes together with the first Roe linearization for Lagrangian gas dynamics have been proposed by Munz in \cite{munz94}. Munz found that in the Lagrangian 
framework, the Roe--averaged velocity for the equations of gas dynamics is simply given by the arithmetic mean, while in the Eulerian frame the Roe average is a more complicated function of the left 
and right densities and velocities. A cell-centered Godunov scheme has been proposed by  Carr\'e et al. \cite{Després2009} for Lagrangian gas dynamics on general multi-dimensional unstructured meshes and  
in \cite{DepresMazeran2003} Despr\'es and Mazeran introduce a new formulation of the multidimensional Euler equations in Lagrangian coordinates as a system of conservation laws associated with constraints. 
Furthermore they propose a way to evolve in a coupled manner both the physical and the geometrical part of the system \cite{Despres2005}, writing the two--dimensional equations of gas dynamics in Lagrangian 
coordinates together with the evolution of the geometry as a weakly hyperbolic system of conservation laws. This allows the authors to design a finite volume scheme for the discretization of Lagrangian gas 
dynamics on moving meshes, based on the symmetrization of the formulation of the physical part. In a recent work Despr\'es et al.  \cite{Depres2012} propose a new method designed for cell-centered Lagrangian 
schemes, which is translation invariant and suitable for curved meshes. General polygonal grids have been considered by Maire et al. \cite{Maire2009,Maire2010,Maire2011}, who develop a general formalism to 
derive first and second order cell-centered Lagrangian schemes in multiple space dimensions. By the use of a node-centered solver \cite{Maire2009}, the authors obtain the time derivatives of the fluxes and  hence a second order method in space and time. Lagrangian schemes for multi--material flows have been successfully introduced in \cite{MaireMM1,MaireMM2,MaireMM3} and Lagrangian schemes with 
additional symmetry preserving properties in cylindrical geometries have been proposed for example in \cite{chengshu3,chengshu4,MaireCyl1,MaireCyl2}. A multi--scale cell centered Godunov--type finite volume scheme for Lagrangian hydrodynamics has been introduced in \cite{Maire2009b}. All the Lagrangian schemes listed before are at most second order accurate in space and time. 

Cheng and Shu were the first to introduce a \textit{third order} accurate essentially non-oscillatory (ENO) reconstruction operator into Godunov--type Lagrangian finite volume 
schemes \cite{chengshu1,chengshu2}. Their cell centered Lagrangian finite volume methods achieve also high order in time, either by using the method of lines (MOL) approach based 
on a third order TVD Runge-Kutta time discretization, or by using a high order one--step Lax-Wendroff-type time stepping. Higher order unstructured Lagrangian finite element methods 
have been recently investigated in \cite{scovazzi1,scovazzi2}. In a very recent paper Dumbser et al. \cite{Dumbser2012} propose a new class of high order accurate Lagrangian--type 
one-–step WENO finite volume schemes for the solution of stiff hyperbolic balance laws. In \cite{BoscheriDumbserLag} this class of high order one--step Lagrangian WENO finite volume schemes
has been extended to unstructured triangular meshes for the \textit{conservative case} and for \textit{single--phase} flows. The present paper is concerned with its extension to  \textit{non--conservative} systems and applications to compressible \textit{multi--phase} flows. 
To our knowledge, the work presented in this article is the first better than second order accurate  unstructured Lagrangian--type finite volume scheme for non--conservative hyperbolic systems with  applications to compressible multi--phase flows. 

To round--off this introduction, we briefly refer also to other Lagrangian--type schemes and alternative Eulerian schemes with improved resolution of material interfaces.  
The following list of references does not pretend to be complete. Among alternative Lagrangian schemes 
one should also mention for example meshless particle schemes, such as the smooth particle hydrodynamics (SPH) method \cite{Monaghan1994,SPHLagrange,SPH3D,Dambreak3D}, which has been 
successfully used to simulate flows in complex domains with large deformations. Since the SPH approach is a fully Lagrangian method the mesh moves with the local fluid velocity, whereas 
in Arbitrary Lagrangian Eulerian (ALE) schemes, see e.g. \cite{Hirt1974,Peery2000,Smith1999,Feistauer1,Feistauer2,Feistauer3,Feistauer4}, the mesh moves with an arbitrary mesh velocity 
that does not necessarily coincide with the real fluid velocity. This adds an extra degree of flexibility and generality to the scheme. The classical SPH method is a truly meshless scheme,
but exhibits also a series of problems, such as the need for artificial stabilization terms and the lack of zeroth order consistency. These problems have been overcome by the so--called 
Particle Finite Element Method (PFEM), which is a Lagrangian (or Arbitrary--Lagrangian--Eulerian) finite element scheme on moving point clouds, \cite{PFEM1,PFEM2,PFEM3,PFEM4,PFEM5,PFEM6}, 
and has been successfully applied to \textit{incompressible} multi--material flows and fluid--solid interaction problems. 
Furthermore, the reader will also find so--called Semi-Lagrangian schemes in literature, which are mainly used for solving transport equations \cite{ALE2000Belgium,ALE1996FV}. Here, the 
numerical solution at the new time is computed from the known solution at the current time by following backward in time the Lagrangian trajectories of the  fluid to the foot--point of the 
trajectory. Semi--Lagrangian schemes are therefore nothing else than a  numerical scheme based on the \textit{method of characteristics}. Since in general the end-point does not coincide with 
a grid point, an interpolation formula is required in order to evaluate the unknown  solution, see e.g.  \cite{Casulli1990,CasulliCheng1992,LentineEtAl2011,HuangQiul2011,QuiShu2011,CIR,BoscheriDumbser}. 
Note that in Semi-Lagrangian algorithms the mesh is \textit{fixed}, like in a classical Eulerian methods. An alternative to Lagrangian methods for the accurate resolution of material interfaces 
has been developed in the Eulerian framework on fixed meshes under the form of ghost--fluid and level--set methods \cite{FedkiwEtAl1,FedkiwEtAl2,FerrariLevelSet,levelset1,levelset2}, as well 
as the volume of fluid method \cite{HirtNichols,Rieber,Loehner}. 

The rest of the paper is organized as follows: in Section \ref{sec.numethod} the high--order 
path--conservative Arbitrary--Lagrangian--Eulerian (ALE) one--step WENO finite volume scheme is
described on unstructured triangular meshes and a numerical convergence study on a smooth unsteady
test problem is carried out, to assess the designed order of accuracy of the scheme in space and 
time. 
In Section \ref{sec.appl} some classical academic benchmark problems with exact or quasi--exact reference 
solution are solved to verify the robustness of the method in the presence of shock waves and significant
mesh distorsion. Finally, in Section \ref{sec.freesurface} our new high order path--conservative ALE method 
is applied to a free--surface flow problem concerning sloshing in a moving tank. For this purpose, a reduced 
version of the Baer--Nunziato model is solved, which has been introduced for the simulation of free surface 
flows on fixed grids in \cite{DIM2D,DIM3D}. 
The paper is rounded--off by some concluding remarks and an outlook to future research in Section \ref{sec.concl}. 

\section{Numerical Method}
\label{sec.numethod} 

The time--dependent computational domain is denoted by $\Omega(t) \subset \mathds{R}^2$ and is discretized at a given time $t^n$ 
by a set of conforming triangles $T^n_i$, the union of which is the \textit{current triangulation} $\mathcal{T}^n_{\Omega}$ 
of the domain $\Omega(t^n)=\Omega^n$ and can be expressed as  
\begin{equation}
\mathcal{T}^n_{\Omega} = \bigcup \limits_{i=1}^{N_E}{T^n_i}.  
\label{trian}
\end{equation}
Here, $N_E$ is the number of elements used to discretize the domain. 

On triangular meshes it is convenient to introduce also a \textit{local} spatial reference coordinate system $\xi-\eta$, which maps the physical element 
$T_i^n$ in the current configuration to the reference element denoted by $T_e$. The spatial mapping reads 
\begin{equation} 
 \x = \x(\xxi,t^n) = \X^n_{1,i} + \left( \X^n_{2,i} - \X^n_{1,i} \right) \xi + \left( \X^n_{3,i} - \X^n_{1,i} \right) \eta
 \label{xietaTransf} 
\end{equation} 
where $\boldsymbol{\xi} = (\xi, \eta)$ and $\mathbf{x}=(x,y)$ are the vectors of the spatial coordinates in the reference system 
and the physical system, respectively, and $\mathbf{X}^n_{k,i} = (X^n_{k,i},Y^n_{k,i})$ is the vector of physical coordinates of 
the $k$-th vertex of triangle $T^n_i$ at time $t^n$. The unit triangle $T_e$ is defined by the nodes 
$\boldsymbol{\xi}_{e,1}=(\xi_{e,1},\eta_{e,1})=(0,0)$, $\boldsymbol{\xi}_{e,2}=(\xi_{e,2},\eta_{e,2})=(1,0)$ 
and $\boldsymbol{\xi}_{e,3}=(\xi_{e,2},\eta_{e,2})=(0,1)$.  

The cell averages, which represent the data that are stored and evolved in time within a finite volume scheme, are defined at time $t^n$ as usual by   
\begin{equation}
  \Q_i^n = \frac{1}{|T_i^n|} \int_{T^n_i} \Q(\x,t^n) d\x.     
 \label{eqn.cellaverage}
\end{equation}  
Here, $|T_i^n|$ denotes the area of triangle $T_i^n$. Higher order in space can be achieved by reconstructing piecewise higher order polynomials $\mathbf{w}_h(\x,t^n)$ 
from the cell averages defined above in Eqn. \eqref{eqn.cellaverage}. For this purpose, we employ a higher order WENO reconstruction procedure following 
\cite{kaeserjcp,DumbserKaeser06b,DumbserKaeser07}. Other WENO schemes on structured and unstructured meshes can be found, e.g., in  \cite{balsarashu,shu_efficient_weno}
and \cite{friedrich,MixedWENO2D,MixedWENO3D,HuShuTri,ZhangShu3D,AboiyarIske}, respectively. Instead of a WENO method, alternative high order nonlinear reconstruction 
operators can be used as well, see e.g. \cite{sonar,abgrall_eno,MOOD}. For the convenience of the reader, the component--wise WENO reconstruction procedure is briefly 
summarized in the next section, for details we refer to \cite{DumbserKaeser06b}. A possible polynomial WENO reconstruction in characteristic variables can be found in 
\cite{DumbserKaeser07}. 

\subsection{WENO Reconstruction} 
\label{sec.weno} 

The reconstructed solution $\mathbf{w}_h(\x,t^n)$ is represented by piecewise polynomials of degree $M$ and is obtained from the given cell averages in an appropriate neighborhood of 
element $T_i^n$, called the reconstruction \textit{stencil} $\mathcal{S}_i^s$. The number of elements inside each reconstruction stencil is denoted by $n_e$. In 2D we use 
in total 7 reconstruction stencils for each element, one central stencil, three forward stencils and three backward stencils, as proposed in \cite{kaeserjcp,DumbserKaeser06b}, 
hence $1 \leq s \leq 7$. 
According to \cite{barthlsq} the total number of stencil elements $n_e$ must be \textit{larger} than the number of degrees of freedom $\mathcal{M} = (M+1)(M+2)/2$ of 
a polynomial of degree $M$. Typically we take $n_e = 2 \mathcal{M}$ in two space dimensions. 

The reconstruction polynomial for each candidate stencil $s$ for triangle $T_i^n$ is written in terms of spatial basis functions $\psi_l(\xxi)$ as 
\begin{equation}
\label{eqn.recpolydef} 
\w^s_h(\x,t^n) = \sum \limits_{l=1}^\mathcal{M} \psi_l(\xxi) \hat \w^{n,s}_{l,i} := \psi_l(\xxi) \hat \w^{n,s}_{l,i},   
\end{equation}
where the mapping $\x=\x(\xxi,t^n)$ is given by Eqn. \eqref{xietaTransf}. In the rest of the paper we will use classical tensor index notation based on the Einstein summation 
convention, which implies summation over two equal indices. The number of the unknown degrees of freedom to be reconstructed for each element is $\mathcal{M}$. 
The basis functions $\psi_l(\xxi)$ are the \textit{orthogonal} basis functions described in \cite{Dubiner,orth-basis,CBS-book}. 

The reconstruction on each stencil $\mathcal{S}_i^s$ is based on integral conservation, i.e.  
\begin{equation}
\label{intConsRec}
\frac{1}{|T^n_j|} \int \limits_{T^n_j} \psi_l(\xxi) \hat \w^{n,s}_{l,i} d\x = \Q^n_j, \qquad \forall T^n_j \in \mathcal{S}_i^s.      
\end{equation}
Since $n_e > \mathcal{M}$ the above system \eqref{intConsRec} is an \textit{over--determined} linear algebraic system that is solved for $\hat \w^{n,s}_{l,i}$ using a 
constrained least--squares technique, see \cite{DumbserKaeser06b}. The linear constraint is that Eqn. \eqref{intConsRec} holds exactly at least for element $T_i^n$.   
The multi--dimensional integrals appearing in the expression above are evaluated using Gaussian quadrature formulae of suitable order, see \cite{stroud} for details. 
Since the triangles are moving in a Lagrangian scheme, the small linear systems \eqref{intConsRec} are solved for each element at the beginning of each time step. 
However, the choice of the stencils $\mathcal{S}_i^s$ remains \textit{fixed} for all times. 

To obtain a higher order essentially non--oscillatory polynomial the scheme must be \textit{nonlinear}, in order to circumvent the Godunov theorem that states that 
linear monotone schemes are at most of order one. The final nonlinear WENO reconstruction polynomial is therefore computed in the following way. 
First, the \textit{smoothness} of each reconstruction polynomial obtained on stencil $\mathcal{S}_i^s$ is measured by a so--called oscillation indicator $\boldsymbol{\sigma}_s$   
\cite{shu_efficient_weno}. According to \cite{DumbserKaeser06b} the smoothness indicator can be easily computed on the reference element using the (universal) oscillation 
indicator matrix $\boldsymbol{\Sigma}_{lm}$ as follows: 
\begin{equation}
\boldsymbol{\sigma}_s = \boldsymbol{\Sigma}_{lm} \hat \w^{n,s}_{l,i} \hat \w^{n,s}_{m,i}, \qquad 
\boldsymbol{\Sigma}_{lm} = \sum \limits_{ \alpha + \beta \leq M}  \, \, \int \limits_{T_e} \frac{\partial^{\alpha+\beta} \psi_l(\xxi)}{\partial \xi^\alpha \partial \eta^\beta} \cdot 
                                                                         \frac{\partial^{\alpha+\beta} \psi_m(\xxi)}{\partial \xi^\alpha \partial \eta^\beta} d\xi d\eta.   
\end{equation} 
The nonlinear weights $\boldsymbol{\omega}_s$ are defined by
\begin{equation}
\tilde{\boldsymbol{\omega}}_s = \frac{\lambda_s}{\left(\sigma_s + \epsilon \right)^r}, \qquad 
\boldsymbol{\omega}_s = \frac{\tilde{\boldsymbol{\omega}}_s}{\sum_q \tilde{\boldsymbol{\omega}}_q},  
\end{equation} 
where we use $\epsilon=10^{-14}$, $r=8$, $\lambda_s=1$ for the one--sided stencils and $\lambda_s=10^5$ for the central stencil, according to \cite{DumbserKaeser06b}. 
The final nonlinear WENO reconstruction polynomial and its coefficients are then given by 
\begin{equation}
\label{eqn.weno} 
 \w_h(\x,t^n) = \sum \limits_{l=1}^{\mathcal{M}} \psi_l(\xxi) \hat \w^{n}_{l,i}, \qquad \textnormal{ with } \qquad  
 \hat \w^{n}_{l,i} = \sum_s \boldsymbol{\omega}_s \hat \w^{n,s}_{l,i}.   
\end{equation}

In order to reduce the computational cost associated with the nonlinear 
WENO reconstruction procedure outlined above, a high order one--step time discretization is used so that reconstruction has to be performed only once for each time step. 

\subsection{Local Space--Time Predictor on Moving Curved Meshes} 
\label{sec.lst} 

Higher order of accuracy in time is achieved by an element--local predictor stage that \textit{evolves} the reconstructed polynomials $\w_h(\x,t^n)$ locally 
in time within each element $T_i(t)$ during the time interval $[t^n;t^{n+1}]$, see \cite{DumbserEnauxToro,Dumbser2008,USFORCE2,HidalgoDumbser,GassnerDumbserMunz}. Such an element--local time--evolution procedure has also been used within the MUSCL scheme of van Leer \cite{leer5} and the original ENO scheme of Harten et al. \cite{eno},  
who called this element--local predictor with initial data $\w_h(\x,t^n)$ the solution of a Cauchy problem \textit{in the small}, since no information from
neighbor elements is used. The coupling with the neighbor elements occurs only later in the final one--step finite volume scheme.  

The local data evolution step leads for each element to piecewise space--time polynomials of degree $M$, denoted by $\q_h(\x,t)$ in the following. 
While the original ENO scheme of Harten et al. uses a higher order Taylor series in time together with the \textit{strong} differential form of the PDE to 
substitute time--derivatives with space derivatives (the so--called Cauchy--Kovalewski or Lax--Wendroff procedure \cite{laxwendroff}), here a 
\textit{weak} formulation of the PDE in space--time is used. The resulting method does not require the computation of higher order derivatives, but just 
pointwise evaluations of the fluxes, source terms and non--conservative products appearing in the PDE. Such a local space--time Galerkin predictor scheme was
introduced for the Eulerian framework in \cite{DumbserEnauxToro,Dumbser2008,USFORCE2,HidalgoDumbser} and is extended here to the 
Lagrangian framework on moving curved space--time elements for PDE with nonconservative products. 

\begin{figure}[!htbp]
\begin{center}
\vspace{4mm}
\begin{tabular}{lr} 
\includegraphics[width=0.4\textwidth]{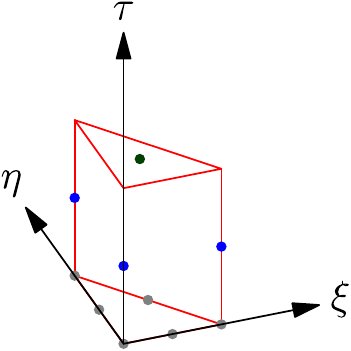}    & 
\includegraphics[width=0.4\textwidth]{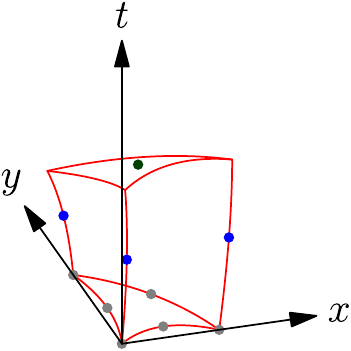}     
\end{tabular}
\vspace{3mm}
\caption{Iso--parametric mapping of the space--time reference element (left) to the physical space--time element (right) used within the local space--time Galerkin predictor.}
\label{fig.stdg.map}
\end{center}
\end{figure}

Let $\mathbf{\tilde{x}}=(x,y,t)$ denote the physical space--time coordinate vector and $\boldsymbol{\tilde{\xi}}=(\xi,\eta,\tau)$ the reference space--time coordinate 
vector, while $\mathbf{x}=(x,y)$ and $\boldsymbol{\xi}=(\xi,\eta)$ are the purely spatial coordinate vectors already 
introduced previously. 
Let furthermore $\theta_l=\theta_l(\boldsymbol{\tilde{\xi}})=\theta_l(\xi,\eta,\tau)$ be a space--time basis function defined by the Lagrange interpolation polynomials passing through the space--time nodes $\boldsymbol{\tilde{\xi}}_m=(\xi_m,\eta_m,\tau_m)$ specified according to \cite{Dumbser2008} and also depicted for the case $M=2$ in 
Figure \ref{fig.stdg.map}. The Lagrange interpolation polynomials  define a \textit{nodal} basis which satisfies the interpolation property 
\begin{equation}
 \theta_l(\boldsymbol{\tilde{\xi}}_m) = \delta_{lm}, 
\end{equation} 
with the usual Kronecker symbol $\delta_{lm}$. 
The element--local space–-time predictor solution $\q_h$, the fluxes $\F_h = (\f_{h}, \g_h)$, the source term $\S_h$ and the non--conservative product 
$\P_h = \B(\q_h) \cdot \nabla \q_h$ are approximated within the space--time element $T_i(t) \times [t^n;t^{n+1}]$ as  
\begin{equation*}
\q_h=\q_h(\boldsymbol{\tilde{\xi}}) = \theta_{l}(\boldsymbol{\tilde{\xi}}) \, \widehat{\q}_{l,i}, \qquad \qquad  
\F_h=\F_h(\boldsymbol{\tilde{\xi}}) = \theta_{l}(\boldsymbol{\tilde{\xi}}) \, \widehat{\F}_{l,i},  
\end{equation*}
\begin{equation} 
\S_h=\S_h(\boldsymbol{\tilde{\xi}}) = \theta_{l}(\boldsymbol{\tilde{\xi}}) \, \widehat{\S}_{l,i}, \qquad \qquad  
\P_h=\P_h(\boldsymbol{\tilde{\xi}}) = \theta_{l}(\boldsymbol{\tilde{\xi}}) \, \widehat{\P}_{l,i}.
\label{thetaSol}
\end{equation}

Since a \textit{nodal} basis is used, the degrees of freedom for $\F_h$, $\S_h$ and $\P_h$ can be simply computed \textit{pointwise} from $\q_h$ as 
\begin{equation}
  \widehat{\F}_{l,i} = \F(\widehat{\q}_{l,i}), \quad 
  \widehat{\S}_{l,i} = \S(\widehat{\q}_{l,i}), \quad 
  \widehat{\P}_{l,i} = \P(\widehat{\q}_{l,i},\nabla \widehat{\q}_{l,i}), \quad 
  \widehat{\nabla \q}_{l,i} = \nabla \theta_{m}(\boldsymbol{\tilde{\xi}}_l) \widehat{\q}_{m,i}.   
\end{equation}  
The degrees of freedom $\widehat{ \nabla \q}_{l,i}$ represent the gradient of $\q_h$ in node $\boldsymbol{\tilde{\xi}}_l$. 

In the present paper an \textit{isoparametric} mapping is used between the physical space--time coordinate vector $\mathbf{\tilde{x}}$ and the reference 
space--time coordinate vector $\boldsymbol{\tilde{\xi}}$, i.e. the mapping is also represented by the nodal basis functions $\theta_l$. Hence,  
\begin{equation}
 \x(\boldsymbol{\tilde{\xi}}) = \theta_l(\boldsymbol{\tilde{\xi}}) \, \widehat{\x}_{l,i}, \qquad 
 t(\boldsymbol{\tilde{\xi}}) = \theta_l(\boldsymbol{\tilde{\xi}}) \, \widehat{t}_l,  
 \label{eqn.isoparametric} 
\end{equation} 
where the degrees of freedom $\widehat{\mathbf{x}}_{l,i} = (\widehat{x}_{l,i},\widehat{y}_{l,i})$ denote the in general unknown vector of physical 
coordinates in space of the moving space--time control volume and the $\widehat{t}_l$ denote the \textit{known} degrees of freedom of the physical 
time at each space--time node $\tilde{\x}_{l,i} = (\widehat{x}_{l,i}, \widehat{y}_{l,i}, \widehat{t}_l)$.  
The general representation of the mapping in time given above in Eqn. \eqref{eqn.isoparametric} simplifies to 
\begin{equation}
t = t_n + \tau \, \Delta t, \qquad  \tau = \frac{t - t^n}{\Delta t}, \qquad \Rightarrow \qquad \widehat{t}_l = t_n + \tau_l \, \Delta t, 
\label{timeTransf}
\end{equation} 
where $t^n$ is the current time and $\Delta t$ is the time step. Hence, $t_\xi = t_\eta = 0$ and $t_\tau = \Delta t$. 
The isoparametric mapping (\ref{eqn.isoparametric}) allow us to transform the physical space--time element to the unit reference space--time element $T_e \times [0,1]$. 
A sketch of this mapping is depicted in Figure \ref{fig.stdg.map} and the Jacobian matrix of the transformation reads 
\begin{equation}
J_{st} = \frac{\partial \mathbf{\tilde{x}}}{\partial \boldsymbol{\tilde{\xi}}} = \left( \begin{array}{ccc} x_{\xi} & x_{\eta} & x_{\tau} \\ y_{\xi} & y_{\eta} & y_{\tau} \\ 0 & 0 & \Delta t \\ \end{array} \right). 
\label{Jac}
\end{equation}
Its inverse is given by 
\begin{equation}
J_{st}^{-1} = \frac{\partial \boldsymbol{\tilde{\xi}}}{\partial \mathbf{\tilde{x}}} = \left( \begin{array}{ccc} \xi_{x} & \xi_{y} & \xi_{t} \\ \eta_{x} & \eta_{y} & \eta_{t} \\ 0 & 0 & \frac{1}{\Delta t} \\ \end{array} \right).
\label{iJac}
\end{equation}

The local reference system and the inverse of the associated Jacobian matrix (\ref{iJac}) are used to rewrite the governing PDE \eqref{eqn.pde.nc} as 
\begin{equation}
\frac{\partial \Q}{\partial \tau} + \Delta t \left[ \frac{\partial \Q}{\partial \xxi} \cdot \frac{\partial \xxi}{\partial t} + \left( \frac{\partial \xxi}{\partial \x} \right)^T \nabla_{\xxi} \cdot \F  + \B(\Q) \cdot \left( \frac{\partial \xxi}{\partial \x} \right)^T \nabla_{\xxi} \Q \right] = \Delta t \mathbf{S}(\Q),
\label{PDECG}
\end{equation}
using $\tau_x = \tau_y = 0$ and $\tau_t = \frac{1}{\Delta t}$ according to (\ref{timeTransf}) and with the notation 
\begin{equation}
 \nabla_{\xxi} = \left( \begin{array}{c} \frac{\partial}{\partial \xi} \\ \frac{\partial}{\partial \eta}  \end{array} \right), \qquad 
 \nabla        = \left( \begin{array}{c} \frac{\partial}{\partial x  } \\ \frac{\partial}{\partial y   }  \end{array} \right) = 
  \left( \begin{array}{cc} \xi_x & \eta _x \\ \xi_y & \eta_y \end{array} \right) 
  \left( \begin{array}{c} \frac{\partial}{\partial \xi} \\ \frac{\partial}{\partial \eta}  \end{array} \right)  = 
  \left( \frac{\partial \xxi}{\partial \x} \right)^T \nabla_{\xxi}. 
\end{equation} 
The term $\frac{\partial \Q}{\partial \xxi} \cdot \frac{\partial \xxi}{\partial t}$ is due to the Lagrangian mesh motion and is zero in the Eulerian case, i.e. for 
fixed meshes. 

We further introduce the abbreviation 
\begin{equation}
 \mathbf{H} = \Delta t \mathbf{S} - \Delta t \left[ \frac{\partial \Q}{\partial \xxi} \cdot \frac{\partial \xxi}{\partial t} + \left( \frac{\partial \xxi}{\partial \x} \right)^T \nabla_{\xxi} \cdot \F  + \B(\Q) \cdot \left( \frac{\partial \xxi}{\partial \x} \right)^T \nabla_{\xxi} \Q \right],
\end{equation} 
and its numerical approximation by 
\begin{equation}
\mathbf{H}_h = \theta_{l}(\boldsymbol{\tilde{\xi}}) \, \widehat{\mathbf{H}}_{l,i},  \qquad 
\label{eqn.stdg.rhs} 
\end{equation} 
as well as the two operators
\begin{equation}
\left[f,g\right]^{\tau} = \int \limits_{T_e} f(\xxi,\tau) g(\xxi,\tau) d\xi d\eta, \quad \left\langle f,g \right\rangle = \int \limits_{0}^{1} \int \limits_{T_e} f(\xxi,\tau)g(\xxi,\tau) d\xi d\eta d\tau  
\label{intOperators}
\end{equation}
that denote the scalar products of two functions $f$ and $g$ over the spatial reference element $T_e$ at time $\tau$ and over the space-time reference element $T_e\times \left[0,1\right]$, respectively. 

Multiplication of \eqref{PDECG} with space--time test functions $\theta_k(\xxi)$, integration over the space--time reference element $T_e \times [0,1]$ and inserting (\ref{thetaSol}) and \eqref{eqn.stdg.rhs} yields 
\begin{equation}
\left\langle \theta_k,\frac{\partial \theta_l}{\partial \tau} \right\rangle \widehat{\q}_{l,i}  
= \left\langle \theta_k,\theta_l \right\rangle \widehat{\mathbf{H}}_{l,i}. \nonumber\\ 
\end{equation}
The term on the left hand side can be integrated by parts in time, which also allows to introduce the initial condition of the local Cauchy problem in a weak form as 
follows: 
\begin{equation}
 \left[ \theta_k(\xxi,1), \theta_l(\xxi,1)\right]^1 \widehat{\q}_{l,i} - \left\langle \frac{\partial \theta_k}{\partial \tau}, \theta_l \right\rangle \widehat{\q}_{l,i}  
= \left[ \theta_k(\xxi, 0), \psi_l(\xxi) \right]^0 \hat \w^n_{l,i} + \left\langle \theta_k,\theta_l \right\rangle \widehat{\mathbf{H}}_{l,i}. 
\label{LagrSTPDECG}
\end{equation}
With the definitions 
\begin{equation}
\K_{1} = \left[ \theta_k(\xxi,1), \theta_l(\xxi,1)\right]^1 - \left\langle \frac{\partial \theta_k}{\partial \tau}, \theta_l \right\rangle, \quad 
\F_0   = \left[ \theta_k(\xxi, 0), \psi_l(\xxi) \right], \quad 
\mathbf{M} = \left\langle \theta_k,\theta_l \right\rangle, 
\end{equation}
the above expression, which is a nonlinear algebraic equation system for the unknown coefficients $\widehat{\q}_{l,i}$ can be written in a more compact matrix form as 
\begin{equation}
 \K_1 \widehat{\q}_{l,i} = \F_0 \hat \w^n_{l,i} + \mathbf{M} \widehat{\mathbf{H}}_{l,i},
\end{equation}
and is conveniently solved using the following iterative scheme: 
\begin{equation}
 \widehat{\q}_{l,i}^{r+1} = \K_1^{-1} \left( \F_0 \hat \w^n_{l,i} + \mathbf{M} \widehat{\mathbf{H}}_{l,i}^r \right), 
 \label{DGfinal} 
\end{equation}
where $r$ denotes the iteration number. For an efficient initial guess based on a second order MUSCL--type scheme, see \cite{HidalgoDumbser}. 

Since the mesh is moving we also have to consider the evolution of the vertex coordinates of the local space--time element, whose motion is described by the ODE system 
\begin{equation}
\frac{d \mathbf{x}}{dt} = \mathbf{V}(\Q,\x,t),
\label{ODEmesh}
\end{equation}
where $\mathbf{V}=\mathbf{V}(\Q,\x,t)$ is the local mesh velocity. In the \textit{Arbitrary Lagrangian-Eulerian} (ALE) framework used in this article, the mesh 
velocity can be chosen independently of the fluid velocity, hence for $\mathbf{V}=0$ the scheme reduces to a pure Eulerian approach, but if $\mathbf{V}$ coincides 
with the local fluid velocity $\mathbf{v}=\mathbf{v}(\Q)$ one obtains a Lagrangian method. The discrete velocity field inside element $T_i(t)$ can be expressed as  
\begin{equation}
\mathbf{V}_h= \theta_{l}(\xxi,\tau) \widehat{\mathbf{V}}_{l,i}, 
\label{Vdof}
\end{equation}
with $ \widehat{\mathbf{V}}_{l,i} = \mathbf{V}(\mathbf{\hat \q}_{l,i}, \hat{\x}_{l,i}, \hat t_l)$.

As in \cite{Dumbser2012} the ODE (\ref{ODEmesh}) can be solved for the unknown coordinate vector $\widehat{\mathbf{x}}_l$ by using again the local space--time DG method: 
\begin{equation}
\K_1 \widehat{\mathbf{x}}_{l,i} = \left[ \theta_k(\xxi,0), \x(\xxi,t^n) \right]^0 + \Delta t \mathbf{M} \, \widehat{\mathbf{V}}_{l,i},
\label{VCG}
\end{equation}
with $\x(\xxi,t^n)$ given according to the spatial mapping \eqref{xietaTransf} based on the known vertex coordinates of triangle $T_i^n$ at time $t^n$. 
This results in the following iteration scheme for the element--local space--time predictor for the nodal coordinates:  
\begin{equation}
 \widehat{\mathbf{x}}^{r+1}_{l,i} = \K_1^{-1} \left( \left[ \theta_k(\xxi,0), \x(\xxi,t^n) \right]^0 + \Delta t \M \widehat{\mathbf{V}}^r_{l,i} \right).
\label{newVertPos}
\end{equation}
Eqn. (\ref{newVertPos}) is iterated \textit{together} with Eqn. (\ref{DGfinal}). The iteration stops when the residuals of (\ref{DGfinal}) and \eqref{newVertPos} 
are less than a prescribed tolerance, which we set to $10^{-12}$ for all examples shown below.  

Once we have carried out the above procedure for all the elements of the computational domain, we end up with an \textit{element--local \textbf{predictor}} for the numerical solution $\q_h$, as well as for the mesh velocity $\mathbf{V}_h$. 

Next, we have to update the mesh \textit{globally}. Let us denote with $\mathcal{V}_k$ the neighborhood of vertex number $k$, i.e. all those elements that have in common the node number $k$. The number of elements in 
the neighborhood  $\mathcal{V}_k$ is denoted with $N_k$. Since the velocity of each vertex is defined by the local predictor within each element, one has to deal with several, in general different, velocities for the 
same node, since all elements belonging to $\mathcal{V}_k$ will in general give a different velocity contribution, according to their element--local predictor. Since we do not admit the geometry to be discontinuous, 
we decide to fix a \textit{unique} node velocity $\overline{\mathbf{V}}_k^n$ to move the node. The final velocity is chosen to be the \textit{average velocity} considering all the contributions 
$\overline{\mathbf{V}}_{k,j}^n$ of the vertex neighborhood as  
\begin{equation}
\overline{\mathbf{V}}_k^n = \frac{1}{N_k}\sum \limits_{T_j^n \in \mathcal{V}_k}{\overline{\mathbf{V}}_{k,j}^n}, \qquad \textnormal{ with } \qquad 
\overline{\mathbf{V}}_{k,j}^n = \left( \int \limits_{0}^{1} \theta_l(\xxi_{e,m(k)}, \tau) d \tau \right) \widehat{\mathbf{V}}_{l,j}. 
\label{NodesVel}
\end{equation} 
The $\xxi_{e,m(k)}$ are the vertex coordinates of the reference triangle $T_e$ corresponding to vertex number $k$, hence $m=m(k)$ with $1 \leq m \leq 3$ is a 
mapping from the global node number $k$ to the element--local vertex number. Since each node now has its own unique velocity, the vertex coordinates can be 
moved according to 
\begin{equation} 
	\mathbf{X}^{n+1}_{k}	= \mathbf{X}^{n}_{k}	+ \Delta t \, \overline{\mathbf{V}}_k^n,  
	\label{eqn.vertex.update}
\end{equation} 
and we can update all the other geometric quantities needed for the computation, e.g. normal vectors, volumes, side lengths, barycenter positions, \textit{etc.}

\subsection{Path--Conservative One--Step Finite--Volume Scheme}
\label{sec.pc} 

\begin{figure}[!htbp]
	\begin{center} 
	\includegraphics[width=0.85\textwidth]{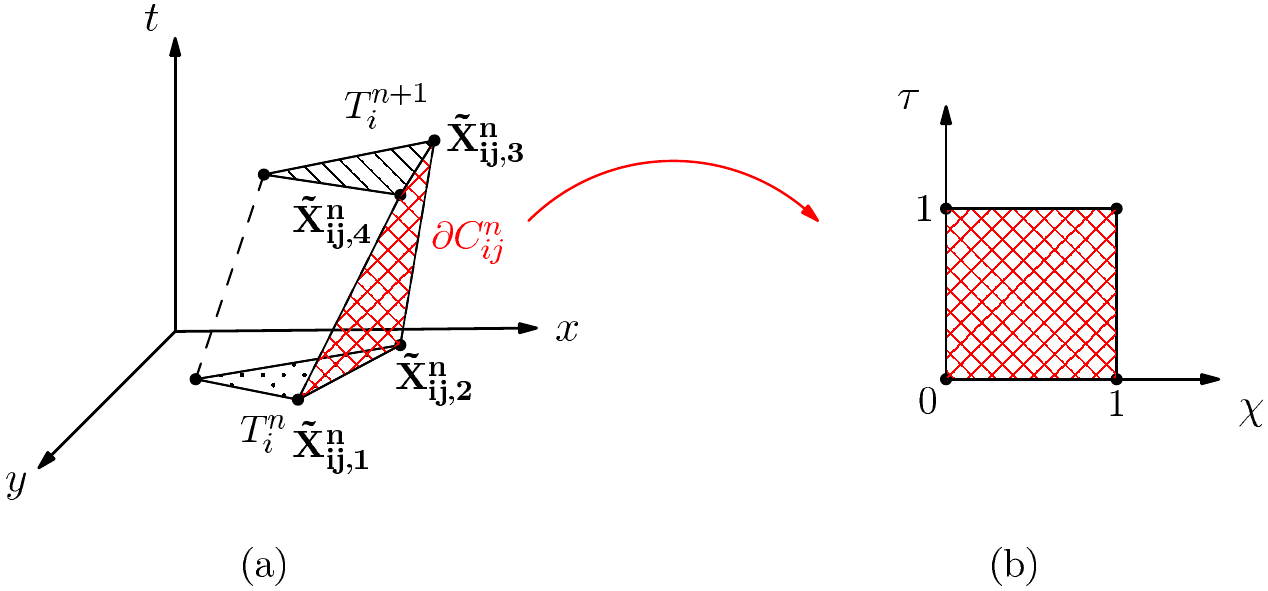}
	\caption{Physical space--time element (a) and parametrization of the lateral 
	         space--time subsurface $\partial C_{ij}^n$  (b).}
	\label{fig:STelem}
	\end{center} 
\end{figure}

The governing PDE system (\ref{eqn.pde.nc}) can be rewritten more compactly using the following space--time divergence form 
\begin{equation}
\tilde \nabla \cdot \tilde{\F} + \tilde \B(\Q) \cdot \tilde \nabla \Q = \mathbf{S}(\Q), 
\label{eqn.st.pde}
\end{equation} 
with the space--time nabla operator 
\begin{equation}
\tilde \nabla  = \left( \frac{\partial}{\partial x}, \, \frac{\partial}{\partial y}, \, \frac{\partial}{\partial t} \right)^T 
\label{eqn.st.nabla} 
\end{equation}
and the space--time flux tensor and system matrices 
\begin{equation}
\tilde{\F}  = \left( \mathbf{f}, \, \mathbf{g}, \, \Q \right), \qquad
\tilde{\B}  = ( \B_1, \B_2, 0), \qquad \tilde{\A} = \frac{\partial \tilde{\F}}{\partial \Q} + \tilde{\B}.    
\label{eqn.st.mat} 
\end{equation}
With \eqref{eqn.st.nabla} and \eqref{eqn.st.mat} the quasi--linear form of the PDE \eqref{eqn.st.pde} reads 
\begin{equation}
 \tilde{\A}(\Q) \cdot \tilde \nabla \Q = \mathbf{S}(\Q). 
\end{equation}
Integration over a space--time control volume $\mathcal{C}^n_i = T_i(t) \times \left[t^{n}; t^{n+1}\right]$ yields 
\begin{equation}
 \int \limits_{\mathcal{C}^n_i} \tilde \nabla \cdot \tilde{\F} \, d\mathbf{x} dt + 
\int\limits_{\mathcal{C}^n_i} \tilde{\B}(\Q) \cdot \tilde \nabla \Q \, d\mathbf{x} dt = 
\int\limits_{\mathcal{C}^n_i} \S(\Q) \, d\mathbf{x} dt.   
\label{STPDE}
\end{equation} 
Application of the theorem of Gauss allows us to write the first space--time volume integral on the left as a flux integral over the space--time 
surface $\partial \mathcal{C}^n_i$. 
Furthermore, the non--conservative product is integrated by using a \textit{path--conservative} approach 
\cite{Toumi1992,Pares2004,Pares2006,Castro2006,Munoz2007,Castro2008,Rhebergen2008,ADERNC,USFORCE2,OsherNC}, which follows the theory of Dal Maso--Le Floch 
and Murat \cite{DLMtheory} and defines the non--conservative term as a Borel measure. For the known limitations and deficiencies of path--conservative schemes 
see \cite{NCproblems,abgrallkarni}. 
Thus 
\begin{equation}
\int \limits_{\partial \mathcal{C}^{n}_i} \left( \tilde{\F} + \tilde{\D} \right) \cdot \ \mathbf{\tilde n} \, dS + \! \! 
\int \limits_{\mathcal{C}^n_i \backslash \partial \mathcal{C}^n_i} \! \! \! \tilde{\B}(\Q) \cdot \tilde \nabla \Q \, d\mathbf{x} dt = 
\int\limits_{\mathcal{C}^n_i} \S(\Q) \, d\mathbf{x} dt,   
\label{I1}
\end{equation}    
where $\mathbf{\tilde n} = (\tilde n_x,\tilde n_y,\tilde n_t)$ is the outward pointing space--time unit normal vector on the space--time surface $\partial C^n_i$, 
and $\tilde{\D}$ is a term that takes into account potential jumps of $\Q$ on the element boundaries according to the path integral 
\begin{equation}
 \tilde{\D} \cdot  \mathbf{\tilde n}  = \int \limits_0^1 \tilde{\B}\left(\Path(\Q^-,\Q^+,s)\right) \cdot \mathbf{\tilde n} \, \frac{\partial \Path}{\partial s} \, ds. 
 \label{eqn.pathint} 
\end{equation} 
Throughout the entire paper and according to \cite{Pares2006,Castro2006,USFORCE2,OsherNC} we use the following straight--line segment path 
\begin{equation}
 \Path = \Path(\Q^-,\Q^+,s) = \Q^- + s (\Q^+ - \Q^-), 
 \label{eqn.segpath} 
\end{equation} 
for which the jump term above \eqref{eqn.pathint} simplifies to 
\begin{equation}
 \tilde{\D} \cdot  \mathbf{\tilde n}  = \left( \int \limits_0^1 \tilde{\B}\left(\Path(\Q^-,\Q^+,s)\right) \cdot \mathbf{\tilde n} \, ds \right) \left( \Q^+ - \Q^- \right). 
 \label{eqn.pathint.seg} 
\end{equation} 
The space--time surface $\partial C^n_i$ above involves overall five space--time sub--surfaces, as 
depicted in  Figure \ref{fig:STelem}:  
\begin{equation}
\partial C^n_i = \left( \bigcup \limits_{T_j(t) \in \mathcal{N}_i} \partial C^n_{ij} \right) 
\,\, \cup \,\, T_i^{n} \,\, \cup \,\, T_i^{n+1},  
\label{dCi}
\end{equation}    
where $\mathcal{N}_i$ denotes the so--called \textit{Neumann neighborhood} of triangle $T_i(t)$, i.e. 
the set of directly adjacent triangles $T_j(t)$ that share a common edge $\partial T_{ij}(t)$ with 
triangle $T_i(t)$. 
The common space--time edge $\partial C^n_{ij}$ during the time interval $[t^n;t^{n+1}]$ is denoted 
above by $\partial C^n_{ij} = \partial T_{ij}(t) \times [t^n;t^{n+1}]$. 

The upper space--time sub--surface $T_i^{n+1}$ and the lower space--time sub--surface $T_i^{n}$ are 
parametrized by $0 \leq \xi \leq 1 \wedge 0 \leq \eta \leq 1-\xi$ and the mapping \eqref{xietaTransf}. 
They are orthogonal to the time coordinate, hence for these faces the space--time unit normal vectors 
simply read $\mathbf{\tilde n} = (0,0,1)$ for $T_i^{n+1}$ and $\mathbf{\tilde n} = (0,0,-1)$ for $T_i^{n}$, 
respectively. 
The lateral space--time sub--faces $\partial C^n_{ij}$ are defined using a simple bilinear 
parametrization, since the old vertex coordinates $\mathbf{X}_{ik}^n$ are given and the new ones 
$\mathbf{X}_{ik}^{n+1}$ are known from \eqref{eqn.vertex.update}.  
\begin{equation}
\partial C_{ij}^n = \mathbf{\tilde{x}} \left( \chi,\tau \right) = 
 \sum\limits_{k=1}^{4}{\beta_k(\chi,\tau) \, \mathbf{\tilde{X}}_{ij,k}^n },	
 \qquad 0 \leq \chi \leq 1,  \quad	0 \leq \tau \leq 1, 										 
\label{SurfPar}
\end{equation}
where $(\chi,\tau)$ represents a side-aligned local reference system according to Figure \ref{fig:STelem}. The $\mathbf{\tilde{X}}_{ij,k}^n$ are the  
physical space--time coordinate vectors for the four vertices that define the lateral space--time sub--surface  $\partial C_{ij}^n$. If $\mathbf{X}^n_{ij,1}$ and $\mathbf{X}^n_{ij,2}$ denote the two spatial nodes at time
$t^n$ that define the common spatial edge $\partial T_{ij}(t^n)$, then the four vectors  $\mathbf{\tilde{X}}_{ij,k}^n$ are  given by 
\begin{equation*}
\mathbf{\tilde{X}}_{ij,1}^n = \left( \mathbf{X}^n_{ij,1}, t^n \right), \qquad 
\mathbf{\tilde{X}}_{ij,2}^n = \left( \mathbf{X}^n_{ij,2}, t^n \right),
\end{equation*} 
\begin{equation}
\mathbf{\tilde{X}}_{ij,3}^n = \left( \mathbf{X}^{n+1}_{ij,2}, t^{n+1} \right), \qquad 
\mathbf{\tilde{X}}_{ij,4}^n = \left( \mathbf{X}^{n+1}_{ij,1}, t^{n+1} \right).  
\label{eqn.lateralnodes} 
\end{equation} 
The $\beta_k(\chi,\tau)$ are a set of bilinear basis functions, which are defined as 
\begin{equation*}
 \beta_1(\chi,\tau) = (1-\chi)(1-\tau), \qquad 
 \beta_2(\chi,\tau) = \chi(1-\tau), 
\end{equation*}
\begin{equation}  
 \beta_3(\chi,\tau) = \chi\tau, \qquad 
 \beta_4(\chi,\tau) = (1-\chi)\tau.
 \label{BetaBaseFunc}
\end{equation}
From \eqref{eqn.lateralnodes} and \eqref{BetaBaseFunc} it follows that the temporal mapping is again simply
$t = t^n + \tau \, \Delta t$, hence $t_\chi = 0$ and $t_\tau = \Delta t$. 
The determinant of the coordinate transformation and the resulting space--time unit normal vector $\mathbf{\tilde n}_{ij}$ of  the sub--surface $\partial C_{ij}^n$ can be computed as follows: 
\begin{equation}
| \partial C_{ij}^n| = \left| \frac{\partial \mathbf{\tilde{x}}}{\partial \chi} \times \frac{\partial \mathbf{\tilde{x}}}{\partial \tau} \right|, 
\quad 
\mathbf{\tilde n}_{ij} = \left( \frac{\partial \mathbf{\tilde{x}}}{\partial \chi} \times \frac{\partial \mathbf{\tilde{x}}}{\partial \tau}\right) / | \partial C_{ij}^n|.
\label{n_lateral}
\end{equation}
Note that the integration over a closed space--time volume as defined above automatically satisfies the geometric conservation law (GCL), 
since from the Gauss theorem follows 
\begin{equation}
 \int_{\partial \mathcal{C}_i^n} \mathbf{\tilde n} \, dS = 0. 
 \label{eqn.gcl} 
\end{equation} 
In \textit{all} numerical simulations shown later in this paper, it has been confirmed for all time steps and all elements that the 
above relation \eqref{eqn.gcl} was always satisfied on the discrete level up to machine precision. 

The final high order ALE one--step finite volume scheme takes the following form: 
\begin{equation}
|T_i^{n+1}| \, \Q_i^{n+1} = |T_i^n| \, \Q_i^n - \sum \limits_{T_j \in \mathcal{N}_i} \,\, {\int \limits_0^1 \int \limits_0^1 
| \partial C_{ij}^n| \tilde{\G}_{ij} \, d\chi d\tau}
+ \int \limits_{\mathcal{C}_i^n \backslash \partial \mathcal{C}_i^n}  \left( \S_h - \P_h \right) \, d\mathbf{x} dt, 
\label{PDEfinal}
\end{equation}
where the term $\tilde{\G}_{ij} \cdot \mathbf{\tilde n}_{ij}$ contains the Arbitrary--Lagrangian--Eulerian numerical flux function as well as the path--conservative jump
term $\tilde{\D}$, to resolve the discontinuity of the predictor solution $\mathbf{q}_h$ at the space--time sub--face 
$\partial C_{ij}^n$. The surface integrals appearing in \eqref{PDEfinal} are approximated using multidimensional Gaussian quadrature rules, see \cite{stroud} for details. 
At the interface $\partial C_{ij}^n$  let us denote the local space--time predictor solution inside element $T_i(t)$ by $\q_h^-$ and 
the element--local predictor solution of the neighbor element $T_j(t)$ by $\q_h^+$, then a simple Rusanov--type scheme \cite{Dumbser2012} is given by 
\begin{equation}
  \tilde{\G}_{ij} =  
  \frac{1}{2} \left( \tilde{\F}(\q_h^+) + \tilde{\F}(\q_h^-)  \right) \cdot \mathbf{\tilde n}_{ij} +  
  \frac{1}{2} \left( \int \limits_0^1 \tilde{\B}(\Path)\cdot \mathbf{\tilde n} \ ds - |\lambda_{\max}| \mathbf{I} \right) \left( \q_h^+ - \q_h^- \right) 
  \label{eqn.rusanov} 
\end{equation} 
where $|\lambda_{\max}|$ is the maximum absolute value of the eigenvalues of the matrix $\tilde{\A} \cdot \mathbf{\tilde n}$ in space--time 
normal direction, which can be expressed in terms of the classical Eulerian system matrix $\A = \partial \F / \partial \Q + \B$ and the normal mesh velocity $\mathbf{V} \cdot \mathbf{n}$ as 
\begin{equation} 
\tilde{\A}_{\mathbf{\tilde n}}  = \tilde{\A} \cdot \mathbf{\tilde n} = \left( \sqrt{\tilde n_x^2 + \tilde n_y^2} \, \right) \left[ \left( \frac{\partial \mathbf{F}}{\partial \Q} + \B \right) \cdot \mathbf{n} - 
(\mathbf{V} \cdot \mathbf{n}) \,  \mathbf{I} \right],  
\end{equation} 
with 
\begin{equation}
\mathbf{n} = \frac{(\tilde n_x, \tilde n_y)^T}{\sqrt{\tilde n_x^2 + \tilde n_y^2}}.  
\end{equation} 
As stated above, $\mathbf{V} \cdot \mathbf{n}$ is the local normal mesh velocity and $\mathbf{I}$ is the $\nu \times \nu$ identity matrix. It can be easily verified that
\begin{equation}
	\mathbf{V} = \frac{1}{\Delta t} \left( \begin{array}{c} x_\tau \\ y_\tau \end{array} \right), \quad  
  \mathbf{\tilde n}_{ij} = \left( \begin{array}{c} 
  \phantom{-} y_\chi \Delta t \\ 
  -x_\chi \Delta t \\    
   x_\chi y_\tau - y_\chi x_\tau 
   \end{array} \right), \quad \textnormal{ hence } \quad 
  \mathbf{V} \cdot \mathbf{n} = -\frac{\tilde n_t}{\sqrt{\tilde n_x^2 + \tilde n_y^2}}.   
\label{eqn.vnormal}    
\end{equation} 
A more sophisticated Osher--type scheme \cite{osherandsolomon} has been introduced in the Eulerian framework for conservative and for non--conservative hyperbolic 
systems in \cite{OsherUniversal,OsherNC}. It has been extended to the Lagrangian framework in one space dimension in \cite{Dumbser2012} and reads 
in the general multi--dimensional case with conservative and non--conservative terms as follows: 
\begin{equation}
  \tilde{\G}_{ij} =  
  \frac{1}{2} \left( \tilde{\F}(\q_h^+) + \tilde{\F}(\q_h^-)  \right) \cdot \mathbf{\tilde n}_{ij}  +  
  \frac{1}{2} \left( \int \limits_0^1 \left( \tilde{\B}(\Path) \cdot \mathbf{\tilde n} -  \left| \tilde{\A}_{\mathbf{\tilde n}}(\Path) \right| \right)  \, ds \, \right) \left( \q_h^+ - \q_h^- \right),  
  \label{eqn.osher} 
\end{equation}
In \eqref{eqn.osher} above, the usual definition of the matrix absolute value operator applies, i.e. 
\begin{equation}
 |\mathbf{A}| = \mathbf{R} |\boldsymbol{\Lambda}| \mathbf{R}^{-1},  \qquad |\boldsymbol{\Lambda}| = \textnormal{diag}\left( |\lambda_1|, |\lambda_2|, ..., |\lambda_\nu| \right),  
\end{equation} 
with the right eigenvector matrix $\mathbf{R}$ and its inverse $\mathbf{R}^{-1}$. According to \cite{OsherNC,OsherUniversal} the path integral 
appearing in \eqref{eqn.osher} is approximated using \textit{Gaussian quadrature rules} of sufficient accuracy. 

\subsection{Numerical Convergence Studies}
\label{sec.conv}

In this section a numerical convergence study is performed for the compressible Baer--Nunziato 
model \eqref{ec.BN} in two space dimensions. This test problem has been proposed in \cite{USFORCE2} and 
has also been used in \cite{OsherNC}. The test problem is similar to the one described  
in \cite{HuShuTri} and \cite{Balsara2004}.  
The exact solution of this smooth unsteady test problem is obtained in two steps: 
 First, an exact 
\textit{stationary} and rotationally symmetric solution of the governing PDE is sought and then the 
problem is made \textit{unsteady} by superimposing a constant, uniform velocity field $\mathbf{\bar v}$ 
using the principle of Galilean invariance of Newtonian mechanics. The exact solution is then simply given 
by the advection of the nontrivial initial condition with the superimposed constant velocity field 
$\mathbf{\bar v}$. The rotationally symmetric solution is found by writing the governing equations 
\eqref{ec.BN} in polar coordinates ($r-\beta$) and by imposing angular symmetry $\partial / \partial \beta = 0$. 
What remains is an ODE system in the radial coordinate $r$ that can be solved analytically, see \cite{USFORCE2}. 

\begin{figure}[!htbp]
\begin{tabular}{cc} 
		\includegraphics[width=0.45\textwidth]{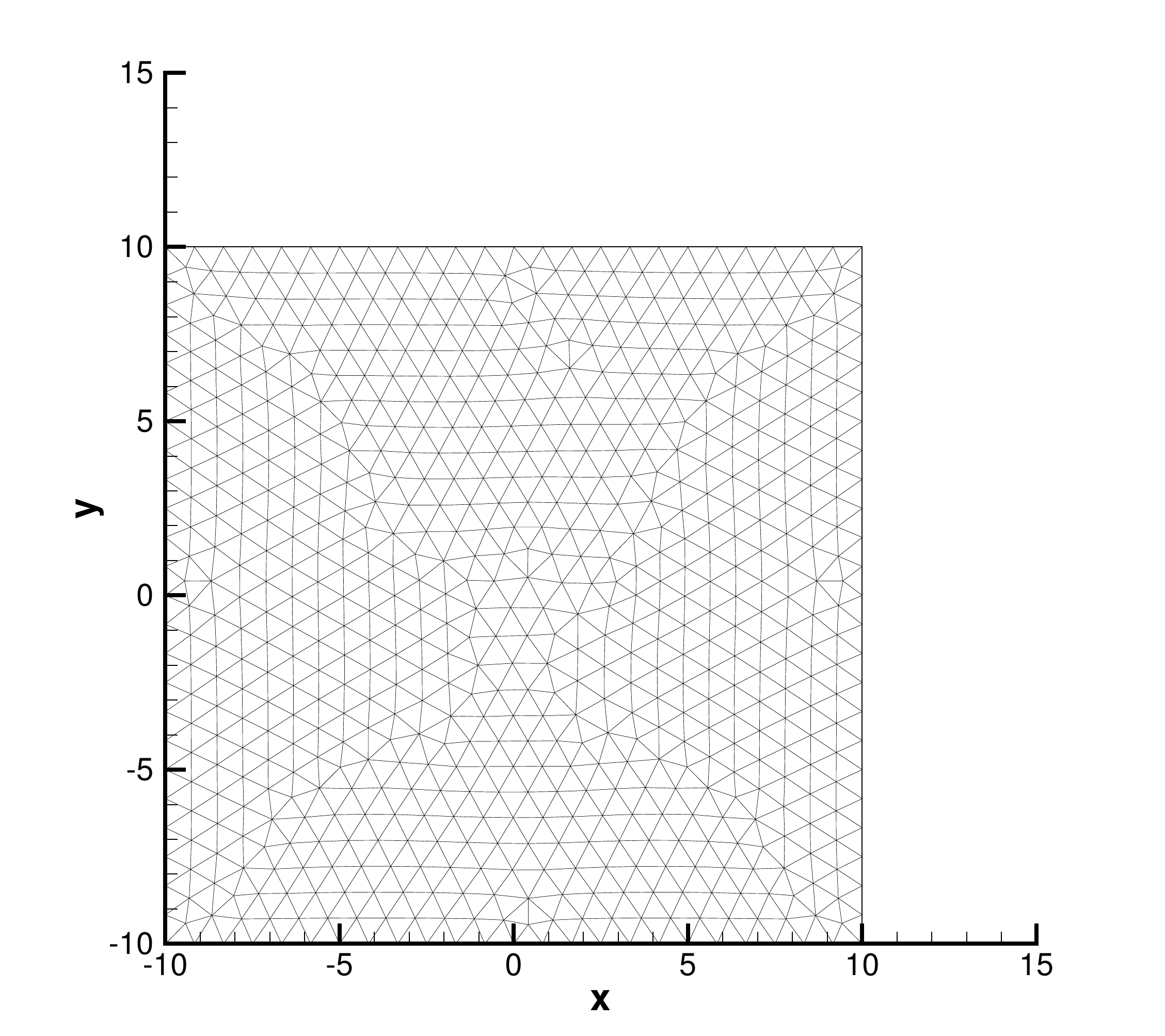} &
		\includegraphics[width=0.45\textwidth]{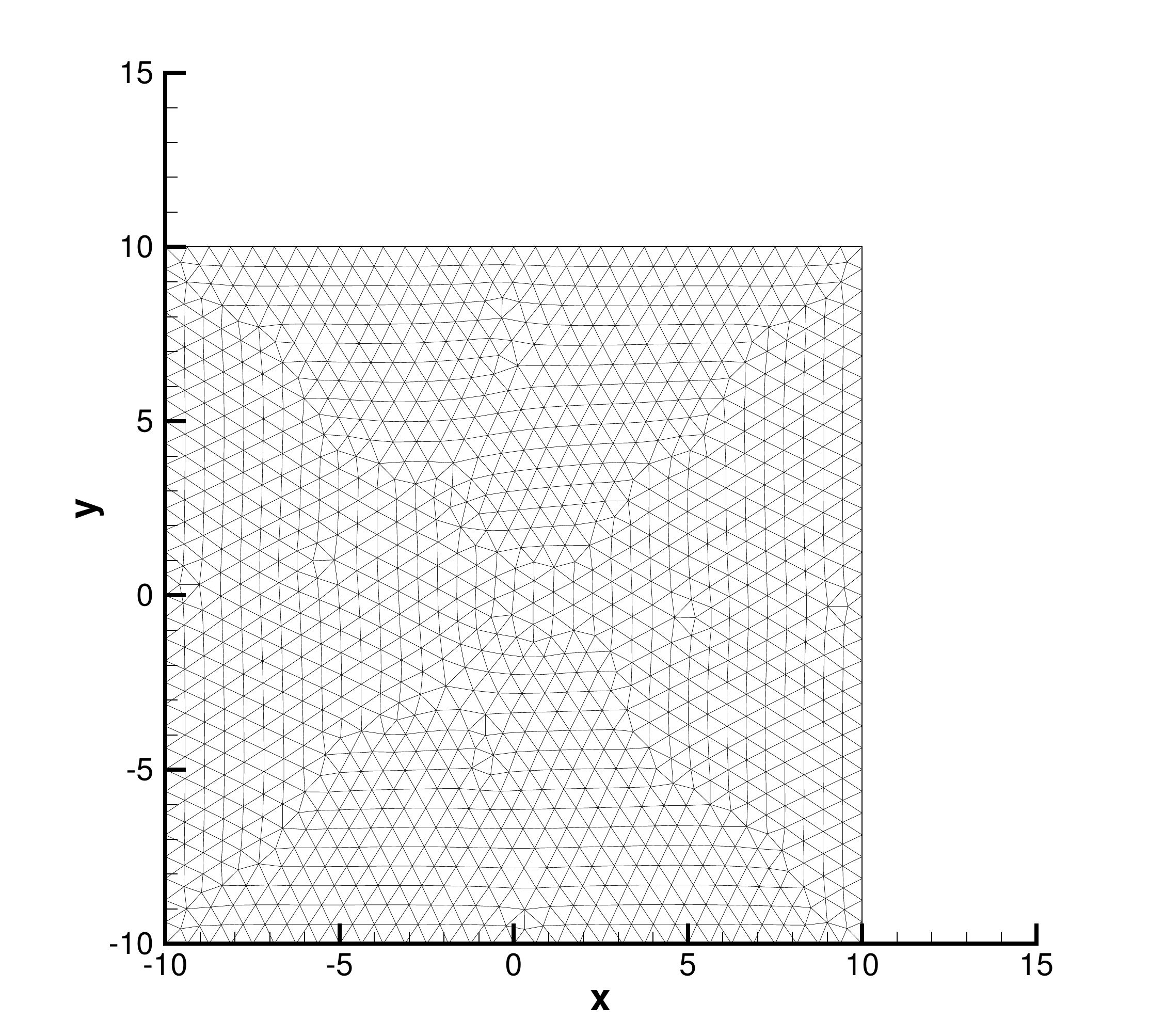} \\
		\includegraphics[width=0.45\textwidth]{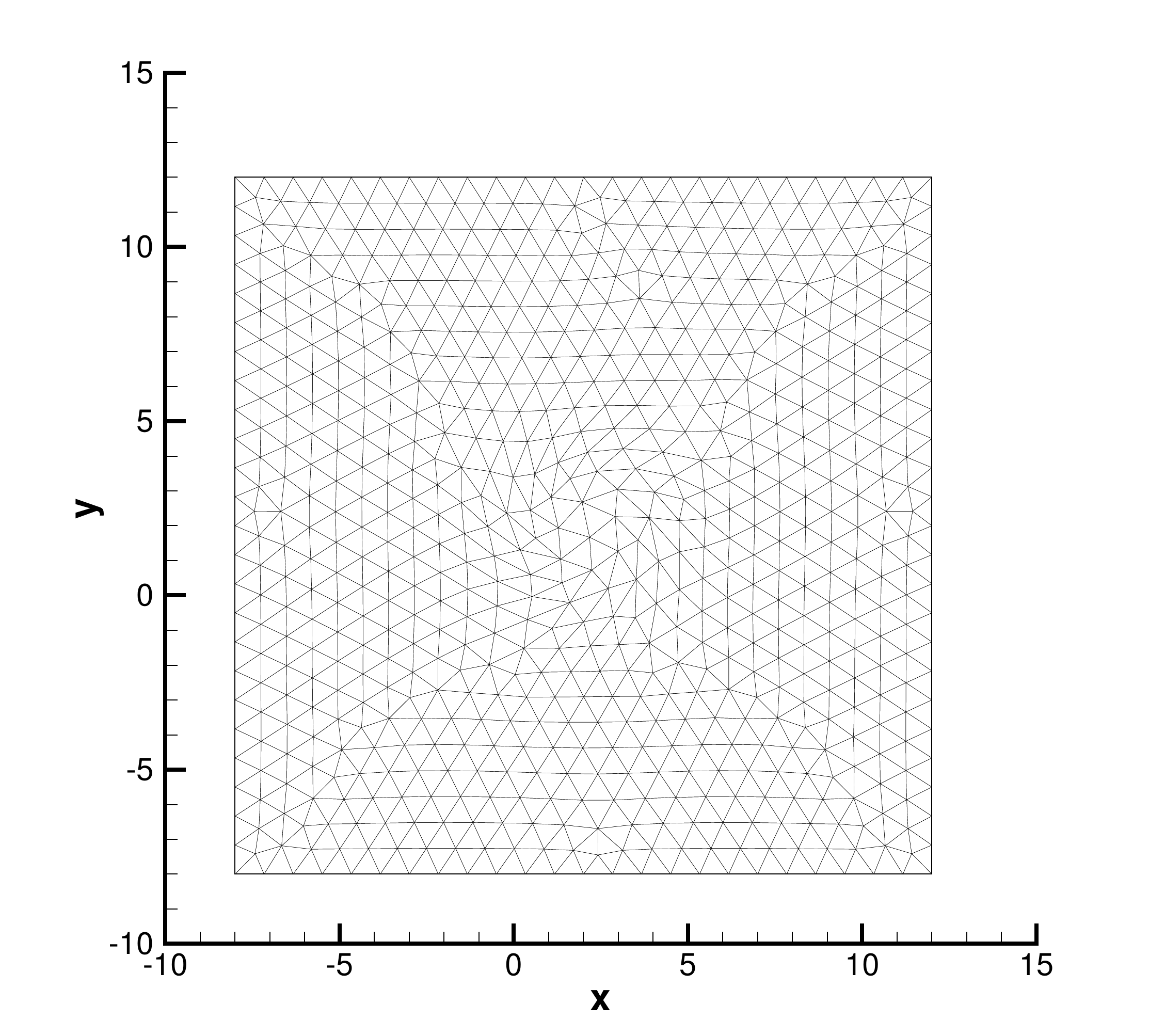} &
		\includegraphics[width=0.45\textwidth]{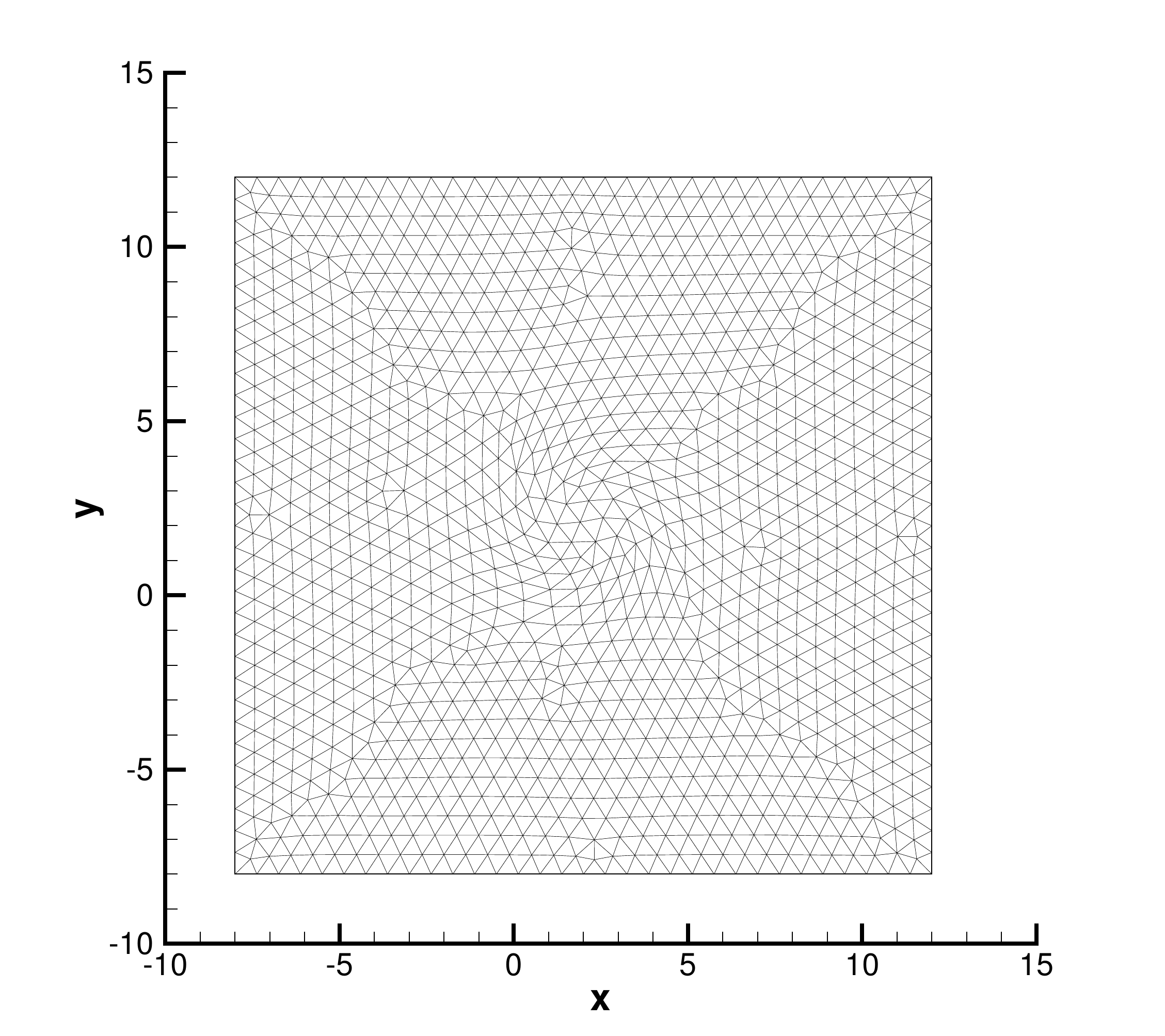} \\
		\includegraphics[width=0.45\textwidth]{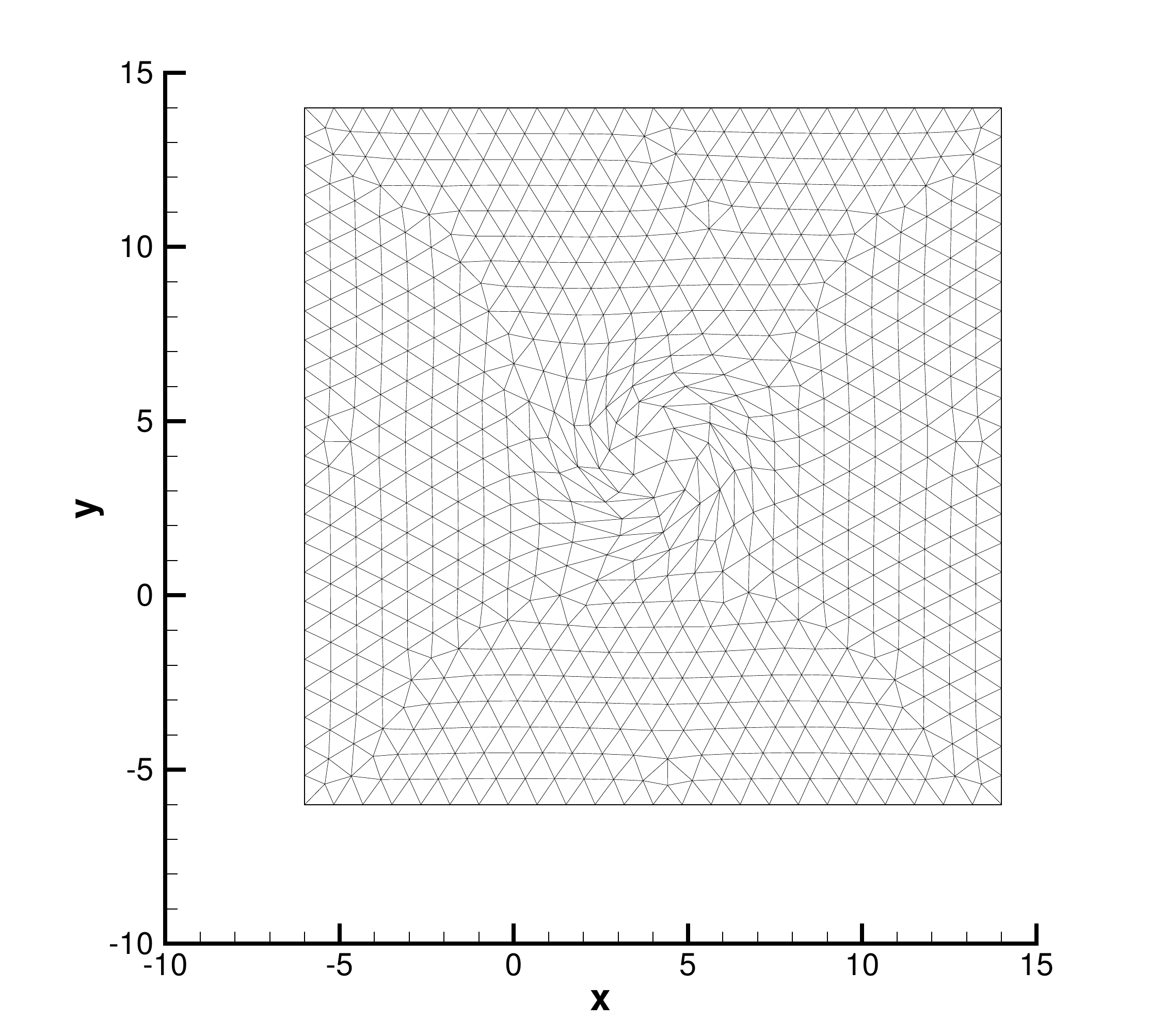} &
		\includegraphics[width=0.45\textwidth]{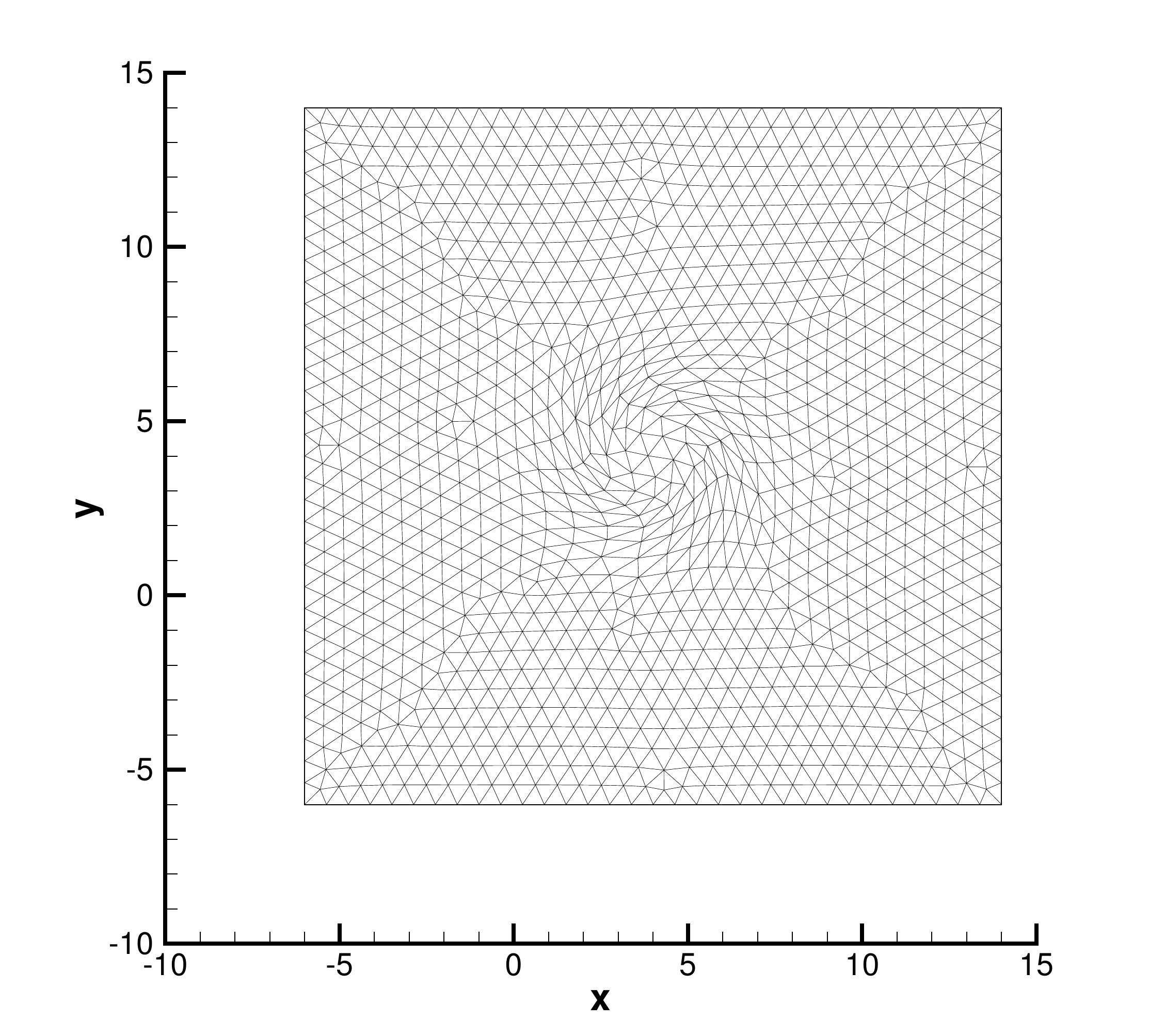} 
\end{tabular} 
		\caption{Moving Lagrangian meshes used for the numerical convergence test at times $t=0$ (top), 
		 $t=1$ (center) and $t=2$ (bottom) with resolution $24 \times 24$ (left) and $32 \times 32$ (right).}
	\label{fig.conv}
\end{figure}

In the following we denote with $u^{\beta}_k$ the angular velocities and with $u^r_k$ the radial velocities. 
Since we are interested in a vortex--type solution, we furthermore suppose that $u^r_k = 0$. From the radial 
momentum equations we then obtain the following ODE system: 
\begin{equation}
\label{eqn.BN.rot.steady}
\begin{array}{ccc} 
\frac{\partial}{\partial r}\left(\phi_1 p_1\right) &=& p_2\frac{\partial}{\partial r}\phi_1+\frac{1}{r}\left(u^{\beta}_1\right)^2 \phi_1\rho_1, \\  
\frac{\partial}{\partial r}\left(\phi_2 p_2\right) &=& p_2\frac{\partial}{\partial r}\phi_2+\frac{1}{r}\left(u^{\beta}_2\right)^2 \phi_2\rho_2.
\end{array} 
\end{equation}
If $\phi_1$, $p_1$ and $p_2$ are know, e.g. by simply \textit{prescribing them}, then \eqref{eqn.BN.rot.steady} is just a 
simple algebraic equation system for the angular velocities $u^{\beta}_k$. As in \cite{USFORCE2} we choose 
\begin{equation}
\label{sol.BN.pres}
p_k=p_{k0}\left(1-\frac{1}{4}e^{\displaystyle{\left(1-r^2/s_k^2\right)}}\right), 
 \qquad (k=1,2)\; ,
\end{equation}
and 
\begin{equation}
\label{sol.BN.alph}
	\phi_1=\frac{1}{3}+\frac{1}{2\sqrt{2\pi}}e^{\displaystyle{-r^2/2}},
\end{equation}
hence the angular velocities of each phase result as
\begin{equation}
\begin{array}{l}
\label{sol.BN.vel}	u_1^\theta=\displaystyle{\frac{1}{2s_1D}}
\sqrt{rD\left[p_{10}\left(4\sqrt{2\pi}F_1+6H_1-12Gs_1^2+3H_1s_1^2\right)+3p_{20}s_1^2\left(4G-H_2\right)\right]}\; ,
\\
\\
u_2^\theta=\displaystyle{ \frac{r\sqrt{2}}{2\rho_2s_2} } 
\sqrt{\rho_2p_{20}F_2}\; ,
\end{array}
\end{equation}
with the auxiliary variables 
\begin{equation*}	
H_k=e^{\displaystyle{-\frac{2r^2+r^2s_k^2-2s_k^2}{2s_k^2}}}, \quad 
F_k=e^{\displaystyle{-\frac{(r-s_k)(r+s_k)}{s_k^2}}}, \quad  (k=1,2),
\end{equation*}
and
\begin{equation*}
G=e^{{-{r^2}/{2}}}, \qquad 
D=\rho_1\left(2\sqrt{2\pi}+3G\right).	
\end{equation*}
To this steady, rotationally symmetric solution of the compressible Baer-Nunziato equations we now add a
constant uniform velocity field $\mathbf{\bar v}=(\bar u, \bar v)$ to make the test problem unsteady, as 
already mentioned above. This can be done since Newtonian mechanics is Galilean invariant. 
With this manufactured analytical solution we can now calculate the convergence rates of the new class of 
high order Lagrangian one--step WENO finite volume schemes presented previously in this section. 
For the computational setup, we use the following parameters: 
\begin{equation*}
	\gamma_1 = 1.4, \quad \gamma_2 = 1.35, \quad \pi_1 = \pi_2 = 0, \quad \bar{u}=\bar{v}=2, \quad \nu=\lambda=0,  
\end{equation*}
\begin{equation}
 \rho_1 = 1, \quad \rho_2 = 2, \quad p_{10} = 1, \quad p_{20} = \frac{3}{2}, \quad s_1=\frac{3}{2}, \quad s_{2}=\frac{7}{5}. 
\label{eqn.BN.param}
\end{equation}
The problem is solved with the Osher--type scheme \eqref{eqn.osher} on a square domain $\Omega=[-10;10] \times [-10;10]$, using 
unstructured triangular meshes with four periodic boundary conditions, see Fig. \ref{fig.conv}, and setting the mesh
velocity to $\mathbf{V}=\mathbf{u}_I=\mathbf{u}_1$. The numerical convergence 
rates are shown for the solid volume fraction $\phi_s$ at time $t=2.0$ in Table \ref{tab.conv}. 
One observes that the schemes reach their designed order of accuracy quite well. To our 
knowledge, this is the first time ever that a better than second order accurate Lagrangian WENO finite volume scheme 
is presented for non--conservative hyperbolic systems on unstructured triangular meshes with applications to the 
Baer--Nunziato model of compressible multi--phase flows. 
\begin{table}  
\caption{Numerical convergence results for the compressible Baer--Nunziato model using the third to sixth order version of the 
Arbitrary--Lagrangian--Eulerian one--step WENO finite volume schemes presented in this article. 
The error norms refer to the variable $\phi_s$ (solid volume fraction) at time $t=2.0$.}  
\begin{center} 
\begin{small}
\renewcommand{\arraystretch}{1.0}
\begin{tabular}{cccccc} 
\hline
  $N_G$ & $\epsilon_{L_2}$ & $\mathcal{O}(L_2)$ & $N_G$ & $\epsilon_{L_2}$ & $\mathcal{O}(L_2)$ \\ 
\hline
  \multicolumn{3}{c}{$\mathcal{O}3$} & \multicolumn{3}{c}{$\mathcal{O}4$} \\
\hline
 24 & 2.6916E-02 & -   &  24 & 1.5993E-02 & -      \\ 
 32 & 1.0906E-02 & 3.1 &  32 & 3.8281E-03 & 5.0    \\ 
 64 & 1.9750E-03 & 2.5 &  64 & 3.0900E-04 & 3.6    \\ 
128 & 2.5442E-04 & 3.0 & 128 & 2.0855E-05 & 3.9    \\ 
\hline 
  \multicolumn{3}{c}{$\mathcal{O}5$} & \multicolumn{3}{c}{$\mathcal{O}6$} \\
\hline
 24 & 1.4493E-02 & -   &  24 & 8.3869E-03 & -      \\ 
 32 & 3.8912E-03 & 4.6 &  32 & 1.9504E-03 & 5.1    \\ 
 64 & 2.5564E-04 & 3.9 &  64 & 6.1843E-05 & 5.0    \\ 
128 & 8.7457E-06 & 4.9 &  96 & 7.4509E-06 & 5.2    \\ 
\hline 
\end{tabular}
\end{small}
\end{center}
\label{tab.conv}
\end{table} 

\section{Test Problems}
\label{sec.appl} 

All the subsequent test problems are solved using the Osher--type method \eqref{eqn.osher} and using 
the interface velocity, i.e. the solid velocity, as mesh velocity, hence $\mathbf{V}=\mathbf{u}_I=\mathbf{u}_1$. 
The CFL number in all test problems is set to $\textnormal{CFL}=0.5$. For all test problems we use a 
third order WENO scheme with reconstruction in characteristic variables \cite{DumbserKaeser07} since 
componentwise reconstruction in conservative variables led to significant spurious oscillations. 
In all cases shown below, friction and pressure relaxation are neglected, hence ($\lambda=\nu=0$). 

\subsection{Riemann Problems} 
\label{sec.rp}

The high order finite volume ALE schemes proposed in this article are subsequently 
validated by applying them to 1D Riemann problems that are solved in a 2D geometry 
on unstructured triangular meshes. The exact solution for these 1D Riemann problems 
can be found in \cite{AndrianovWarnecke,Schwendeman,DeledicquePapalexandris}. 
From the above mentioned articles we have chosen a subset of four Riemann 
problems, whose initial conditions are listed in Table \ref{tab.rpbn.ic}. Some of 
the test cases use the stiffened gas EOS, some of them consider just a mixture of 
two ideal gases. 

The initial two--dimensional computational domain is given by $\Omega(0)=[-0.5;0.5] \times [-0.05; 0.05]$,
which is discretized using an initial characteristic mesh spacing of $h=1/200$, corresponding to an equivalent 
one-dimensional resolution of 200 cells. The initial discontinuity is located at $x=0$ and the final simulation 
times are listed in Table \ref{tab.rpbn.ic}.  In $x$-direction we use transmissive boundaries and in $y$-direction 
periodic boundary conditions are imposed. 

\begin{figure}[!ht]
\begin{center}
\begin{tabular}{cc} 
\includegraphics[width=0.4\textwidth]{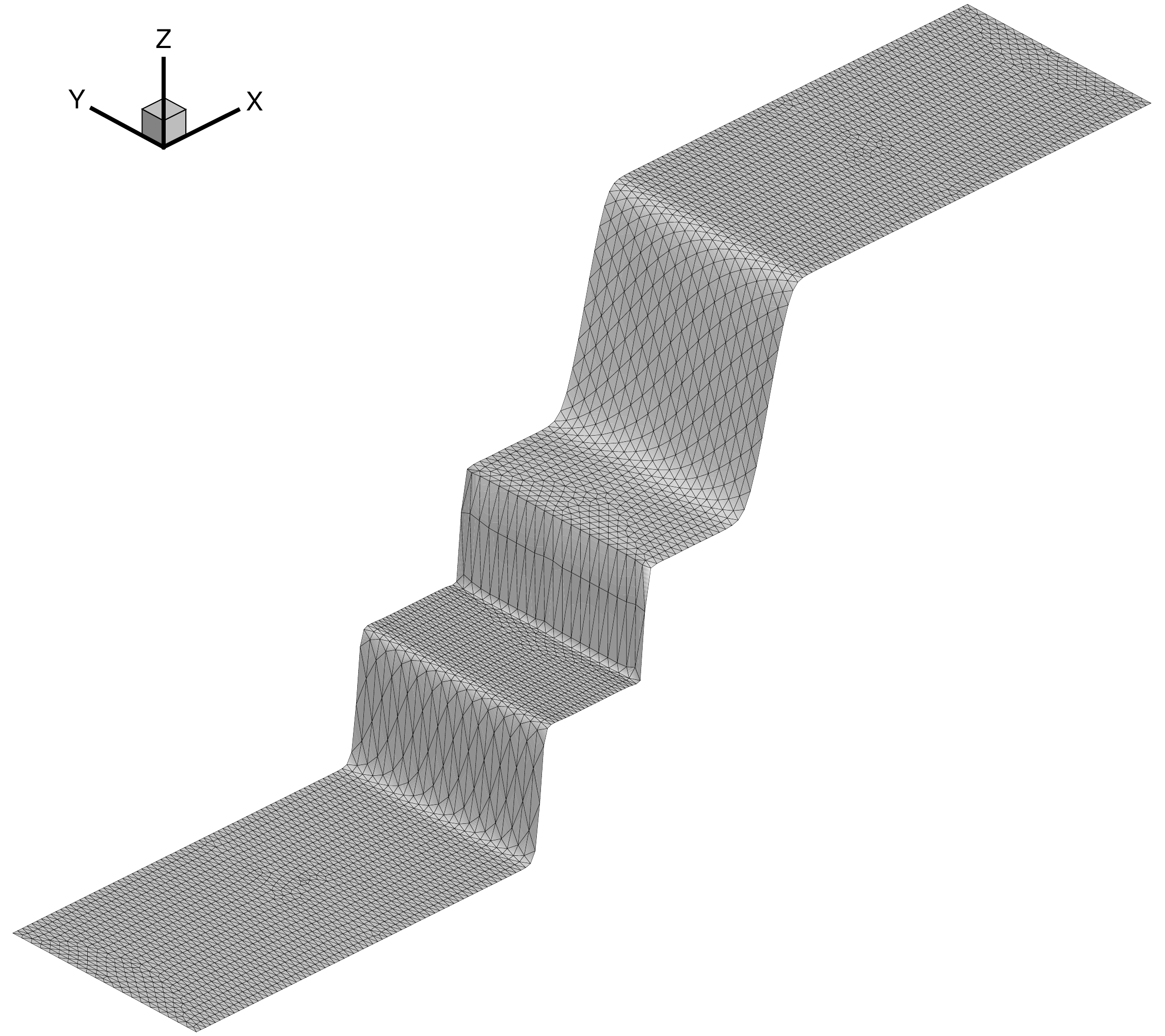}    & 
\includegraphics[width=0.4\textwidth]{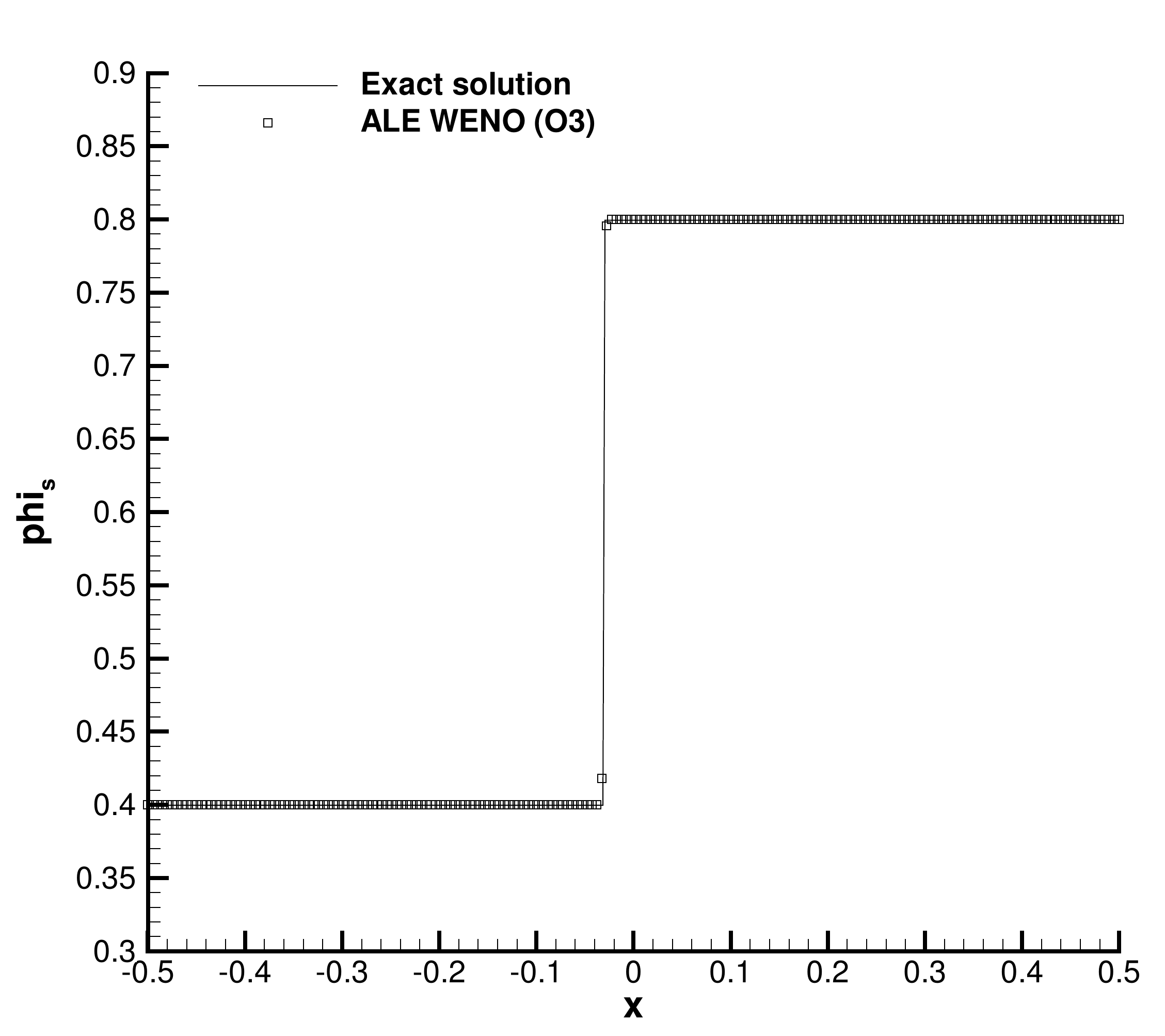}    \\ 
\includegraphics[width=0.4\textwidth]{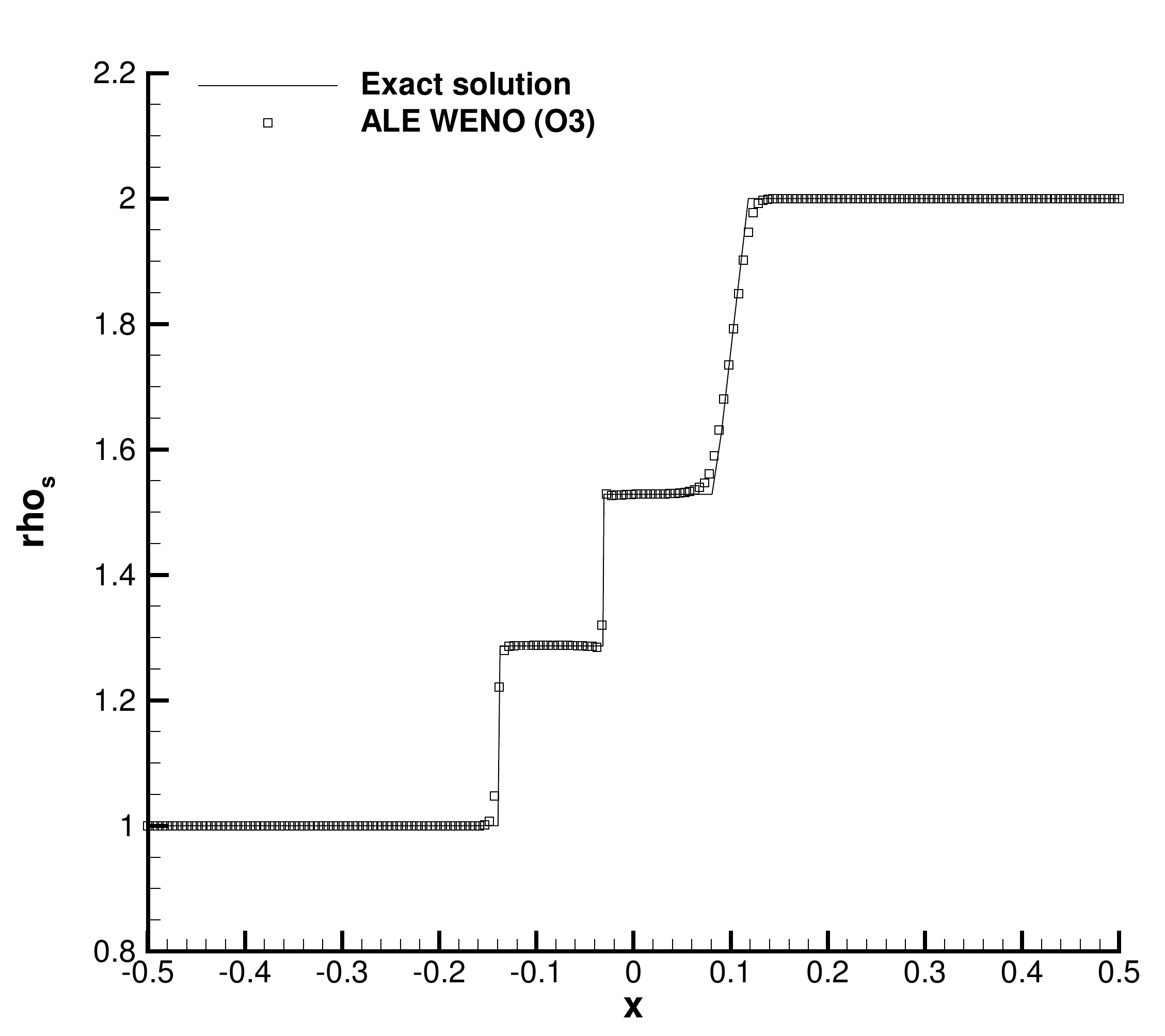}    & 
\includegraphics[width=0.4\textwidth]{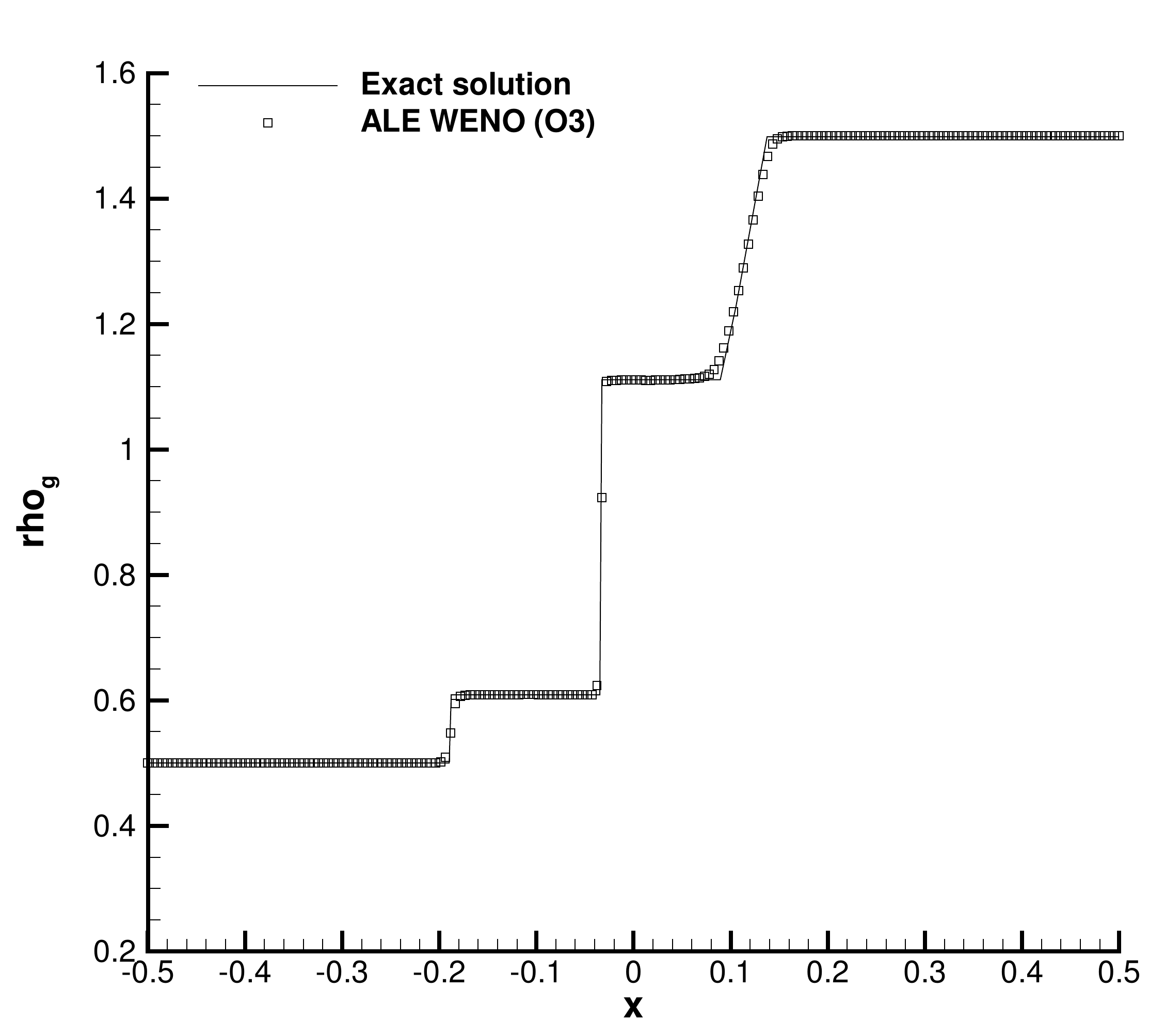}    \\ 
\includegraphics[width=0.4\textwidth]{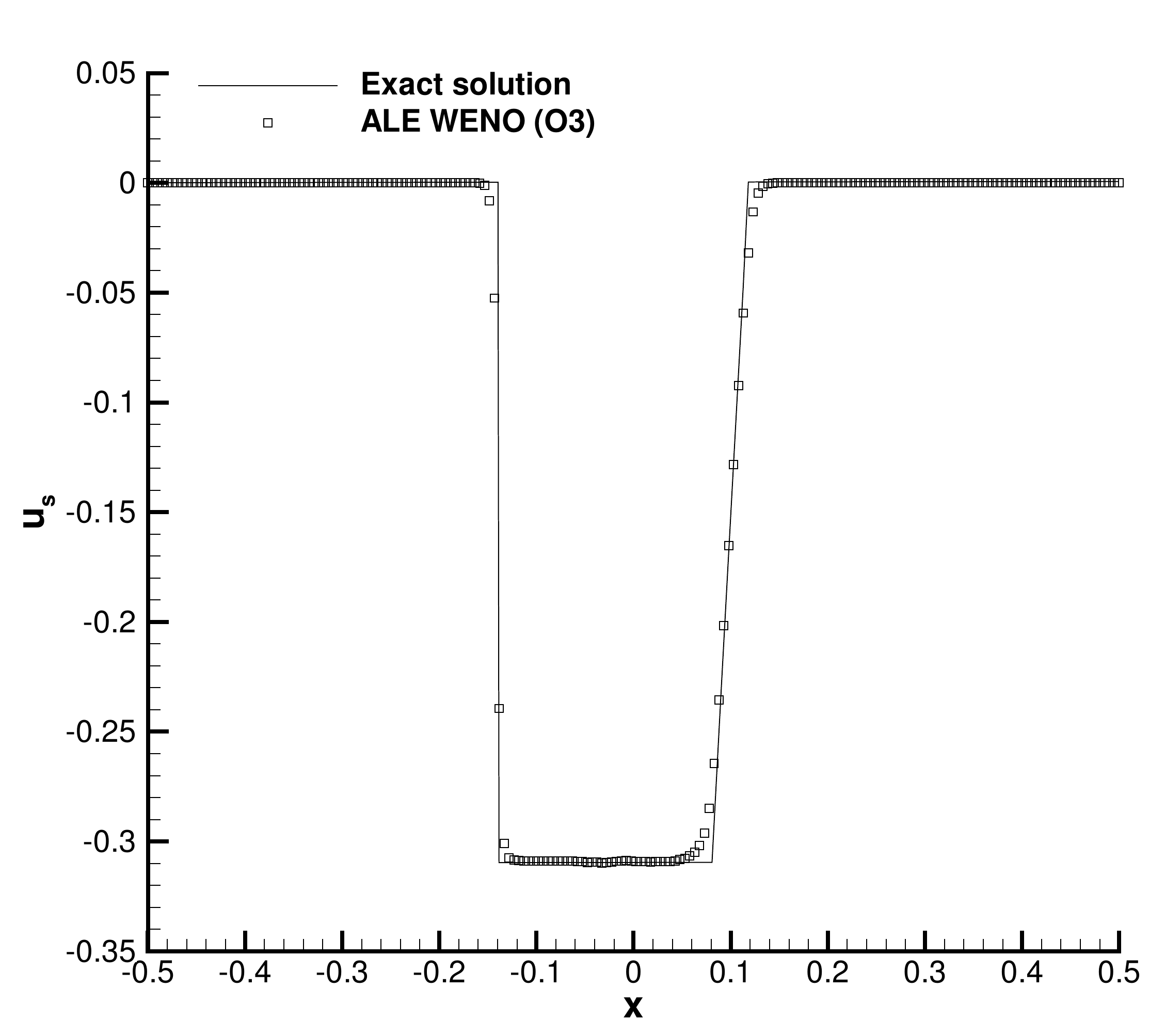}      &  
\includegraphics[width=0.4\textwidth]{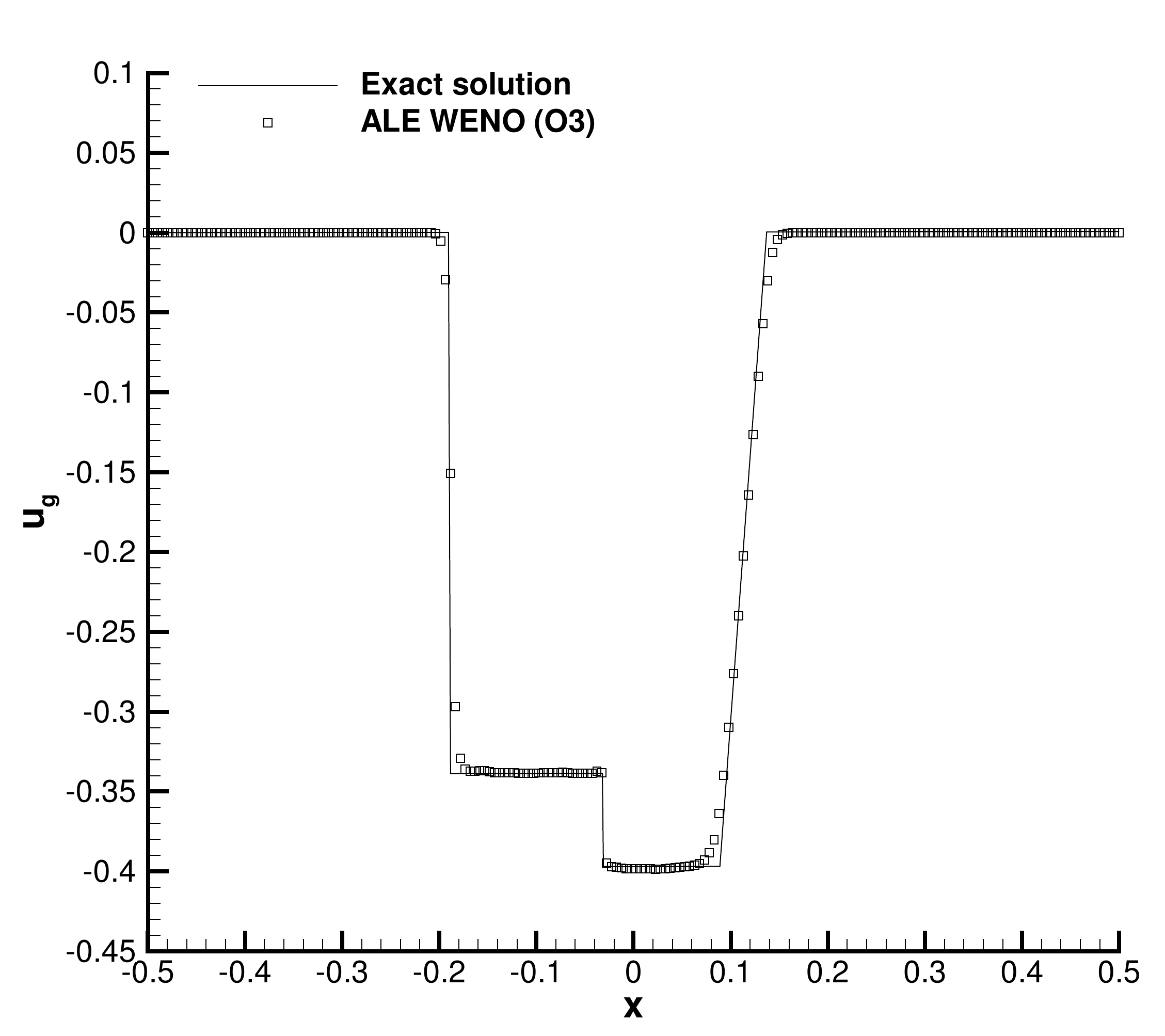}      \\ 
\includegraphics[width=0.4\textwidth]{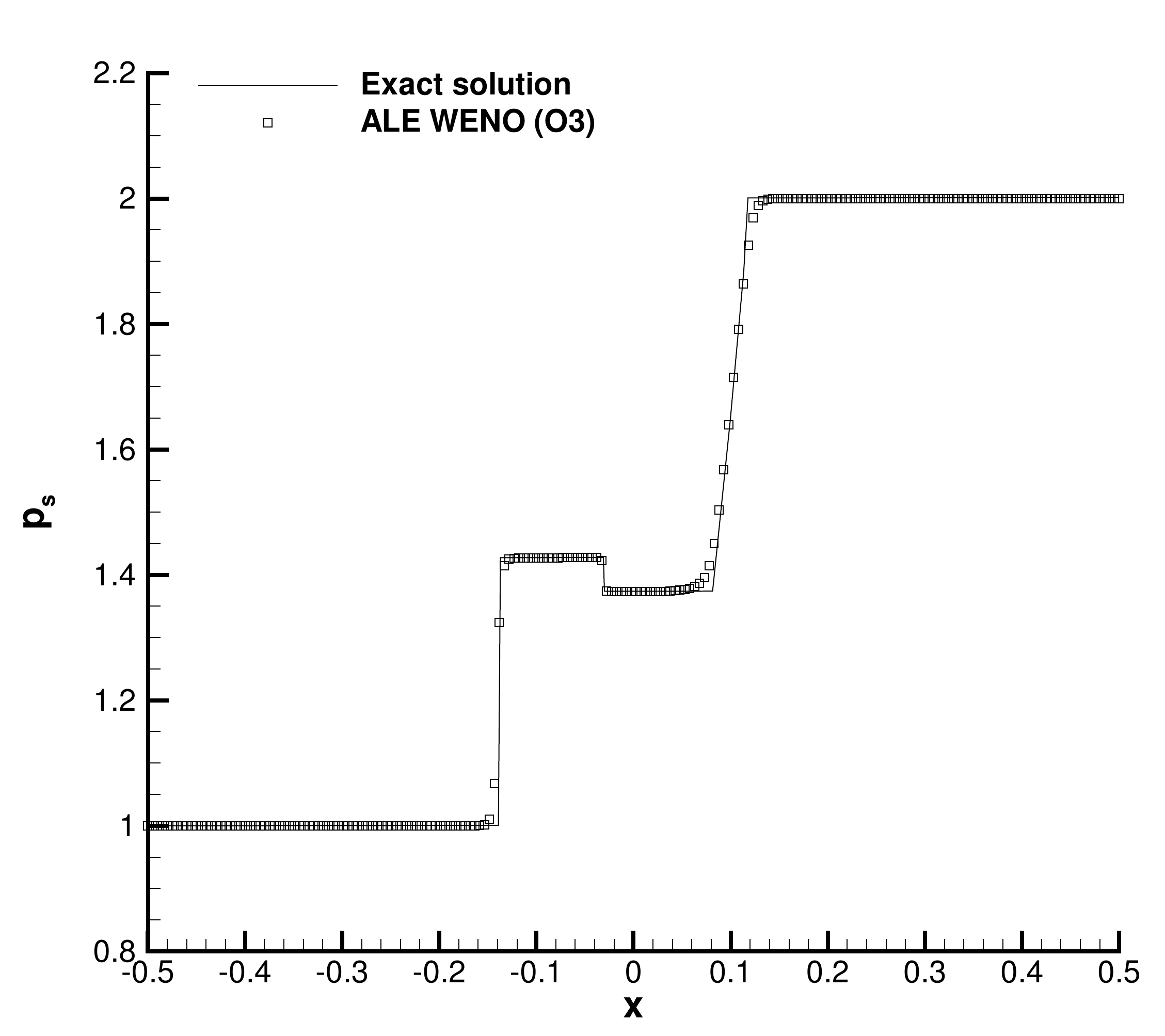}      & 
\includegraphics[width=0.4\textwidth]{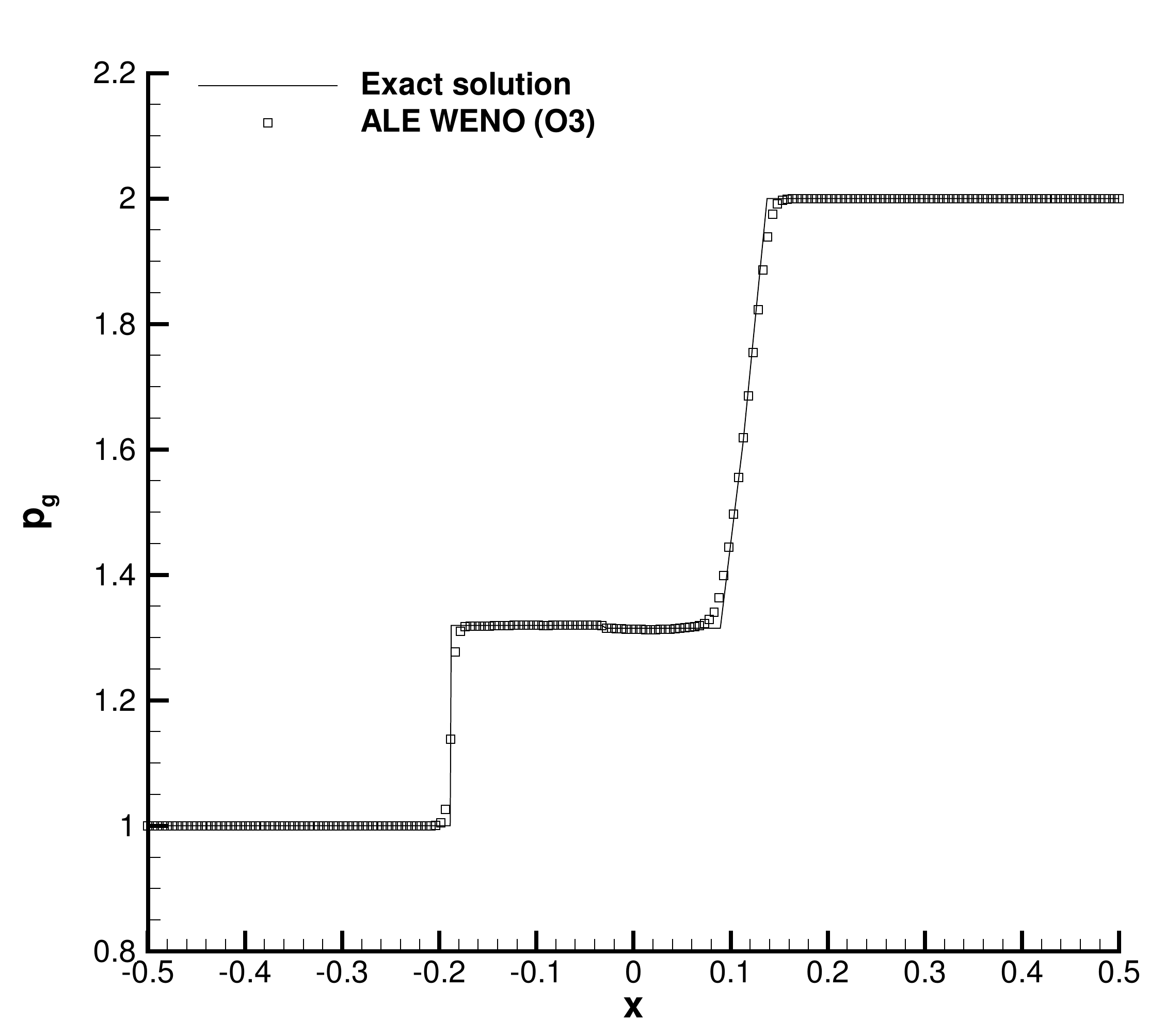}   
\end{tabular}
\caption{Results for Riemann problem RP1 of the seven-equation Baer-Nunziato model at time $t=0.1$.}
\label{fig.bn.rp1}
\end{center}
\end{figure}

\begin{figure}[!ht]
\begin{center}
\begin{tabular}{cc} 
\includegraphics[width=0.4\textwidth]{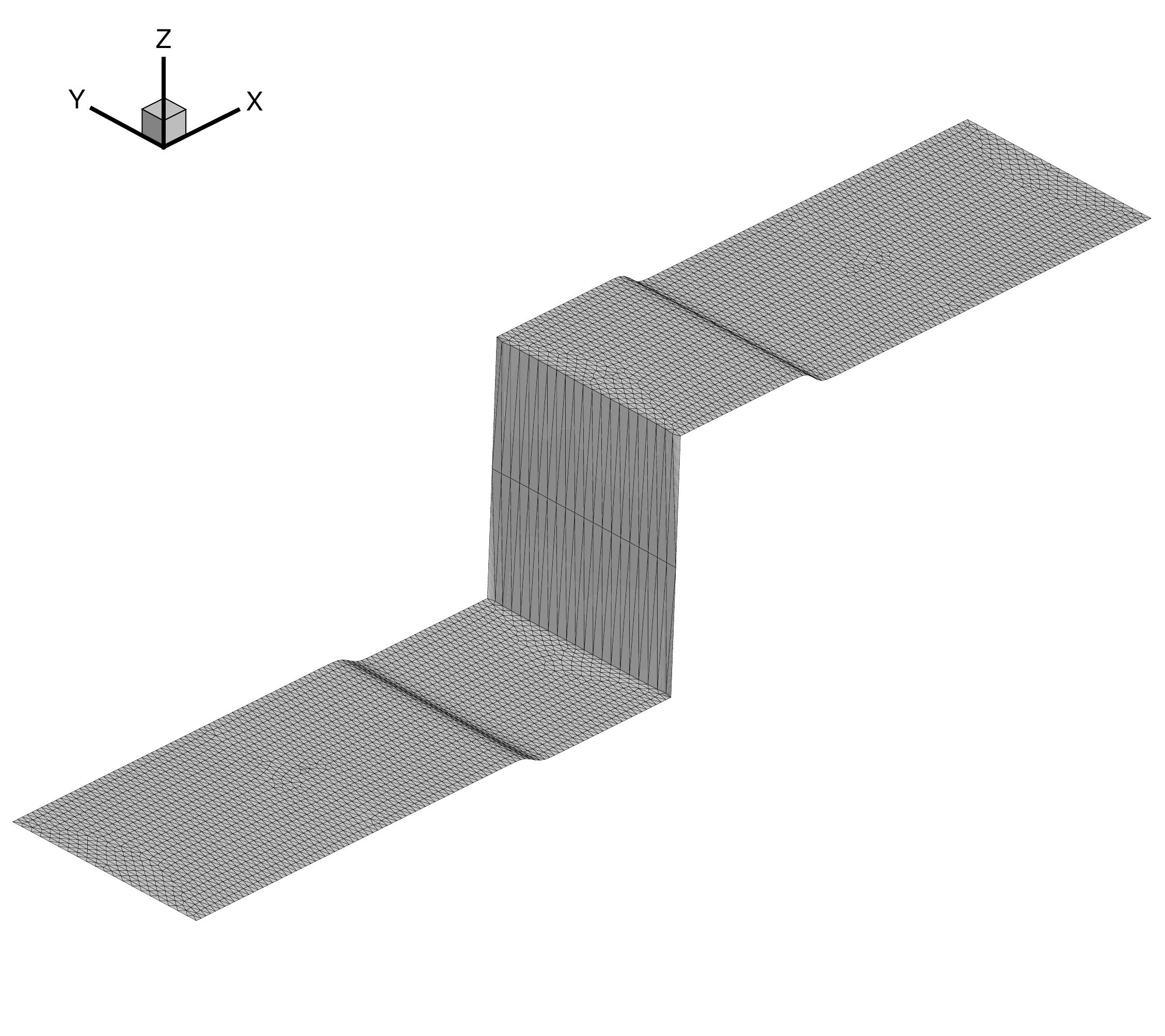}    & 
\includegraphics[width=0.4\textwidth]{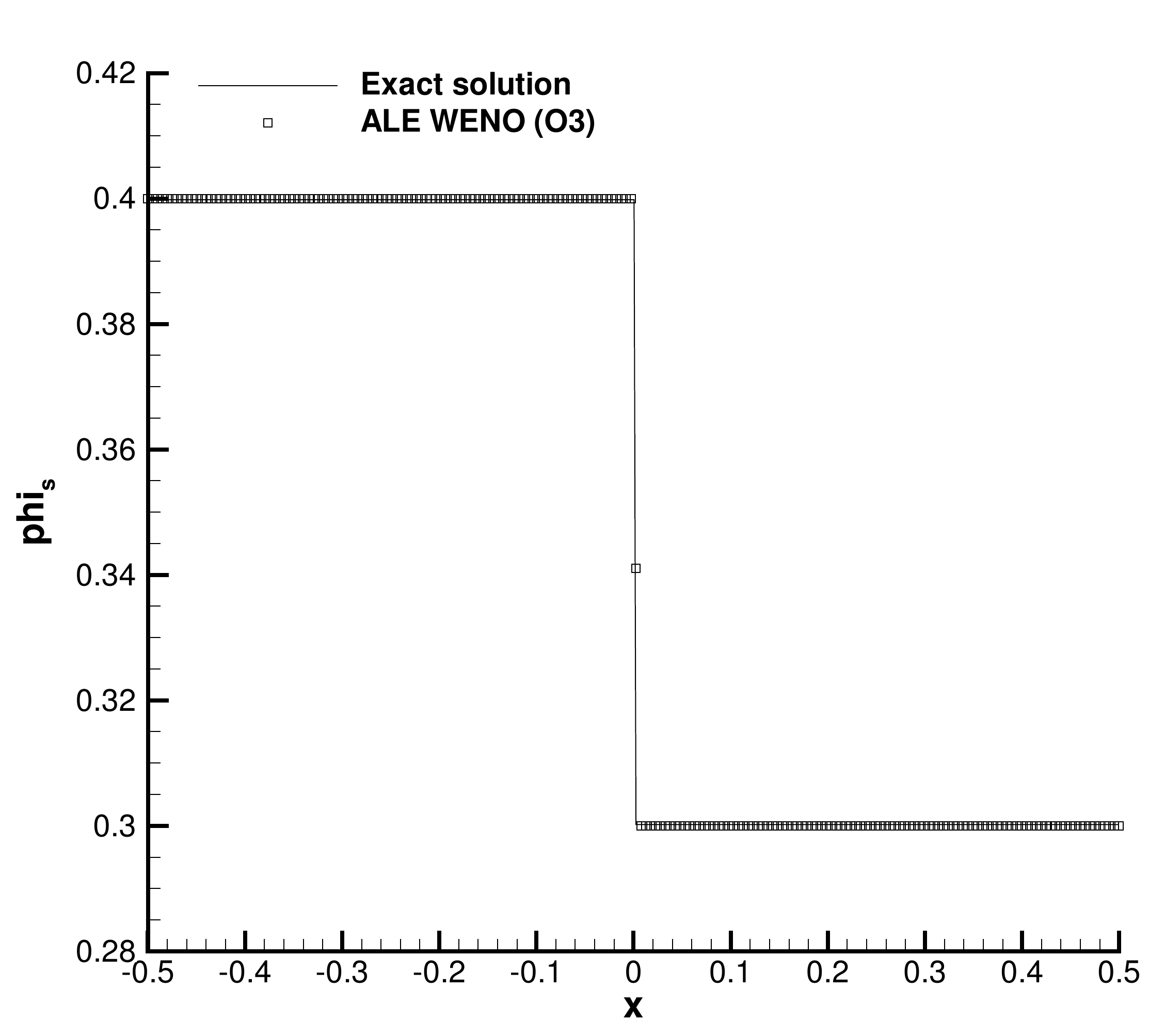}    \\ 
\includegraphics[width=0.4\textwidth]{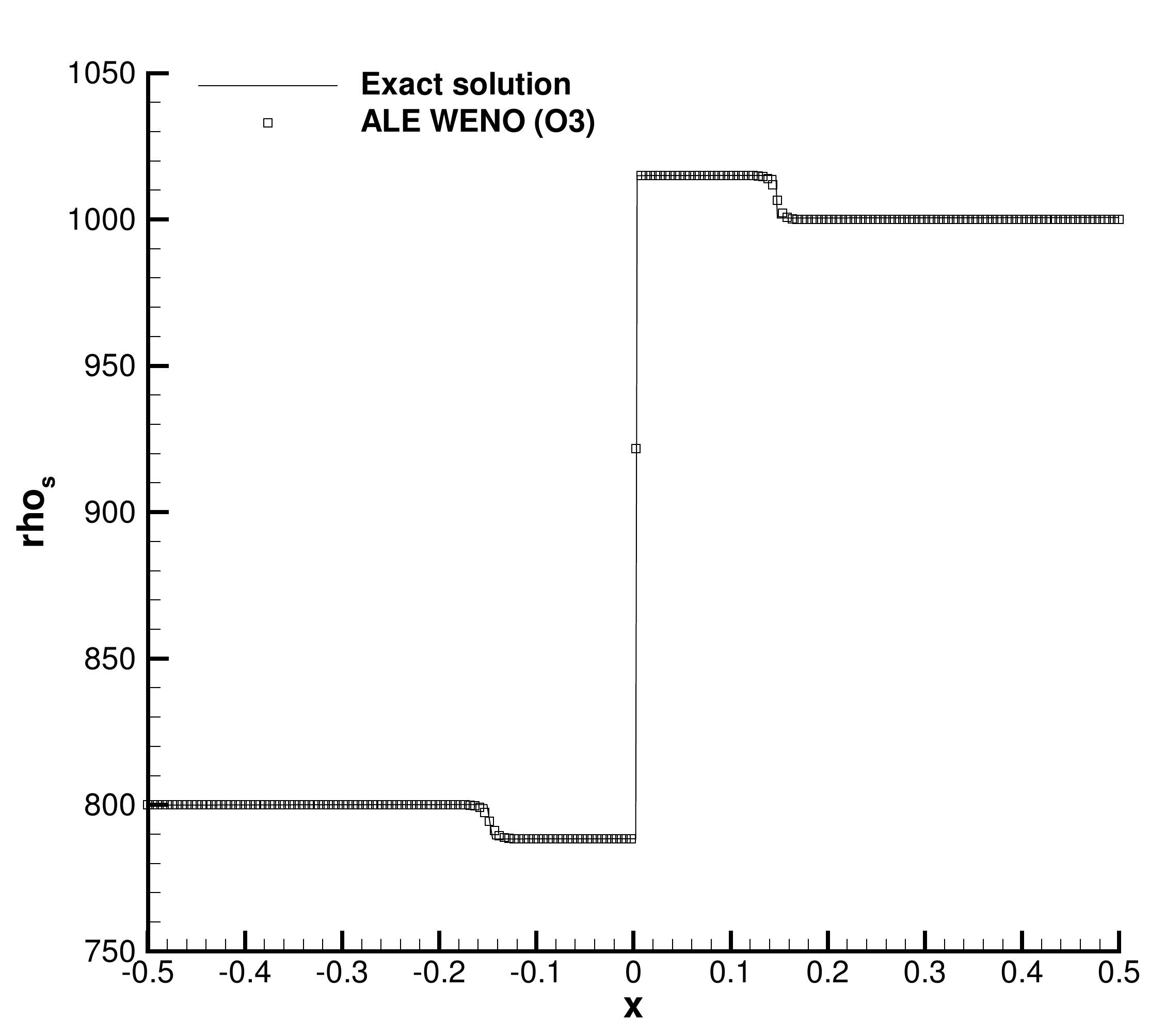}    & 
\includegraphics[width=0.4\textwidth]{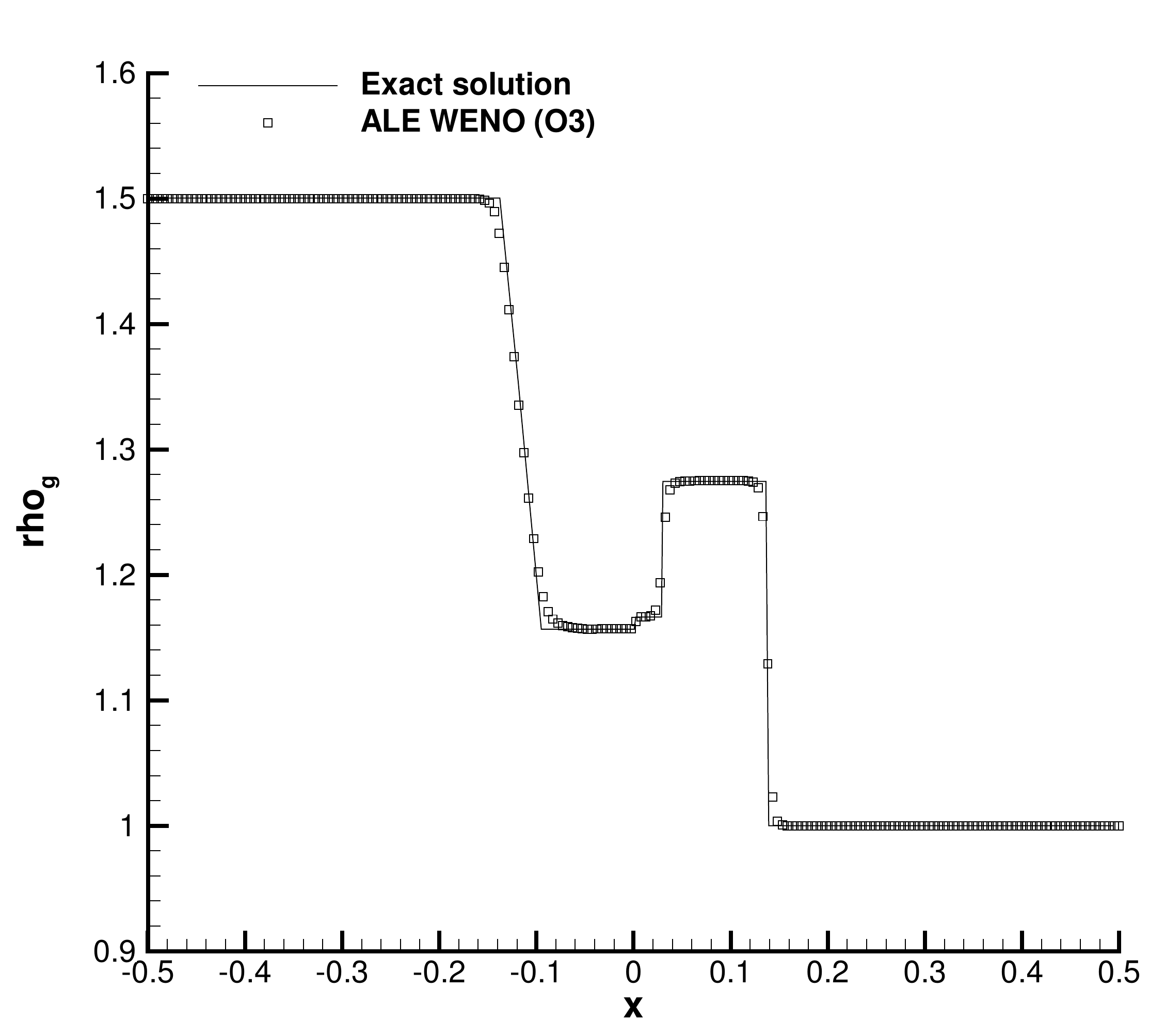}    \\ 
\includegraphics[width=0.4\textwidth]{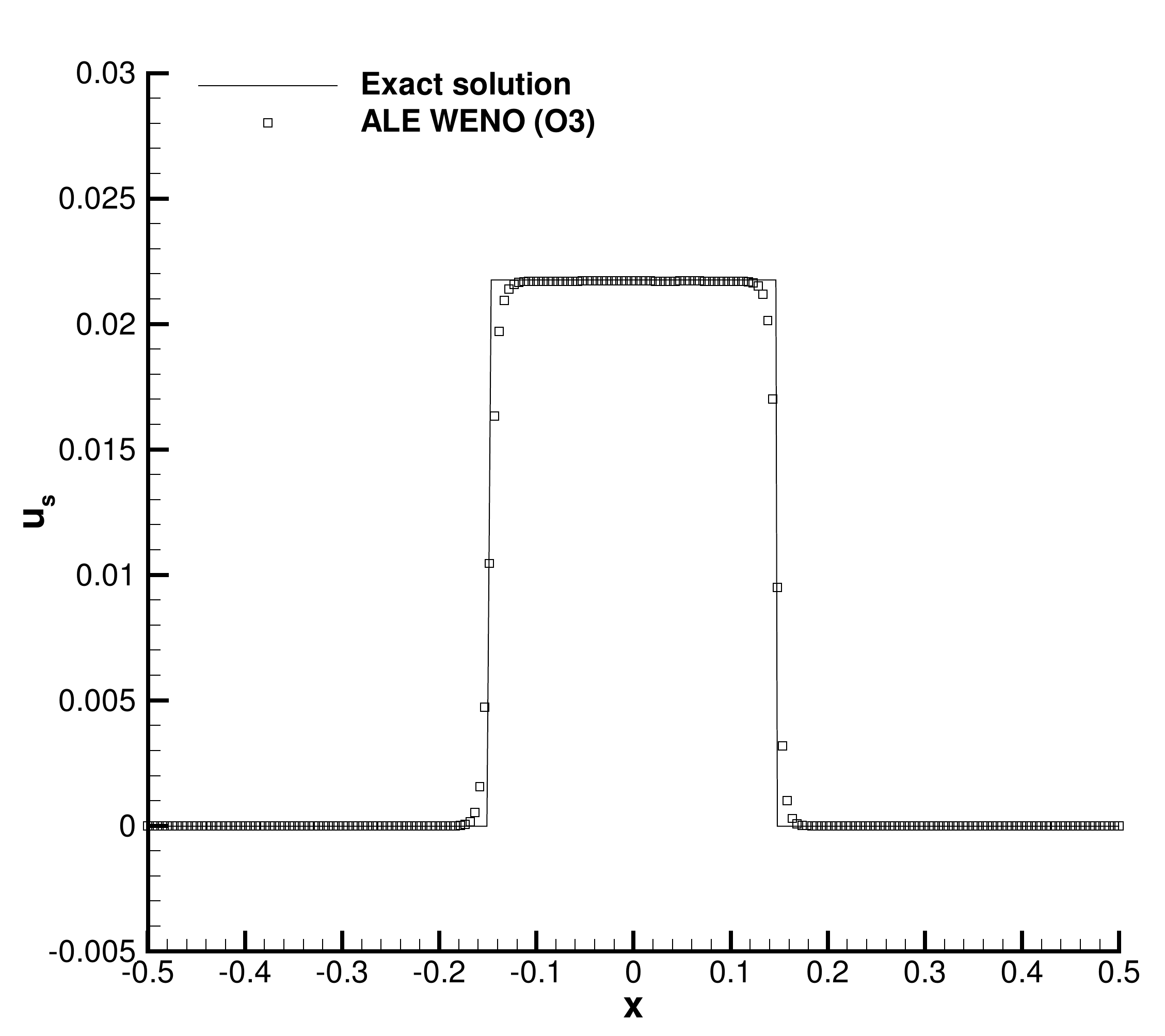}      &  
\includegraphics[width=0.4\textwidth]{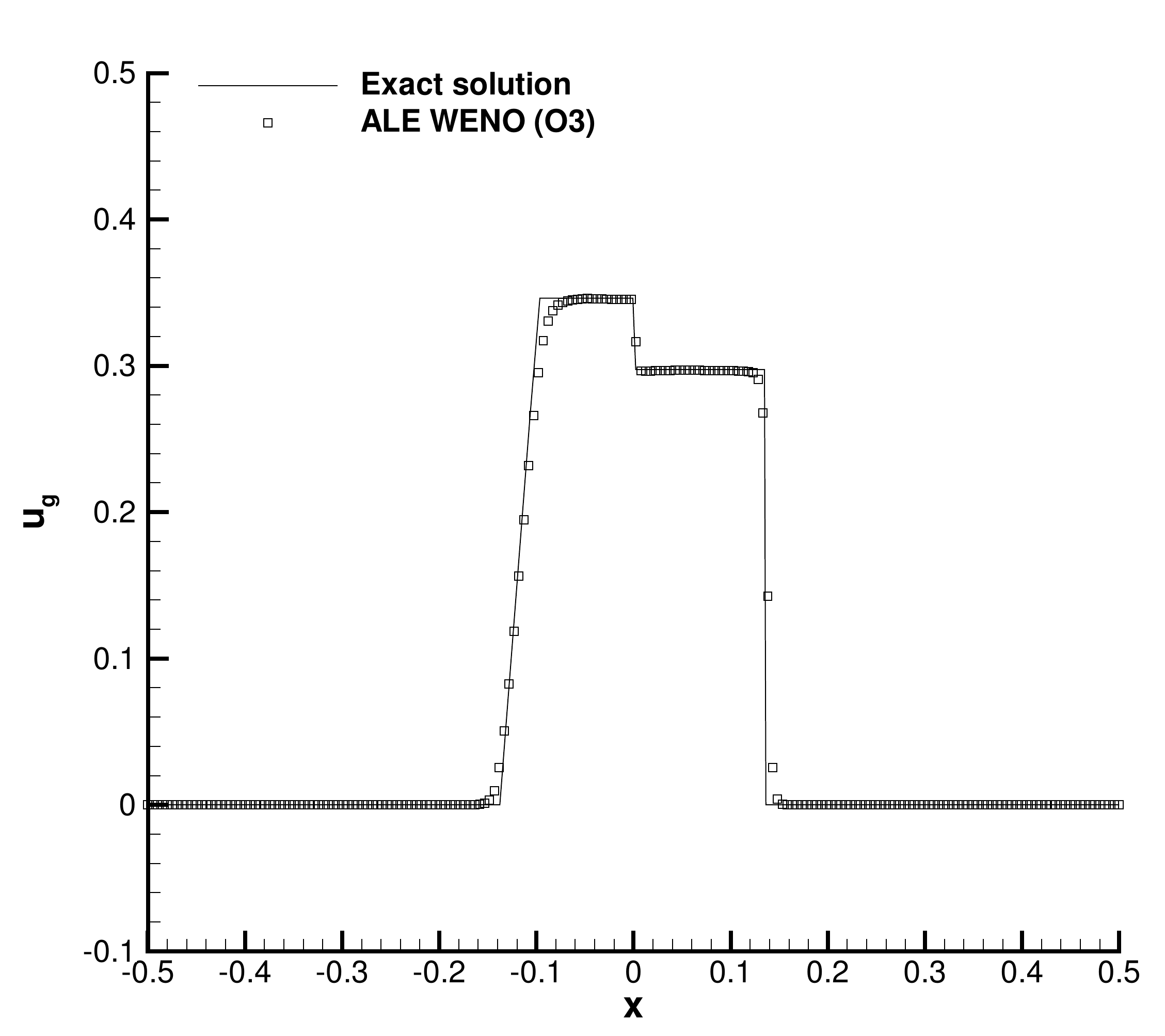}      \\ 
\includegraphics[width=0.4\textwidth]{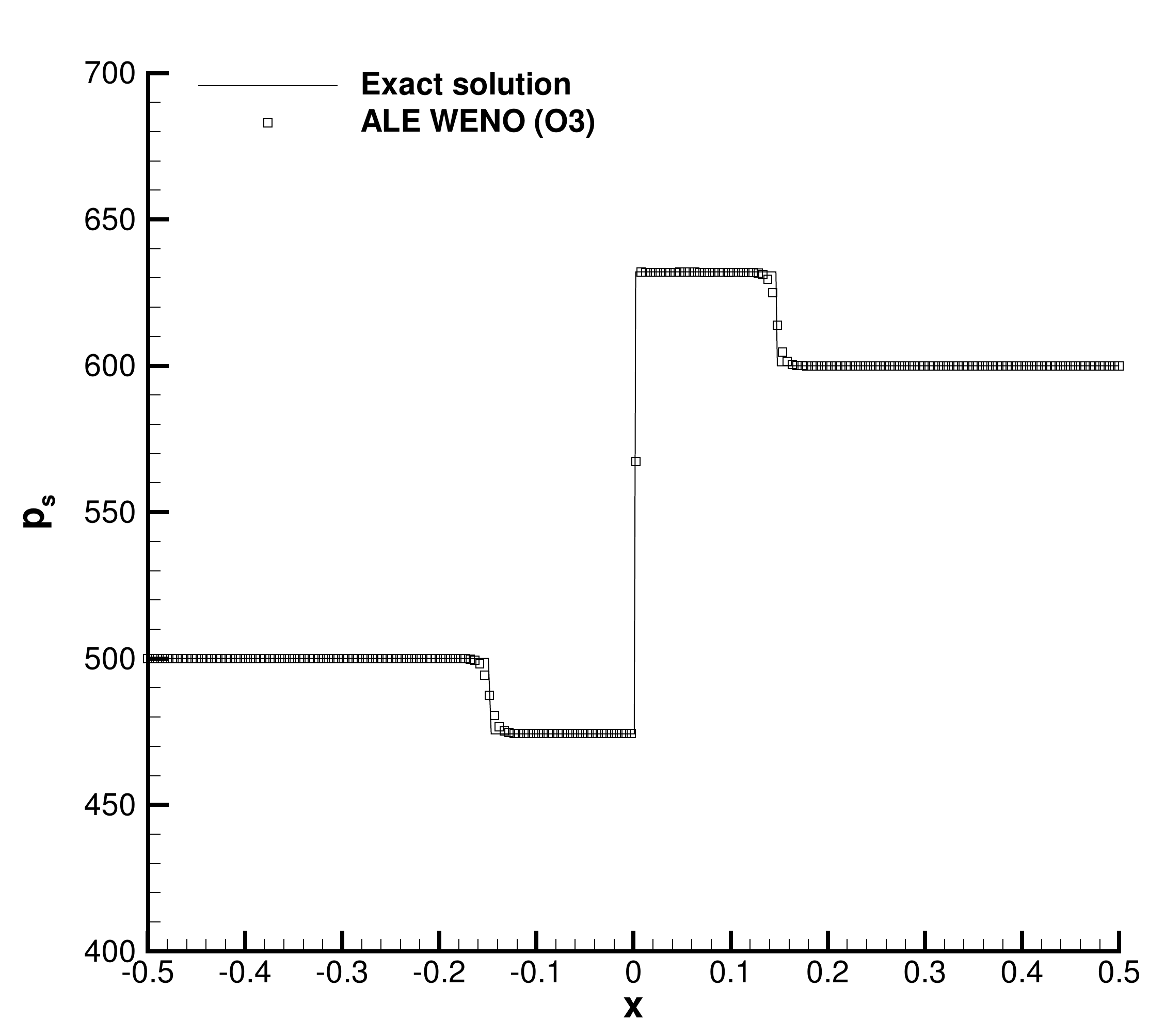}      & 
\includegraphics[width=0.4\textwidth]{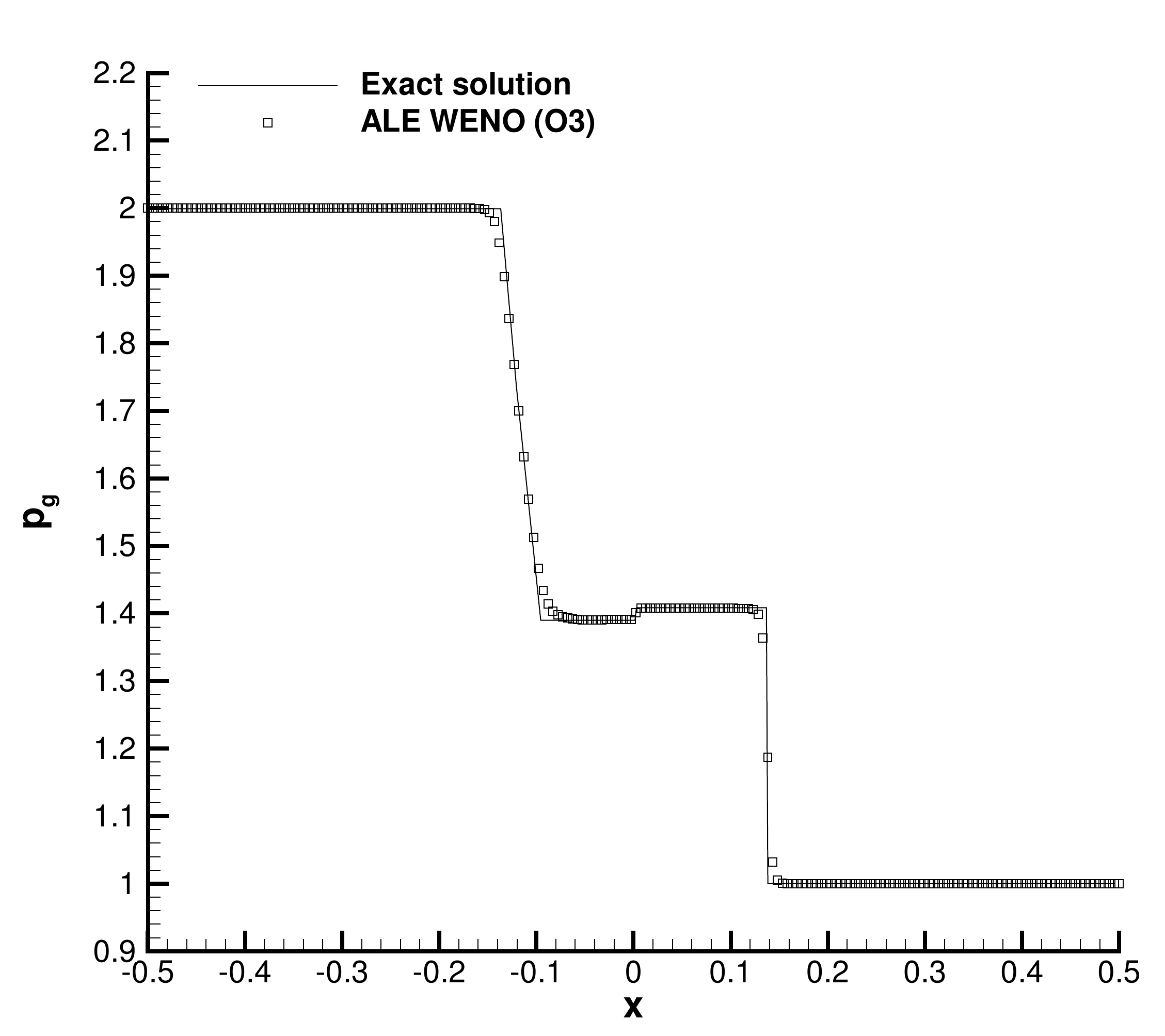}   
\end{tabular}
\caption{Results for Riemann problem RP2 of the seven-equation Baer-Nunziato model at time $t=0.1$.}
\label{fig.bn.rp2}
\end{center}
\end{figure}

\begin{figure}[!ht]
\begin{center}
\begin{tabular}{cc} 
\includegraphics[width=0.4\textwidth]{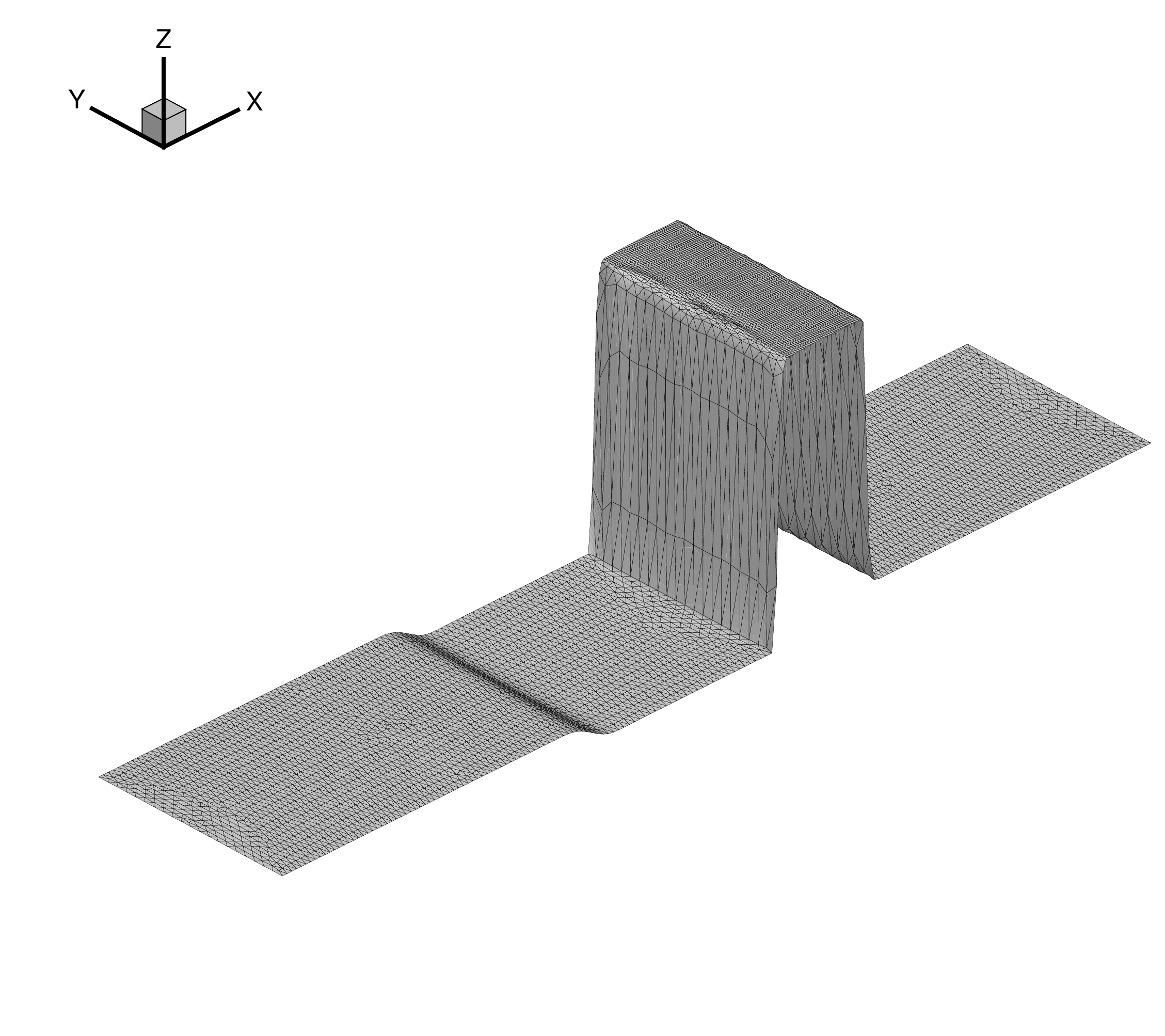}    & 
\includegraphics[width=0.4\textwidth]{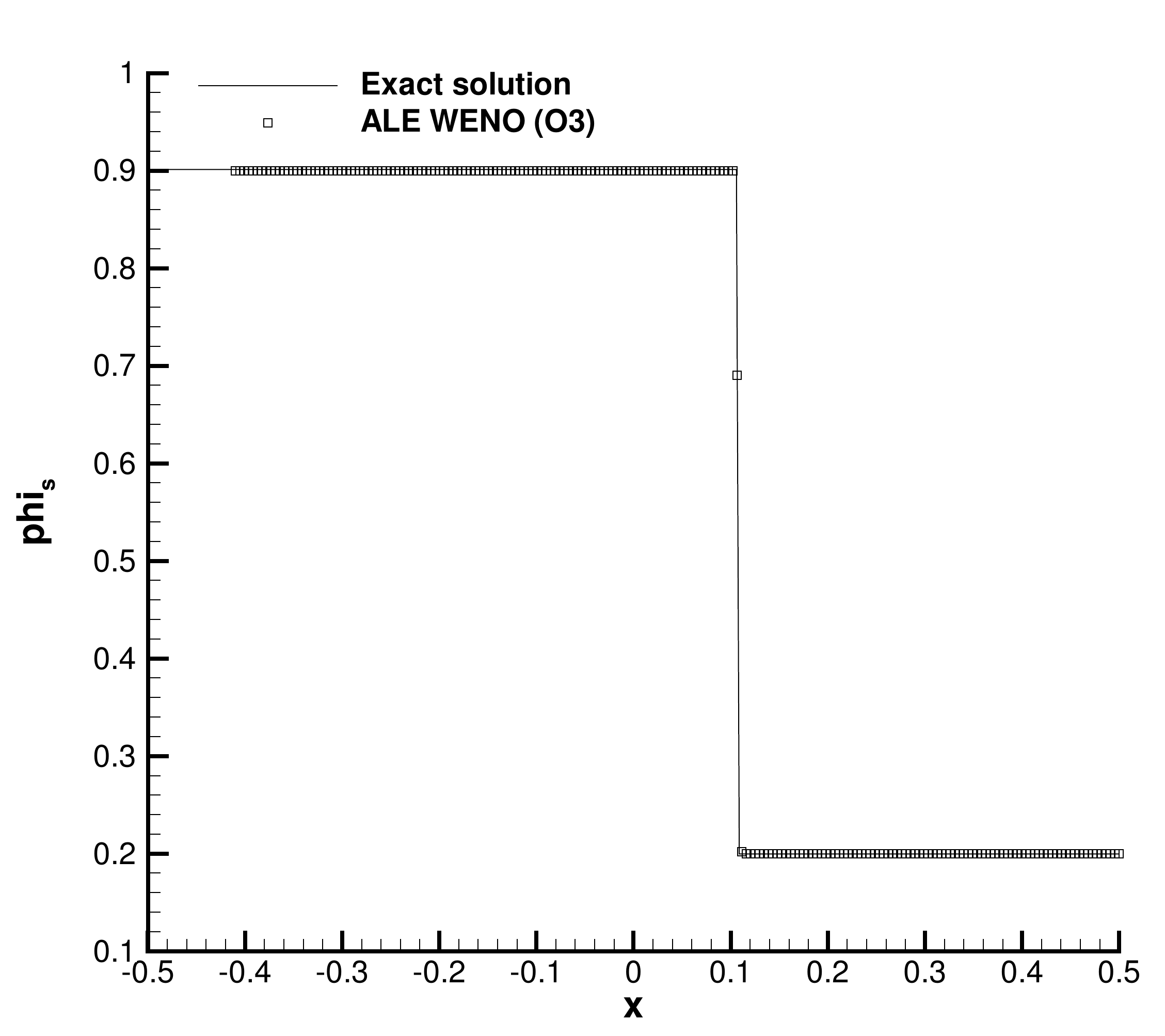}    \\ 
\includegraphics[width=0.4\textwidth]{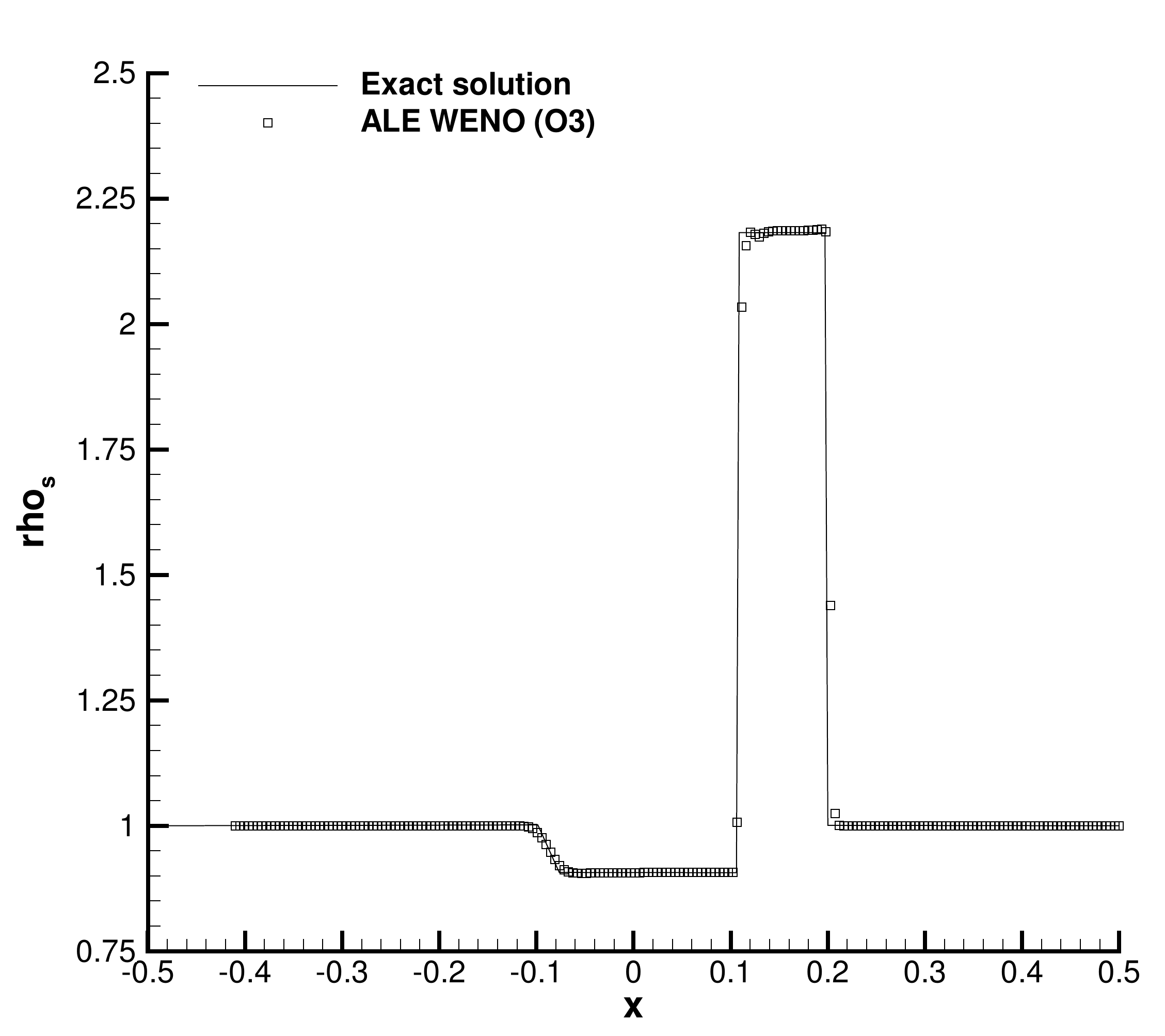}    & 
\includegraphics[width=0.4\textwidth]{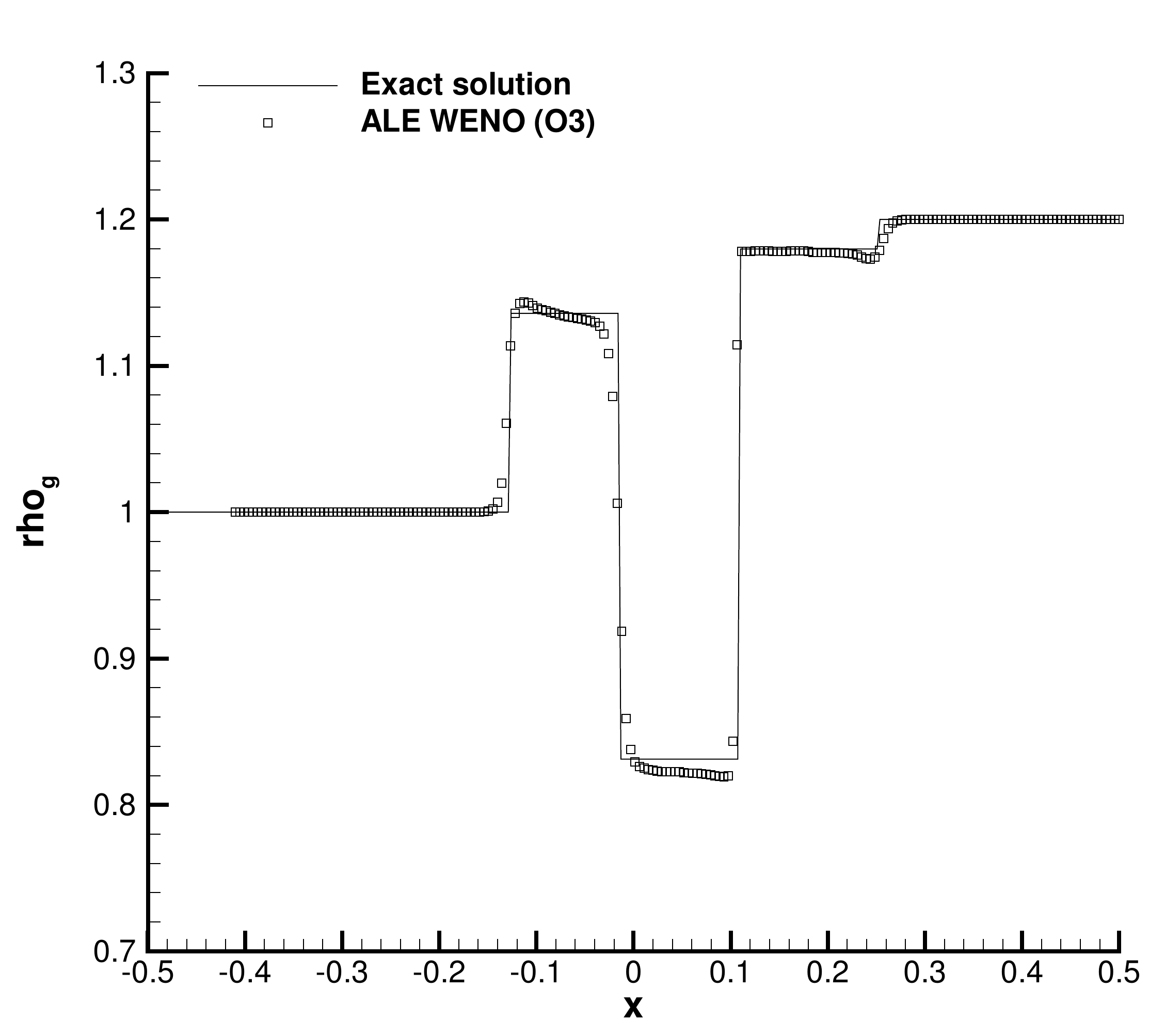}    \\ 
\includegraphics[width=0.4\textwidth]{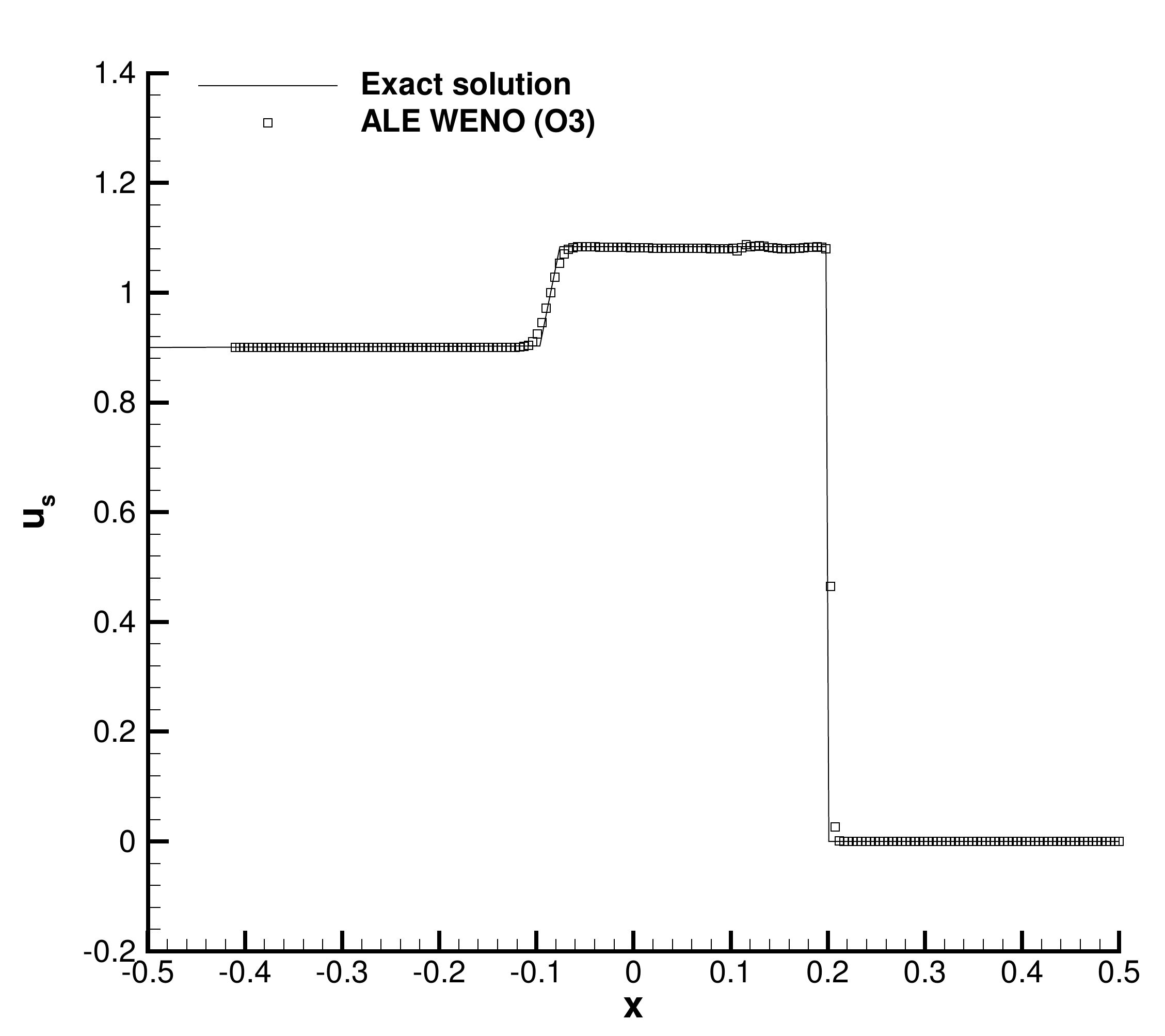}      &  
\includegraphics[width=0.4\textwidth]{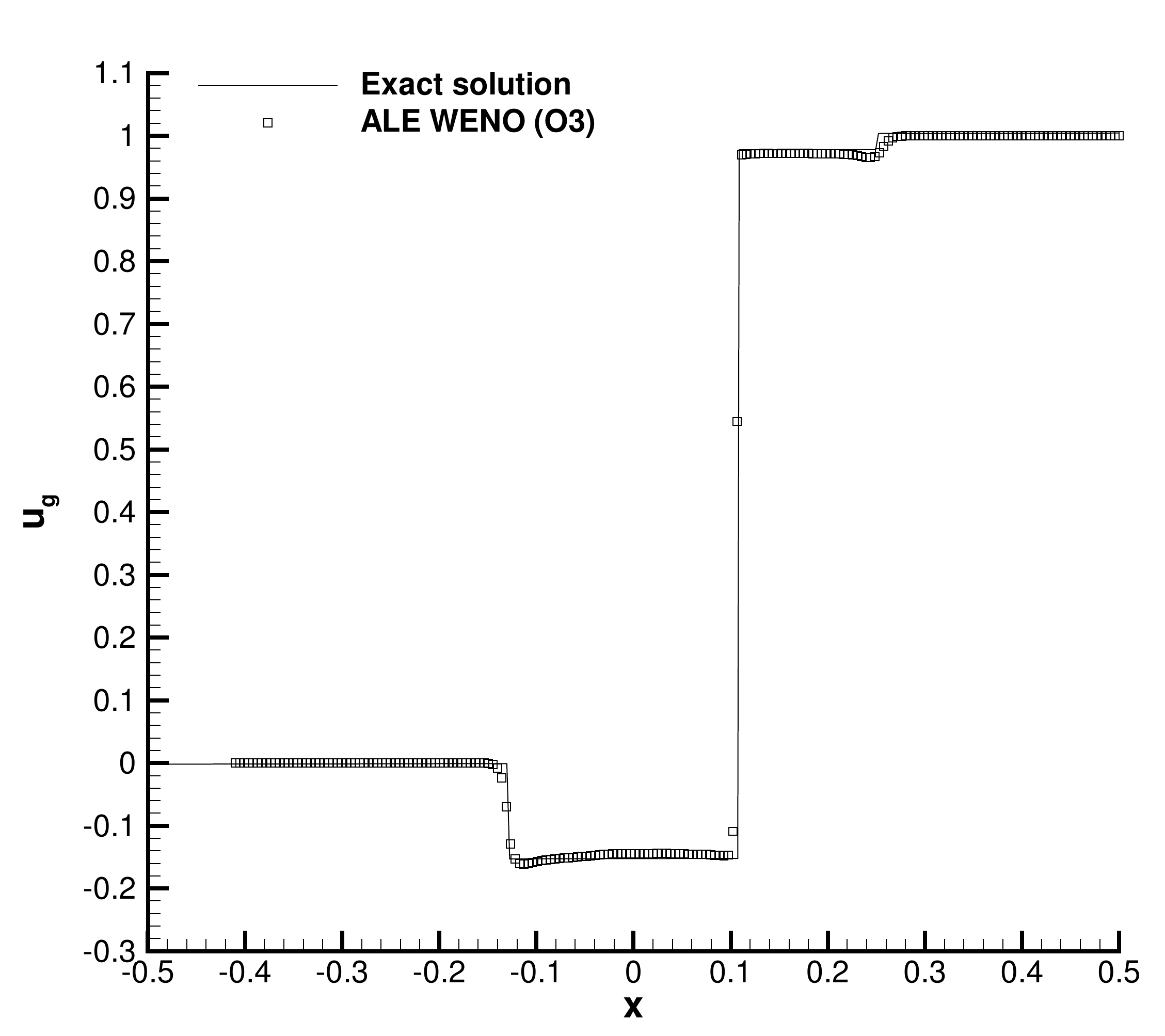}      \\ 
\includegraphics[width=0.4\textwidth]{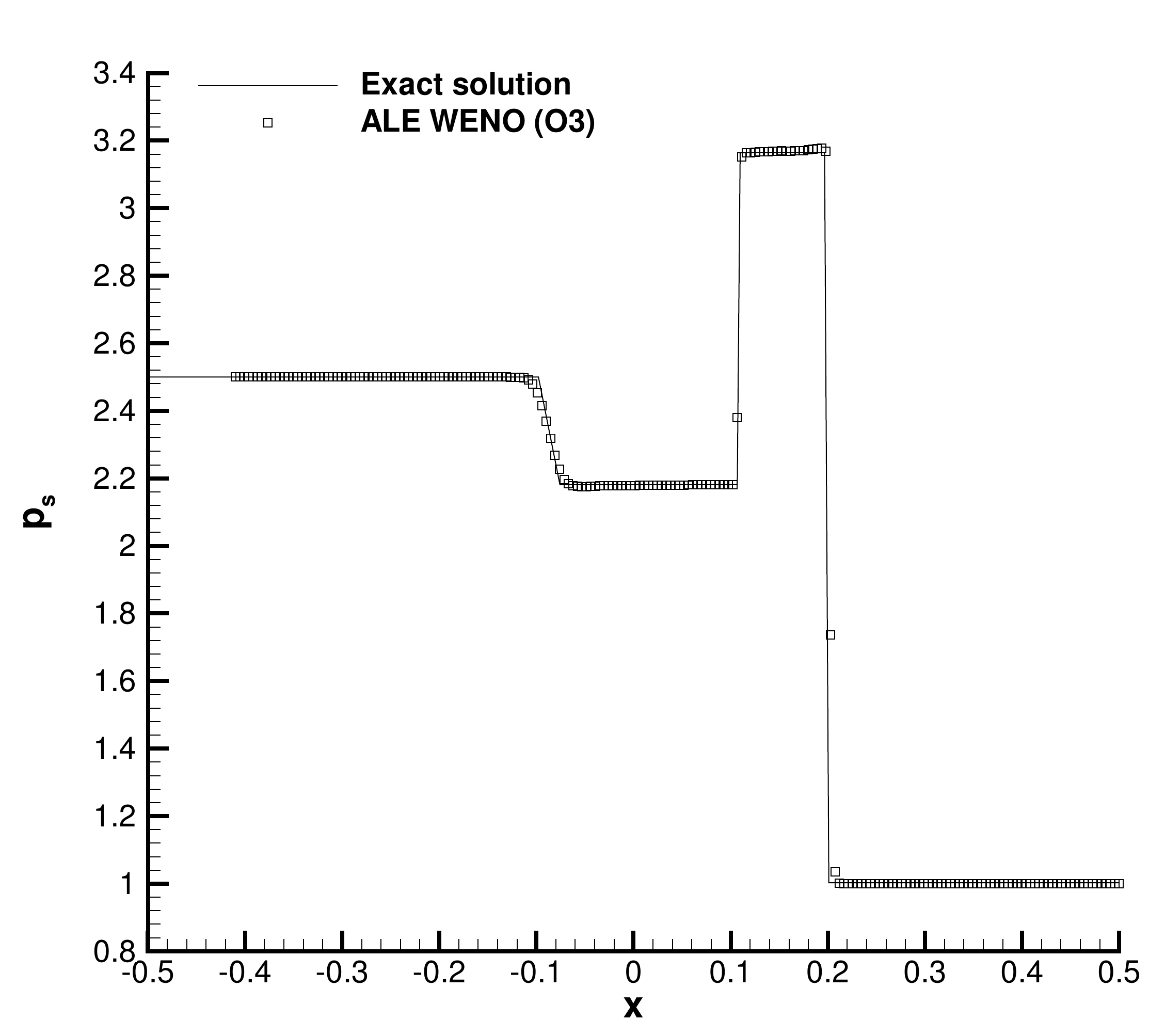}      & 
\includegraphics[width=0.4\textwidth]{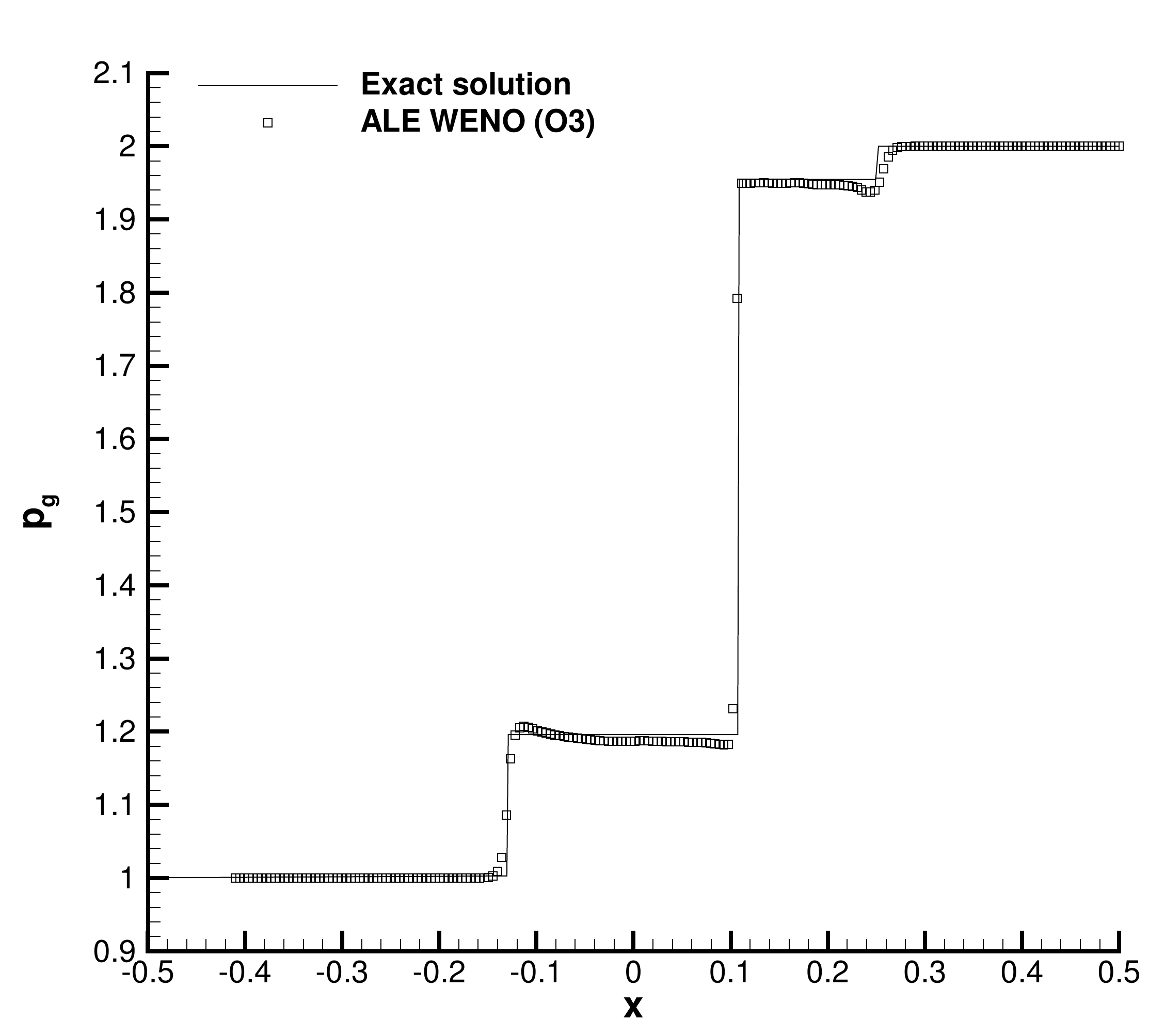}   
\end{tabular}
\caption{Results for Riemann problem RP3 of the seven-equation Baer-Nunziato model at time $t=0.1$.}
\label{fig.bn.rp3}
\end{center}
\end{figure}

\begin{figure}[!ht]
\begin{center}
\begin{tabular}{cc} 
\includegraphics[width=0.4\textwidth]{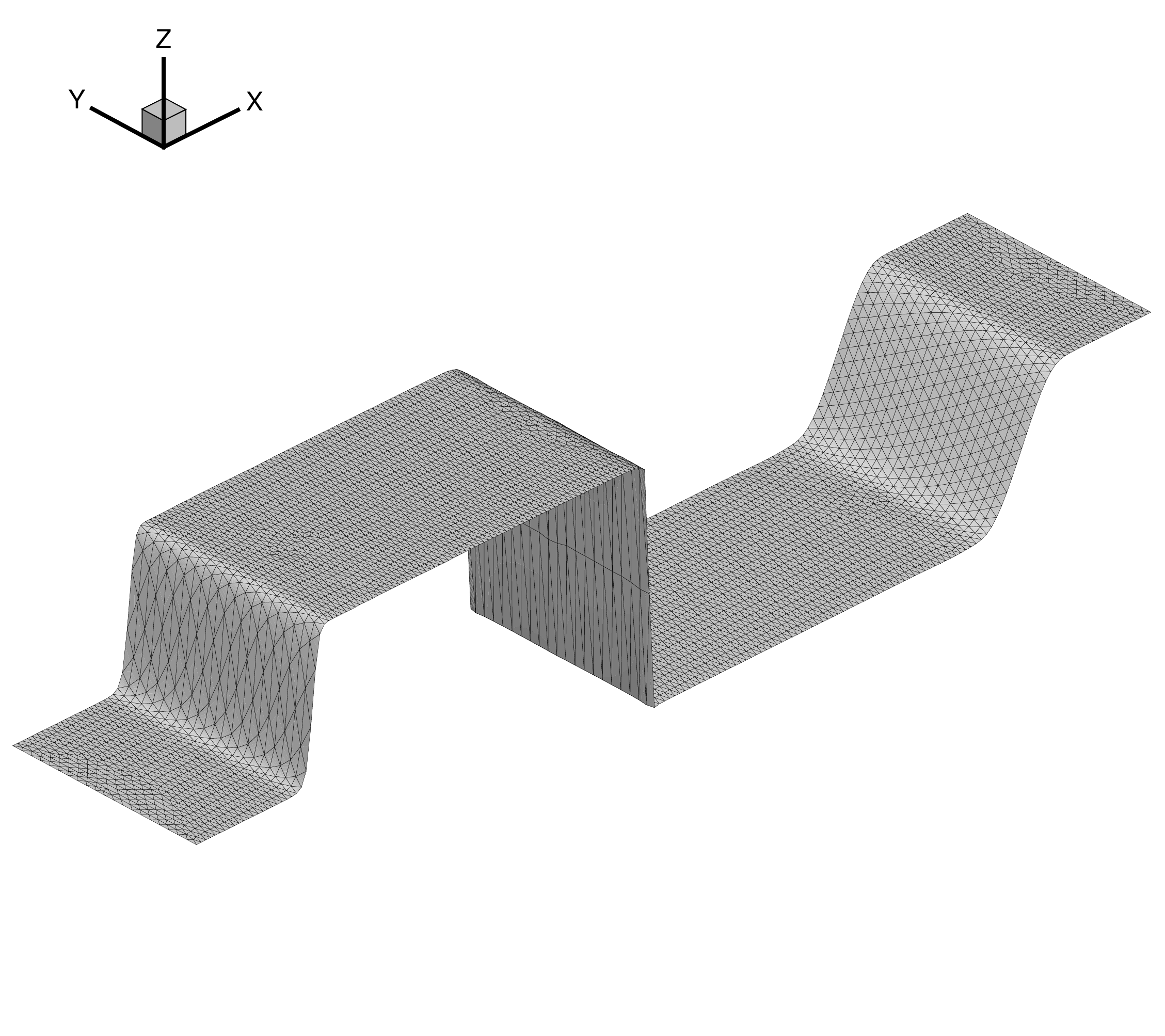}      & 
\includegraphics[width=0.4\textwidth]{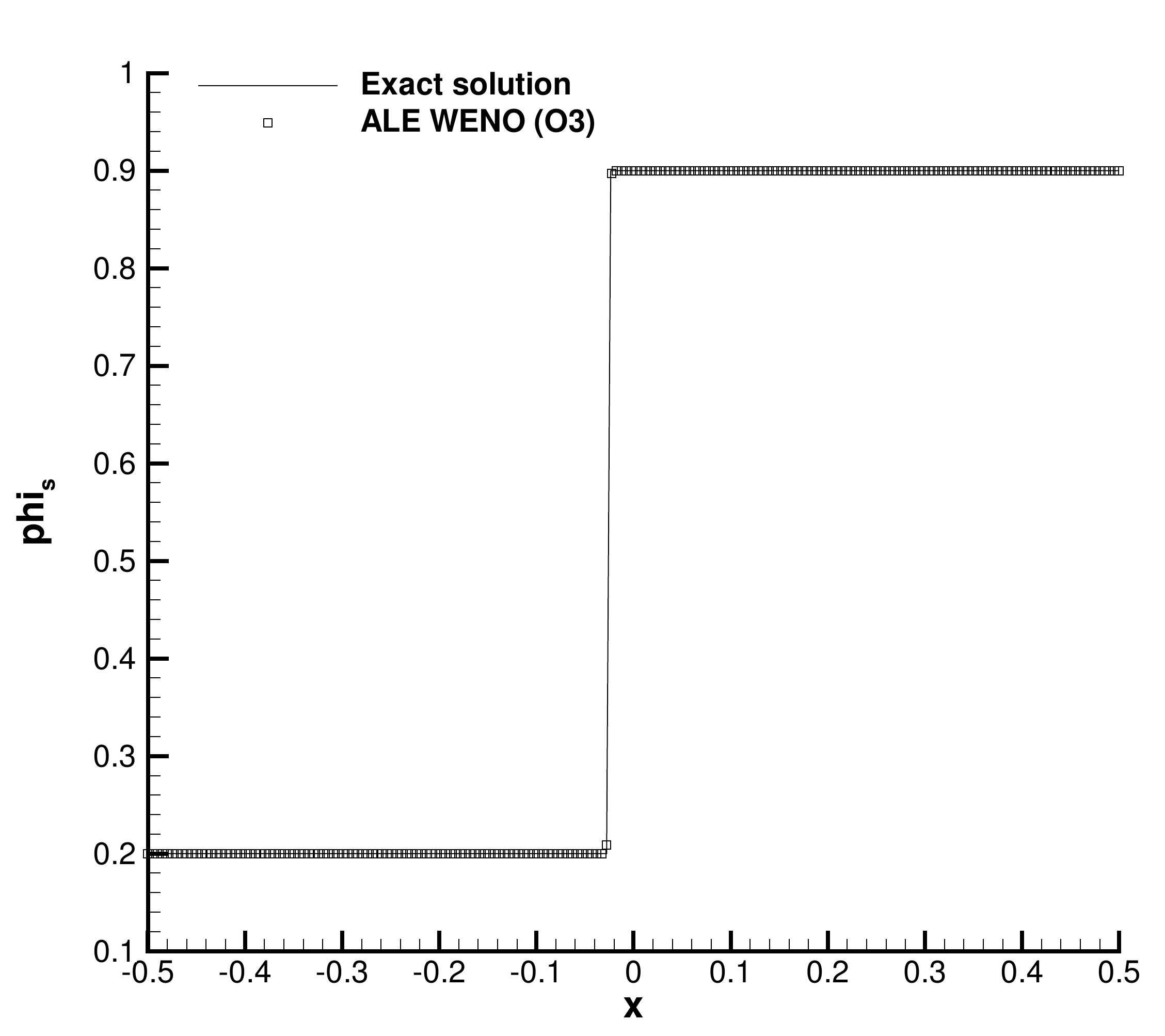}    \\ 
\includegraphics[width=0.4\textwidth]{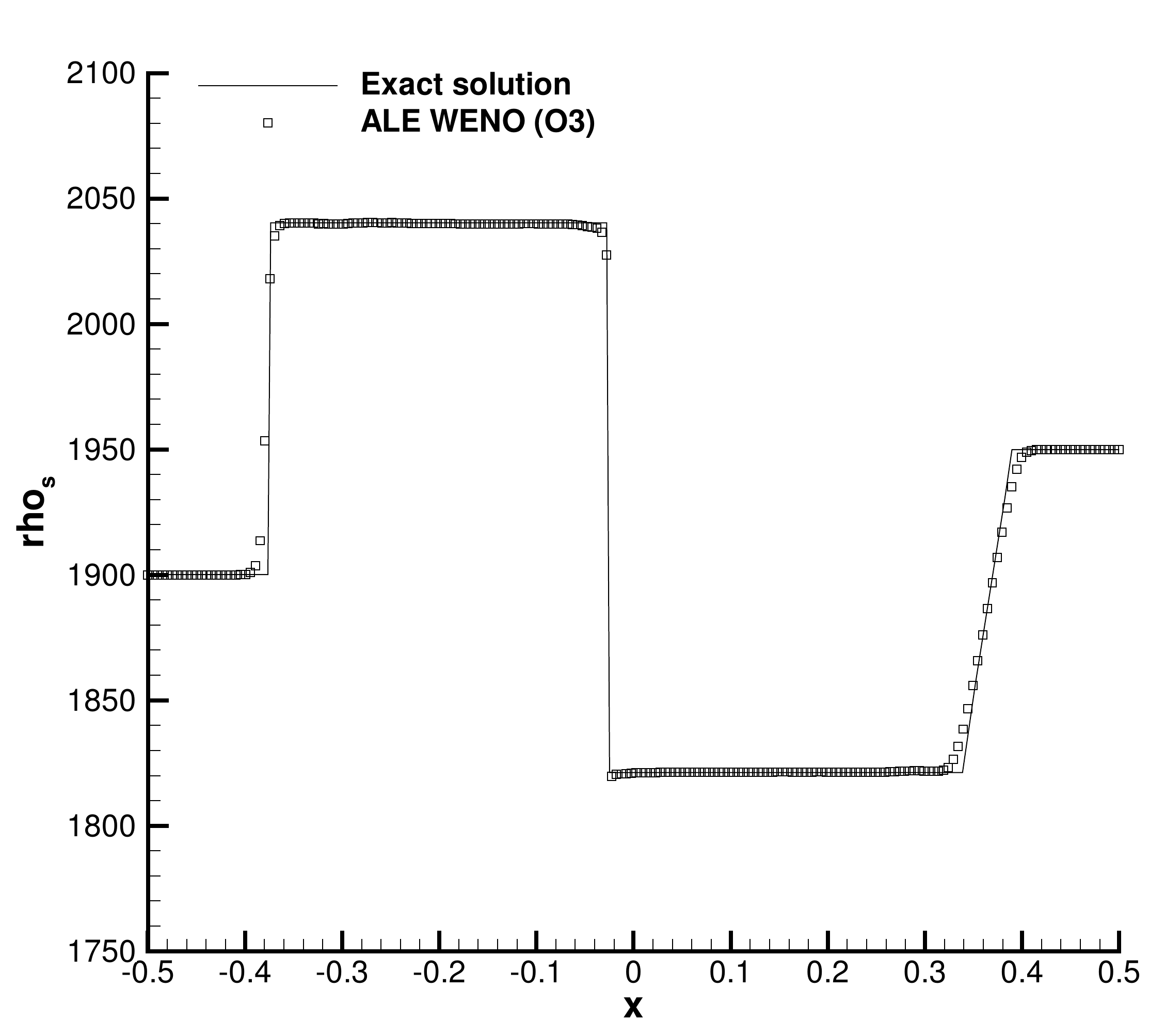}    & 
\includegraphics[width=0.4\textwidth]{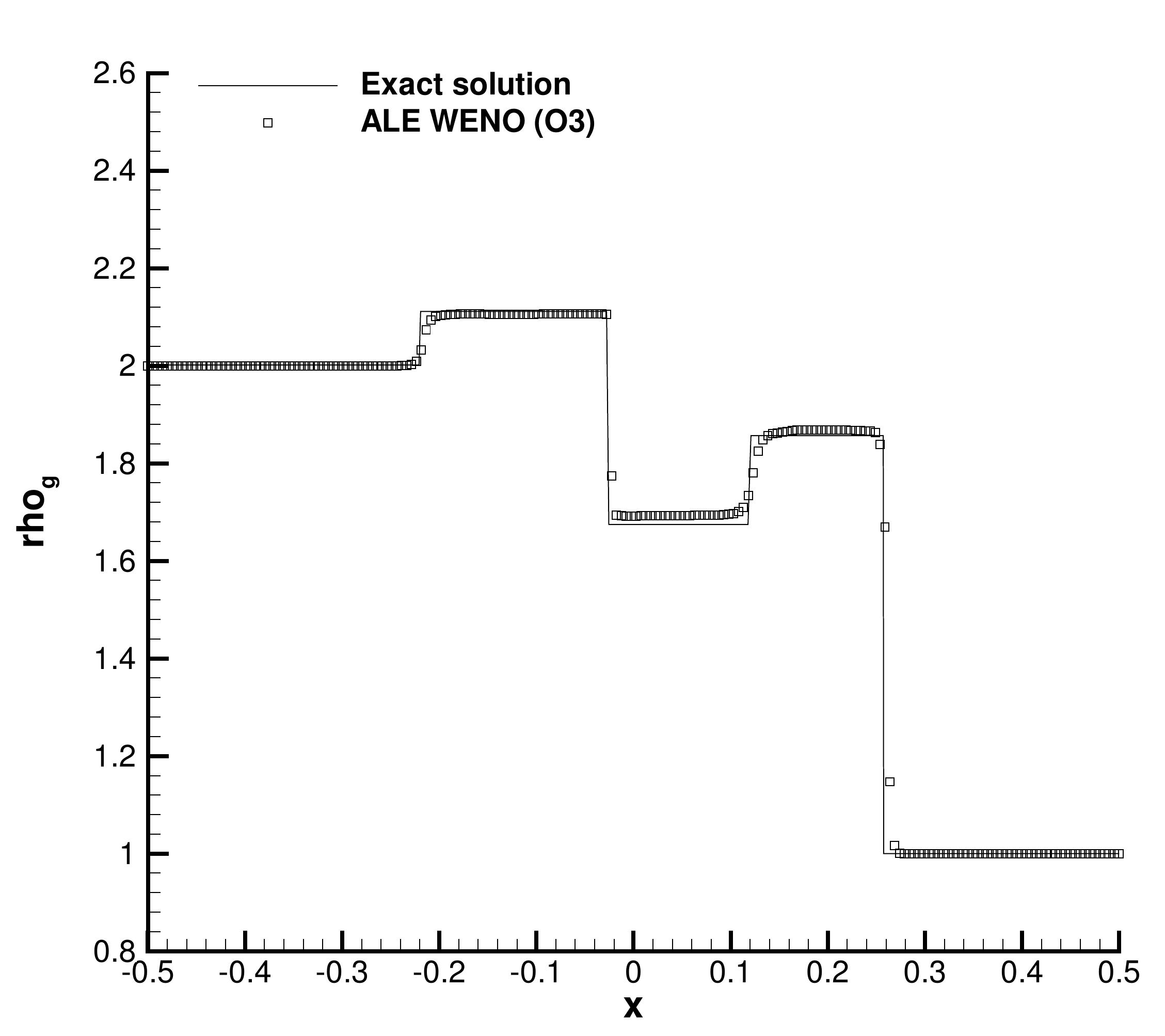}    \\ 
\end{tabular}
\caption{Results for Riemann problem RP4 of the seven-equation Baer-Nunziato model at time $t=0.15$.}
\label{fig.bn.rp4}
\end{center}
\end{figure}

The numerical results are shown in Figs. \ref{fig.bn.rp1} - \ref{fig.bn.rp4} and are compared
with the exact solution. On the top left of each figure a sketch of the mesh is depicted, while 
the other subfigures contain a one--dimensional cut through the reconstructed numerical solution $\w_h$ 
along the $x$-axis, evaluated at the final time on 200 \textit{equidistant} sample points. 
Due to the Lagrangian formulation of the method, the solid contact is resolved in a very sharp manner in all cases, 
which was actually the main aim in the design of a high order Lagrangian--type scheme for the compressible
Baer--Nunziato model. 
Also for the other waves we can note in general a very good agreement between our numerical results and the exact 
reference solutions given in  \cite{AndrianovWarnecke,Schwendeman,DeledicquePapalexandris}, apart for RP3, 
where a visible misfit between the exact solution and numerical solution is obtained for the gas pressure and 
the gas density. However, even a very diffusive central path--conservative FORCE scheme in the Eulerian 
framework \cite{USFORCE2} had difficulties with this test problem and produced visible spurious oscillations 
and a wrong position of the solid contact in the gas phase. Further investigations on this problem are necessary, 
in particular, whether a different choice of the path is able to improve the situation. 

\begin{table}[!b]
\caption{Initial states left (L) and right (R) for the Riemann problems solved in 2D and 3D with the 
Baer-Nunziato model. Values for $\gamma_i$, $\pi_i$ and the final time $t_e$ are also given.}
\renewcommand{\arraystretch}{1.0}
\begin{center}
\begin{tabular}{ccccccccc}
\hline
   & $\rho_s$ & $u_s$  & $p_s$ & $\rho_g$ & $u_g$ & $p_g$ & $\phi_s$ & $t_e$  \\
\hline 
\multicolumn{1}{l}{\textbf{RP1 \cite{DeledicquePapalexandris}:} } & 
\multicolumn{8}{c}{ $\gamma_s = 1.4, \quad \pi_s = 0, \quad \gamma_g = 1.4, \quad \pi_g = 0$}  \\
\hline 
L & 1.0    & 0.0   & 1.0  & 0.5 & 0.0   &  1.0 & 0.4 & 0.10 \\
R & 2.0    & 0.0   & 2.0  & 1.5 & 0.0   &  2.0 & 0.8 &      \\
\hline 
\multicolumn{1}{l}{\textbf{RP2 \cite{DeledicquePapalexandris}:}} & 
\multicolumn{8}{c}{ $\gamma_s = 3.0, \quad \pi_s = 100, \quad \gamma_g = 1.4, \quad \pi_g = 0$}  \\
\hline
L & 800.0   & 0.0   & 500.0  & 1.5 & 0.0   & 2.0 & 0.4 & 0.10  \\
R & 1000.0  & 0.0   & 600.0  & 1.0 & 0.0   & 1.0 & 0.3 &       \\
\hline 
\multicolumn{1}{l}{\textbf{RP3 \cite{DeledicquePapalexandris}:}} & 
\multicolumn{8}{c}{ $\gamma_s = 1.4, \quad \pi_s = 0, \quad \gamma_g = 1.4, \quad \pi_g = 0$}  \\ 
\hline
L & 1.0     & 0.9       & 2.5      & 1.0       & 0.0      &  1.0 & 0.9   & 0.10   \\
R & 1.0     & 0.0       & 1.0      & 1.2       & 1.0      &  2.0 & 0.2   &        \\
\hline 
\multicolumn{1}{l}{\textbf{RP4 \cite{Schwendeman}:}} & 
\multicolumn{8}{c}{ $\gamma_s = 3.0, \quad \pi_s = 3400, \quad \gamma_g = 1.35, \quad \pi_g = 0$}  \\
\hline
L & 1900.0   & 0.0   & 10.0    & 2.0 & 0.0   & 3.0 & 0.2 & 0.15   \\
R & 1950.0   & 0.0   & 1000.0  & 1.0 & 0.0   & 1.0 & 0.9 &       \\
\hline 
\end{tabular}
\end{center}
\label{tab.rpbn.ic}
\end{table}

\subsection{Cylindrical Explosion Problems} 
\label{sec.ep}

\begin{figure}[!htbp]
\begin{center}
\begin{tabular}{cc} 
\includegraphics[width=0.4\textwidth]{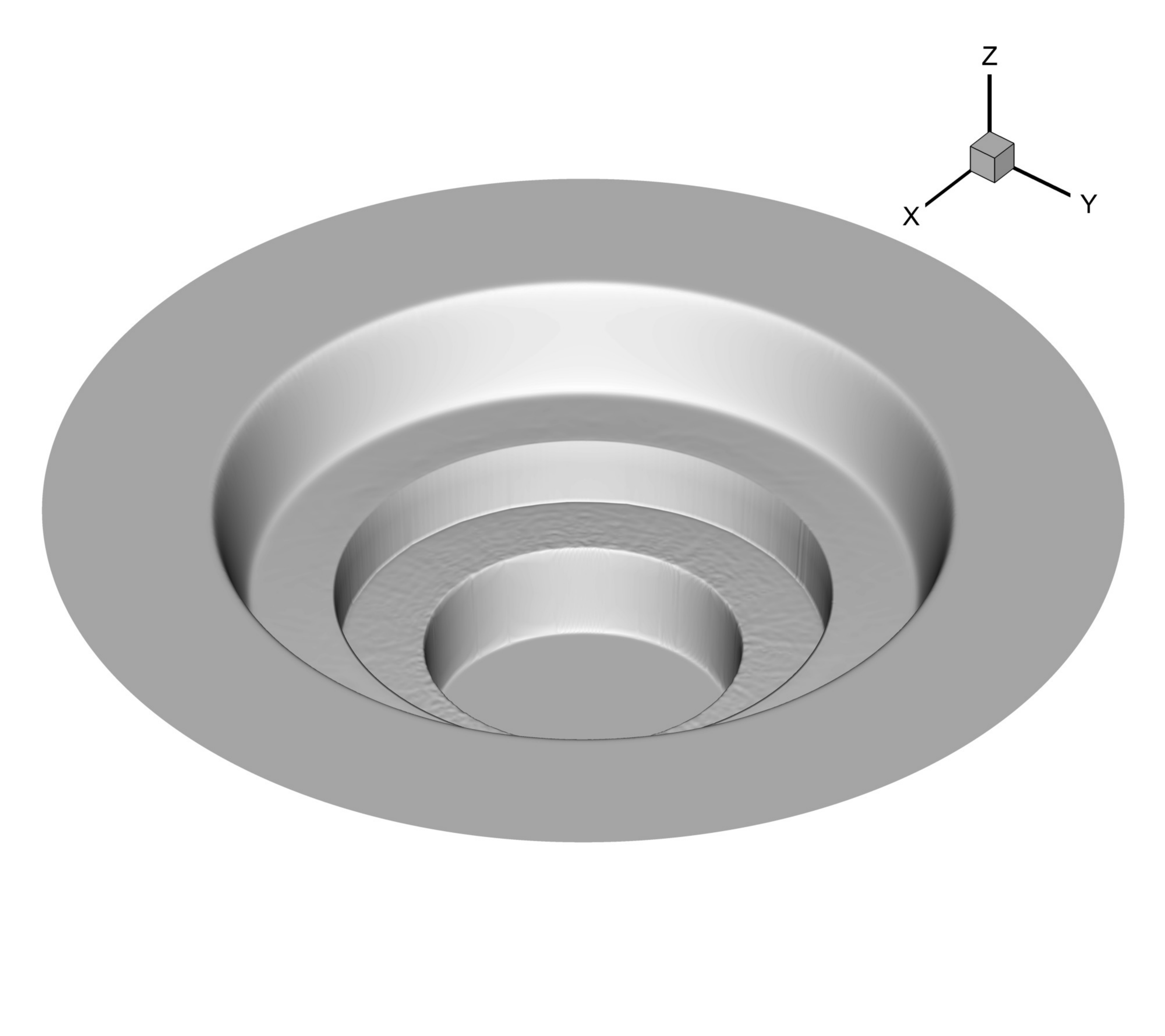}       & 
\includegraphics[width=0.4\textwidth]{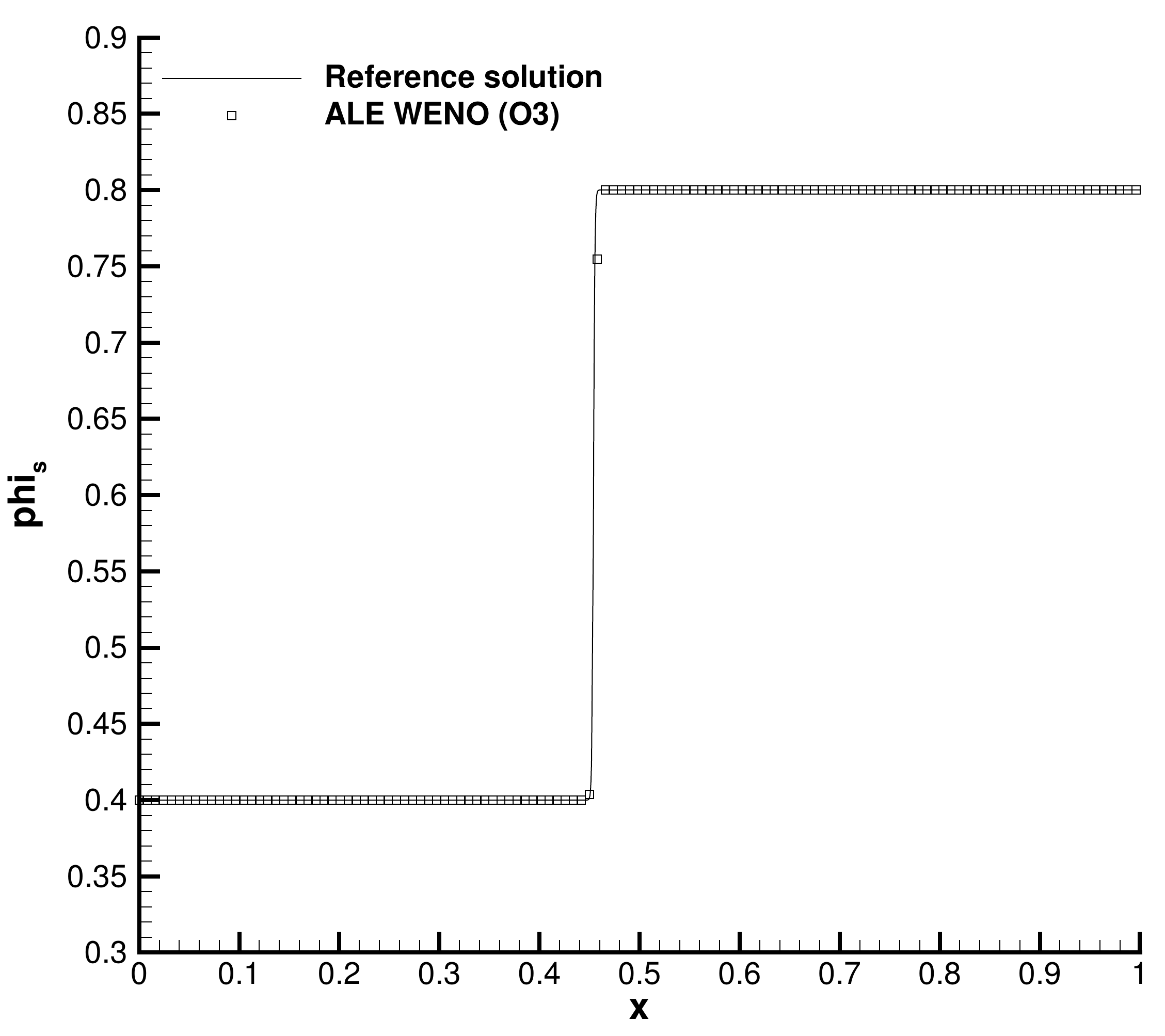}    \\ 
\includegraphics[width=0.4\textwidth]{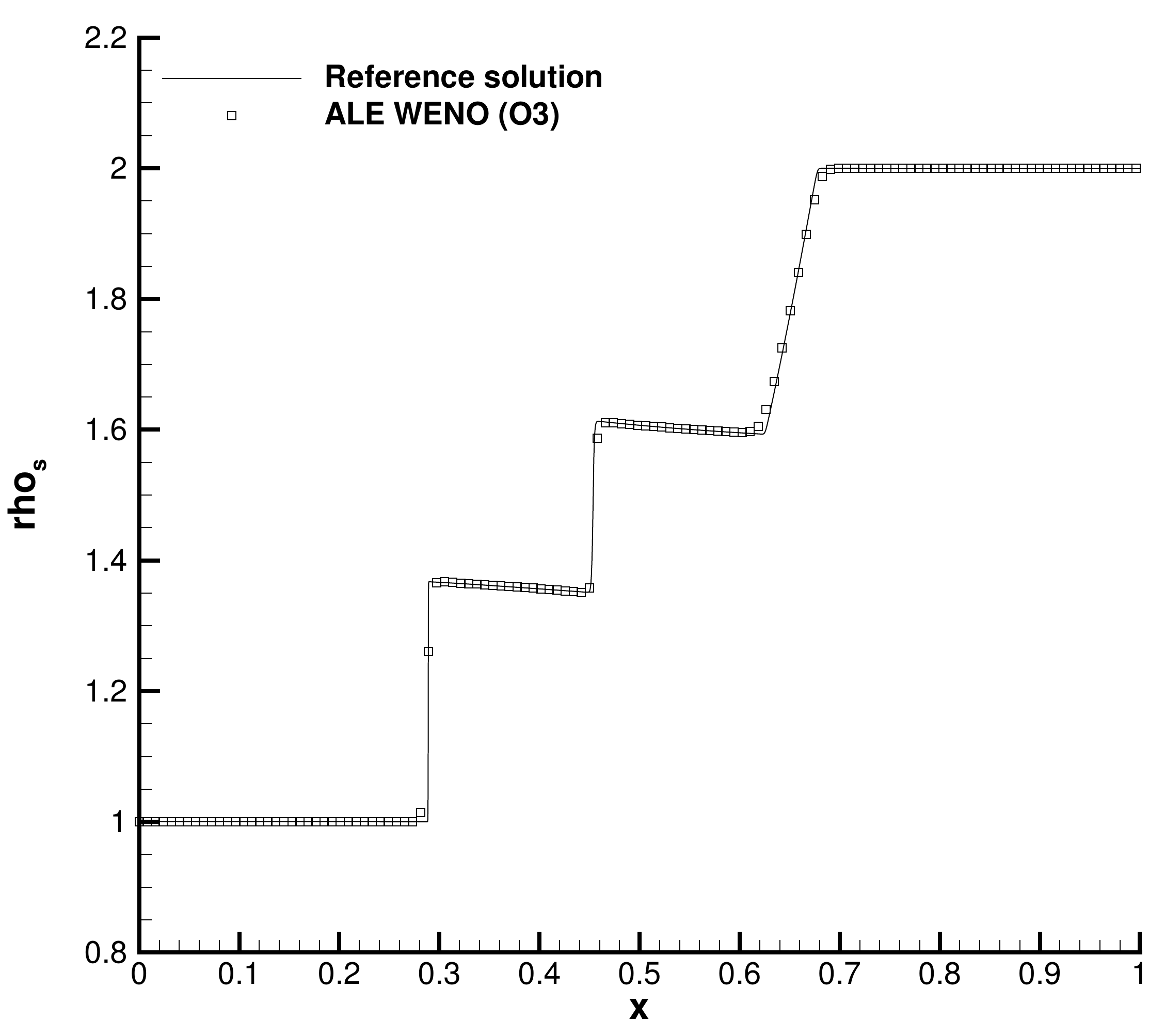}     & 
\includegraphics[width=0.4\textwidth]{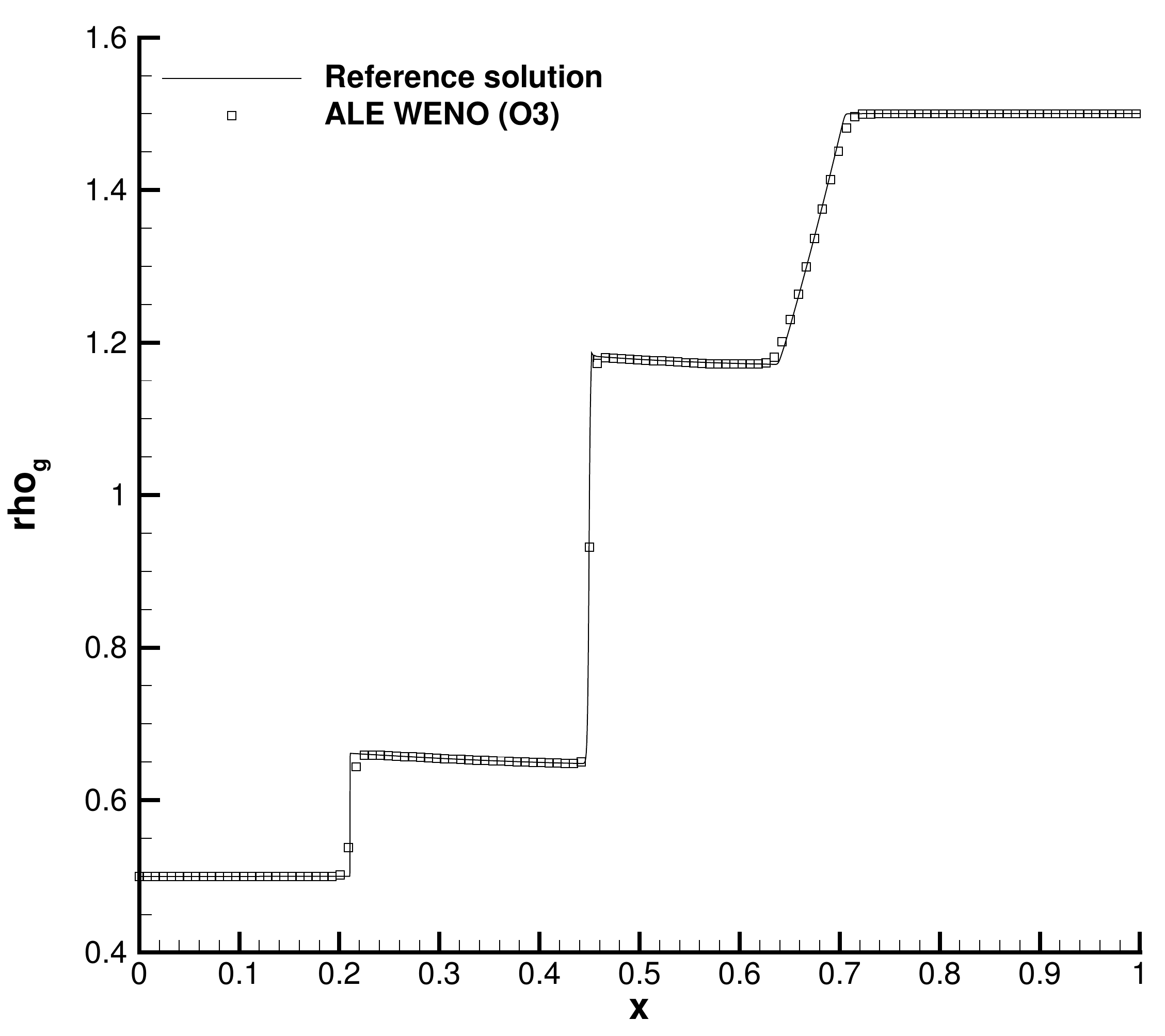}    \\ 
\includegraphics[width=0.4\textwidth]{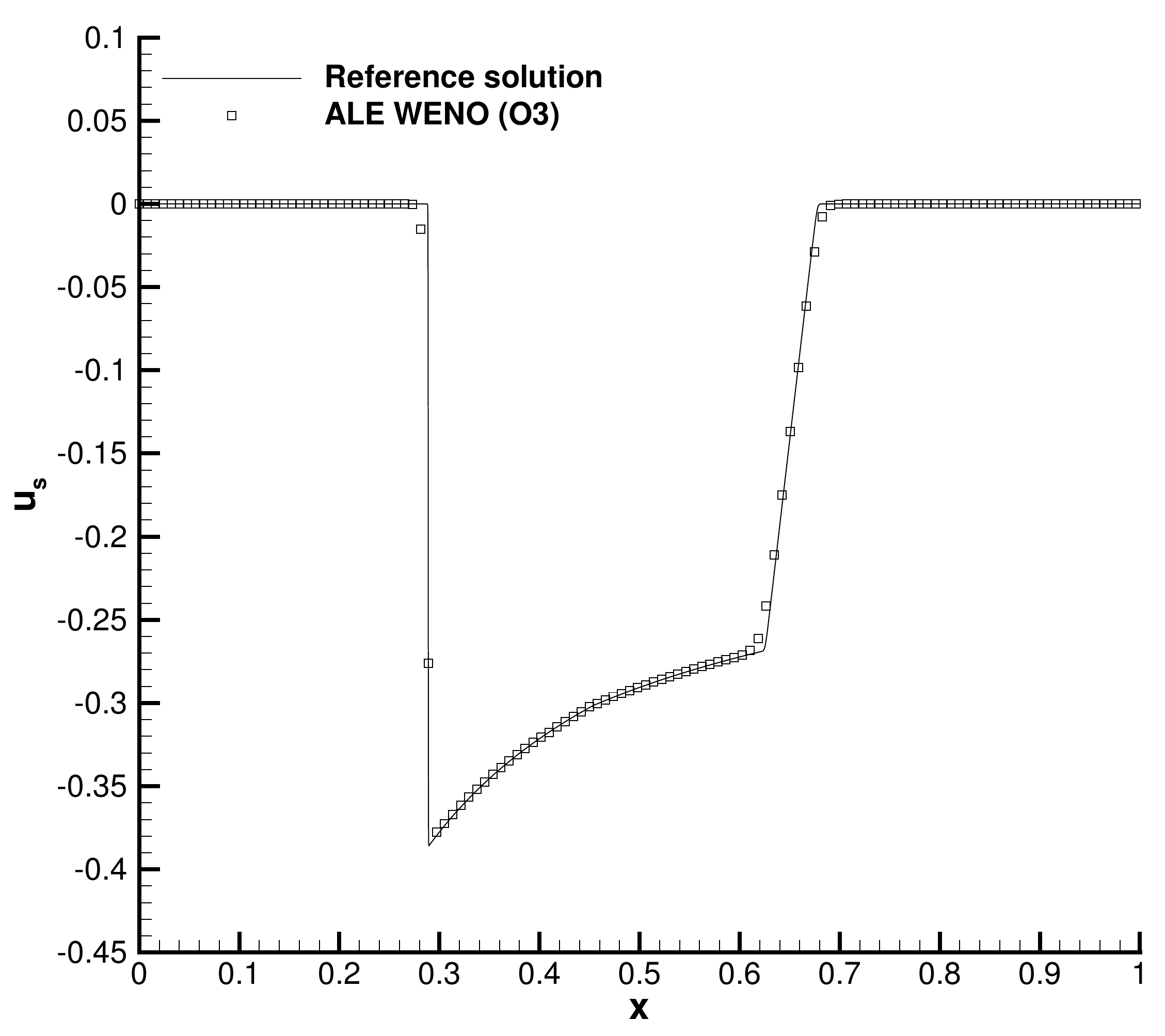}       & 
\includegraphics[width=0.4\textwidth]{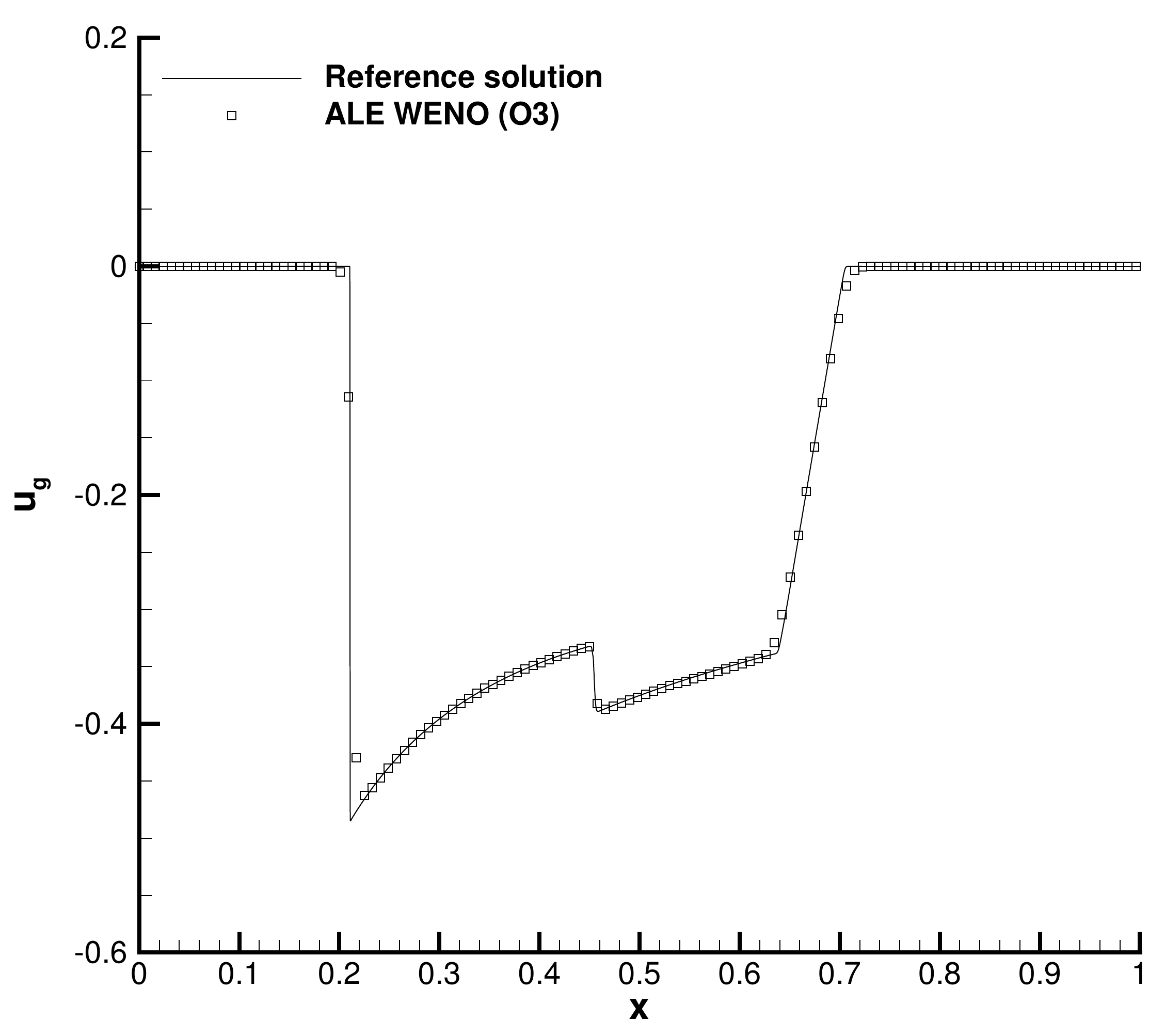}      \\ 
\includegraphics[width=0.4\textwidth]{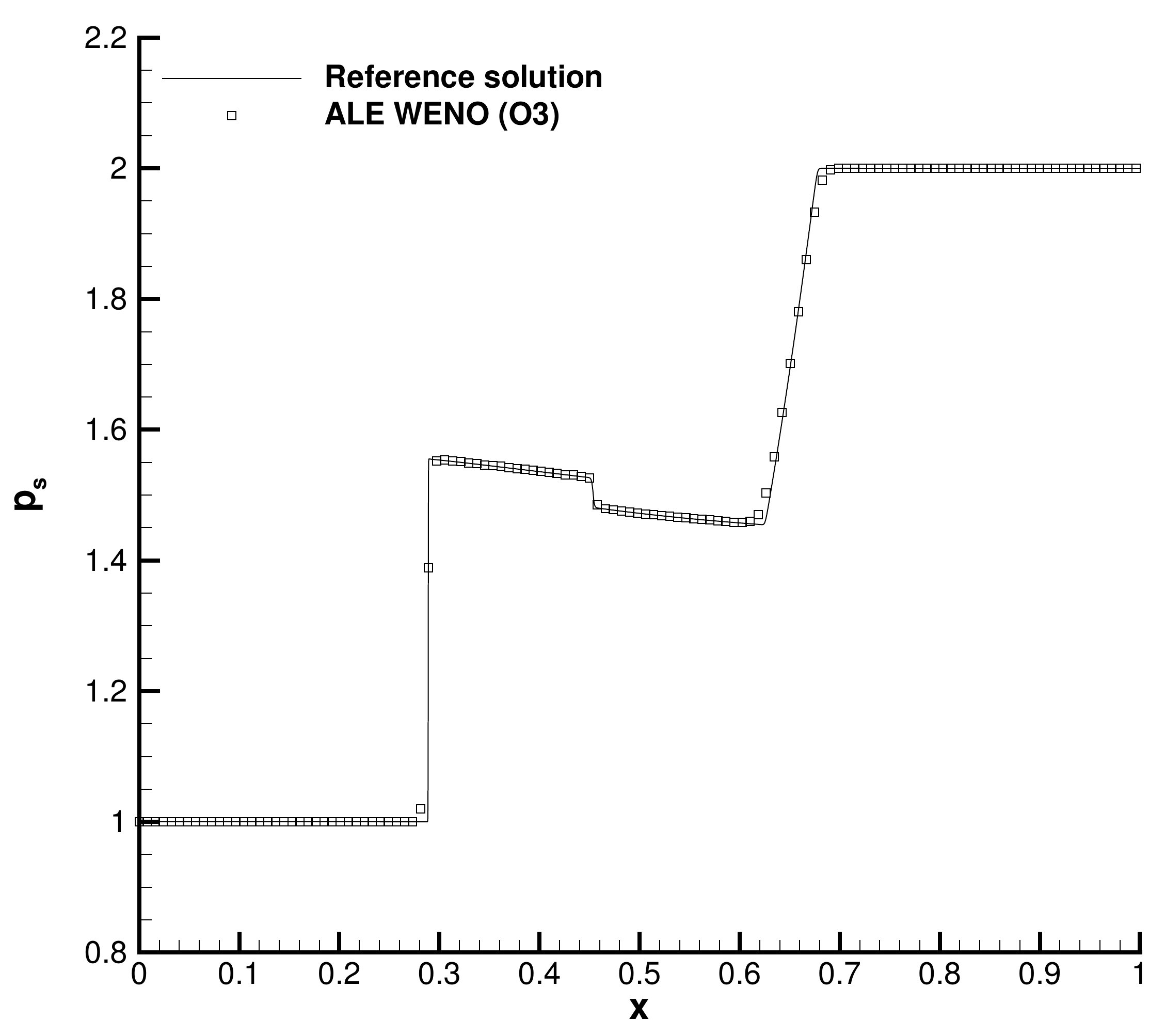}       & 
\includegraphics[width=0.4\textwidth]{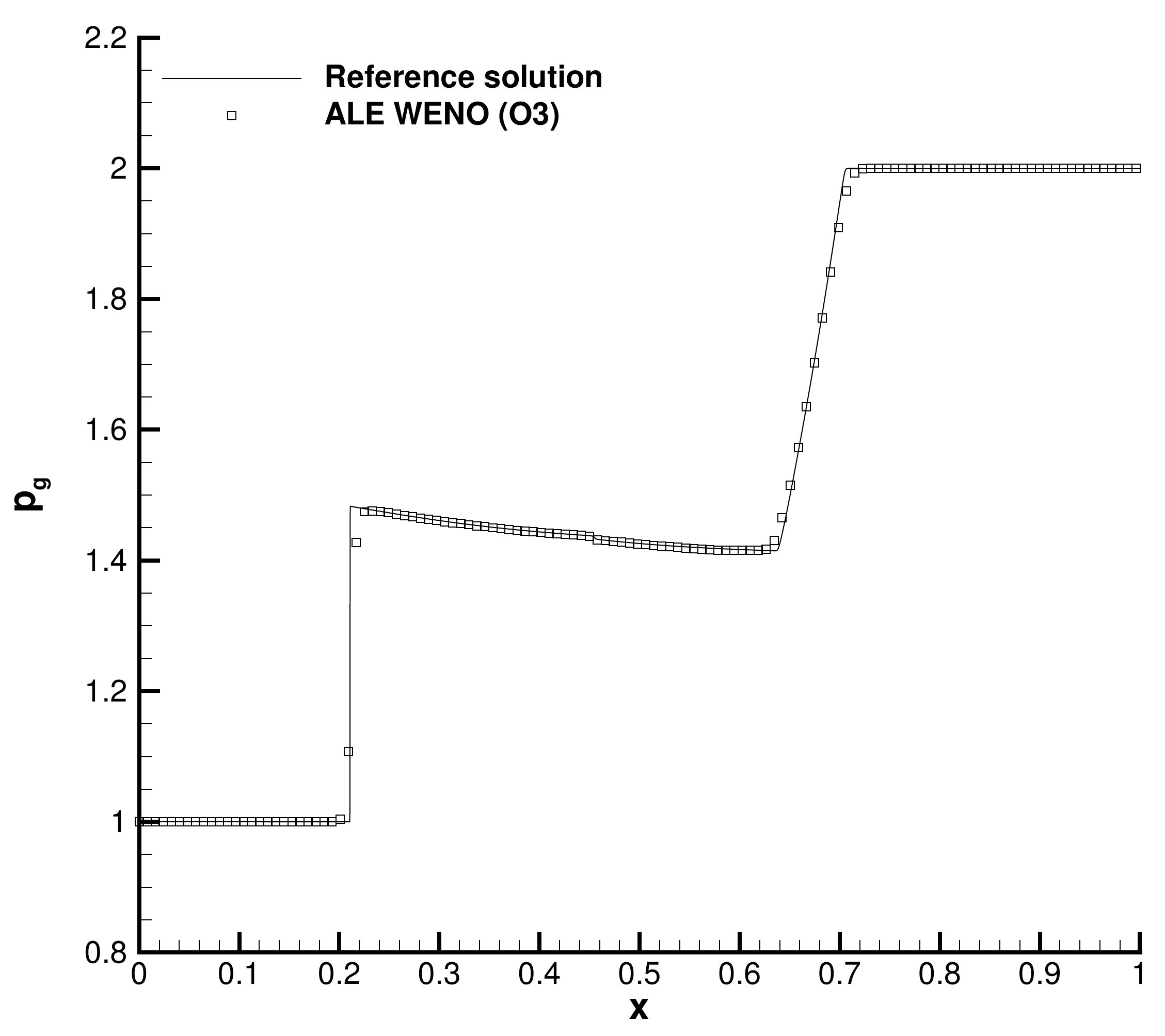}       
\end{tabular}
\caption{Results obtained for the first cylindrical explosion problem EP1 at $t=0.15$ and comparison with the reference solution.} 
\label{fig.bn.ep1}
\end{center}
\end{figure}

\begin{figure}[!htbp]
\begin{center}
\begin{tabular}{cc} 
\includegraphics[width=0.4\textwidth]{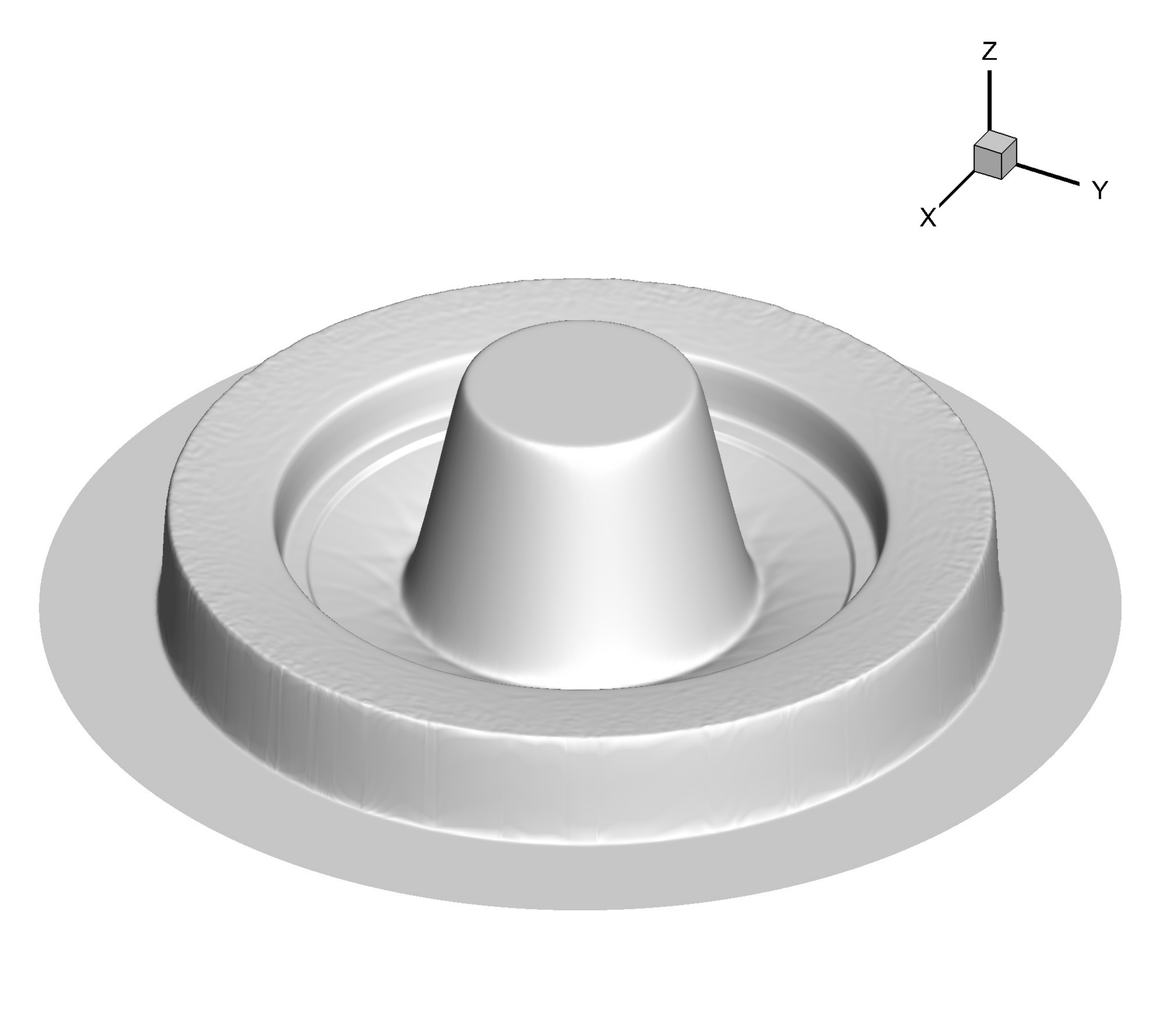}       & 
\includegraphics[width=0.4\textwidth]{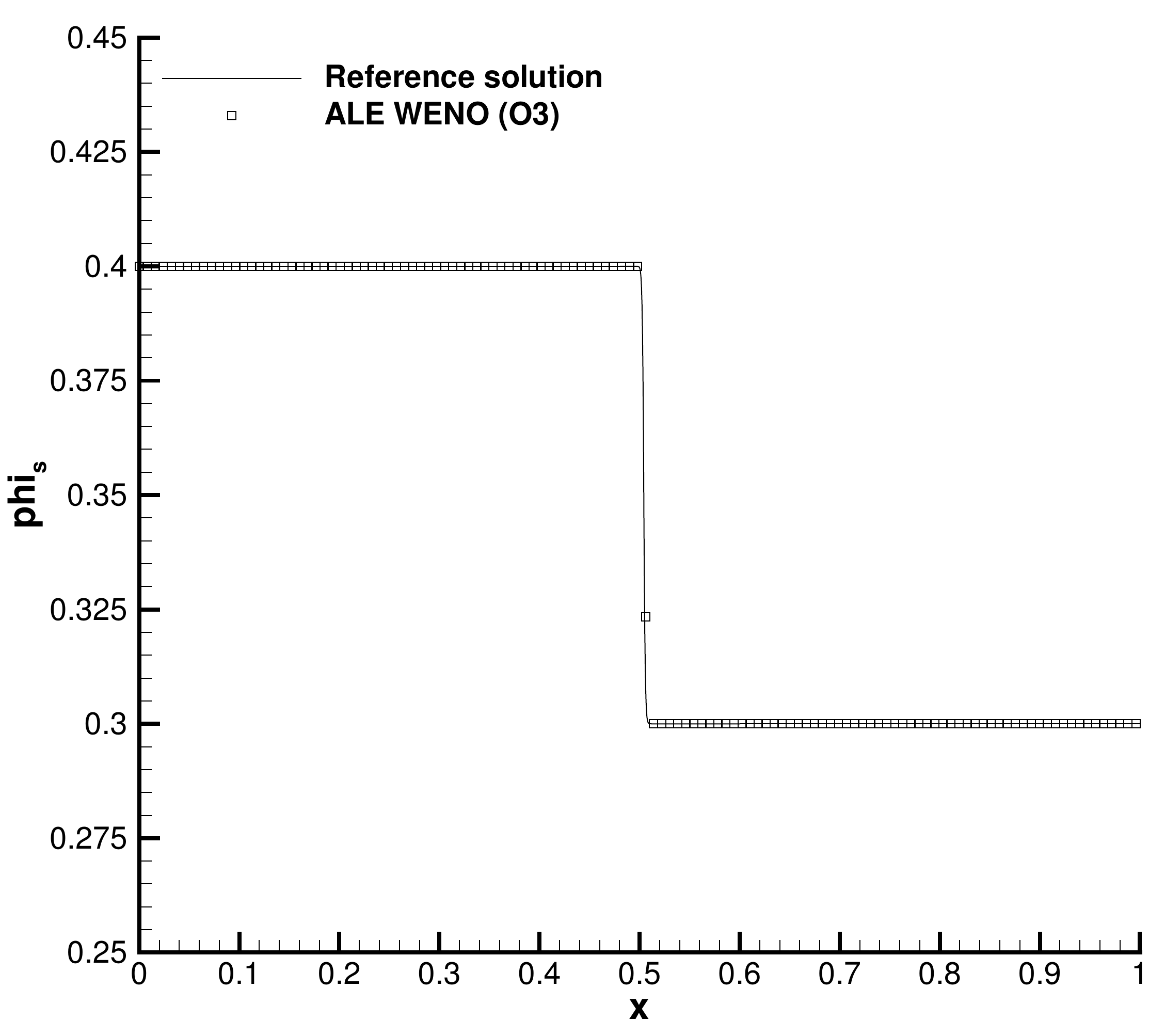}    \\ 
\includegraphics[width=0.4\textwidth]{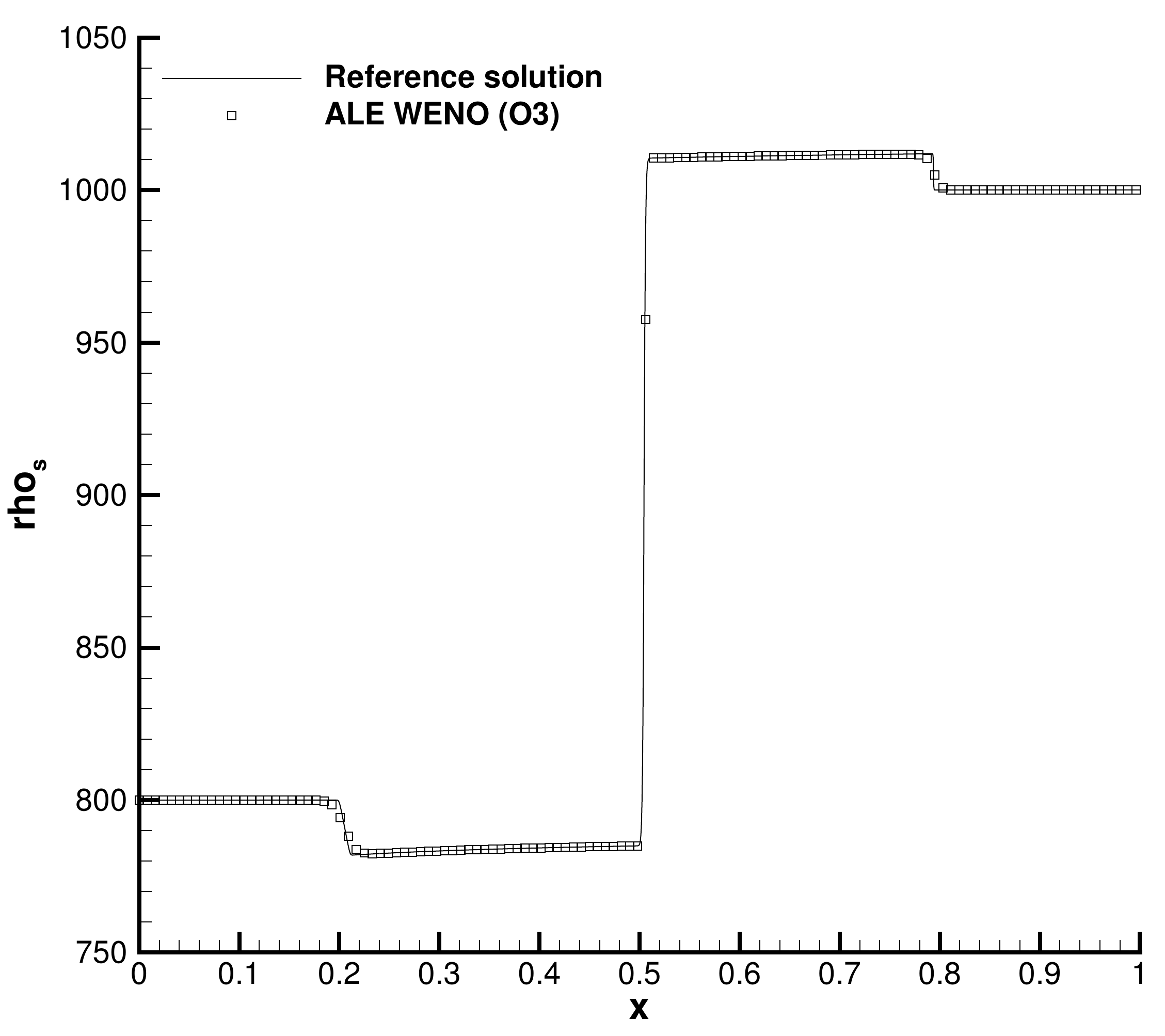}     & 
\includegraphics[width=0.4\textwidth]{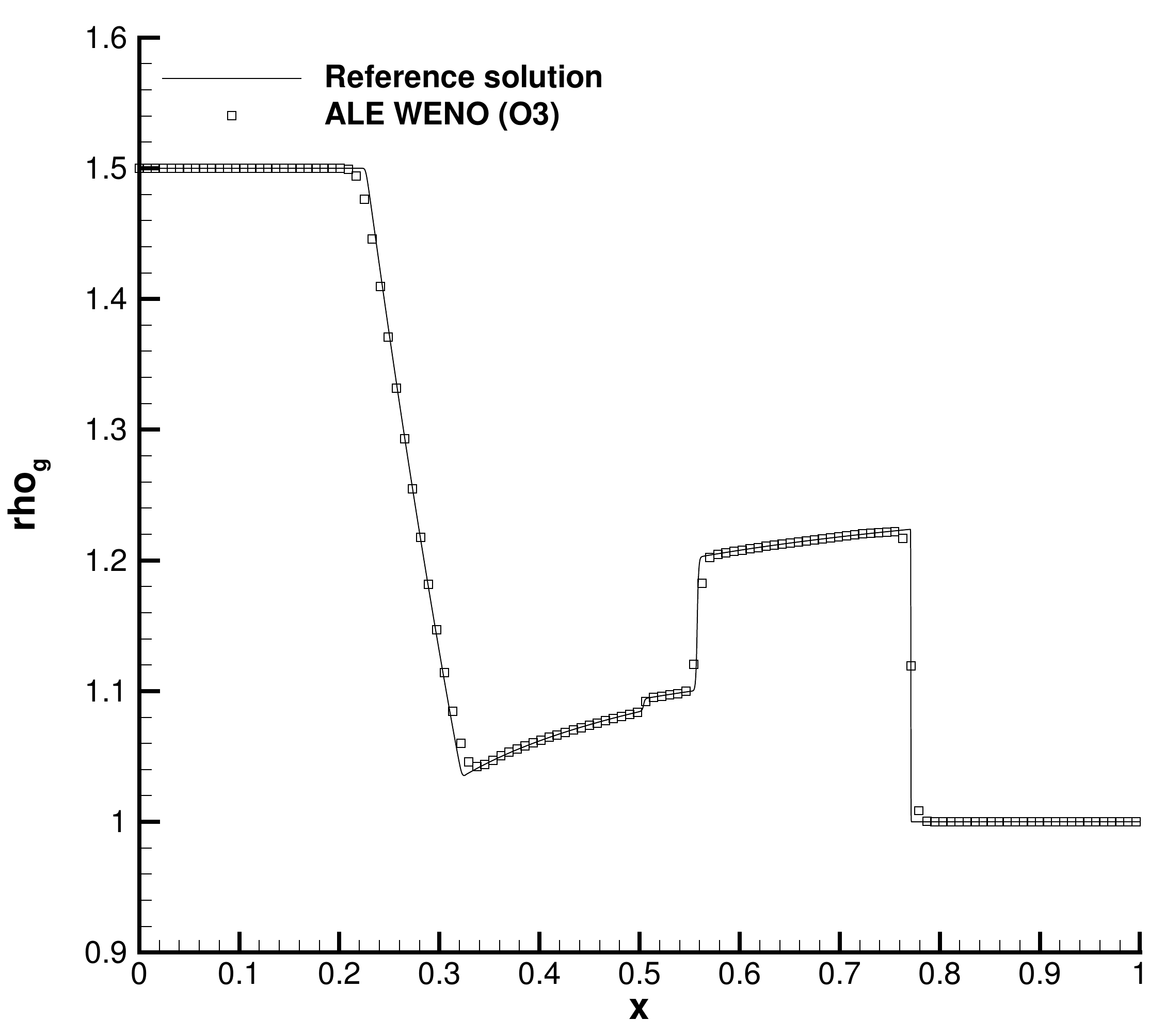}    \\ 
\includegraphics[width=0.4\textwidth]{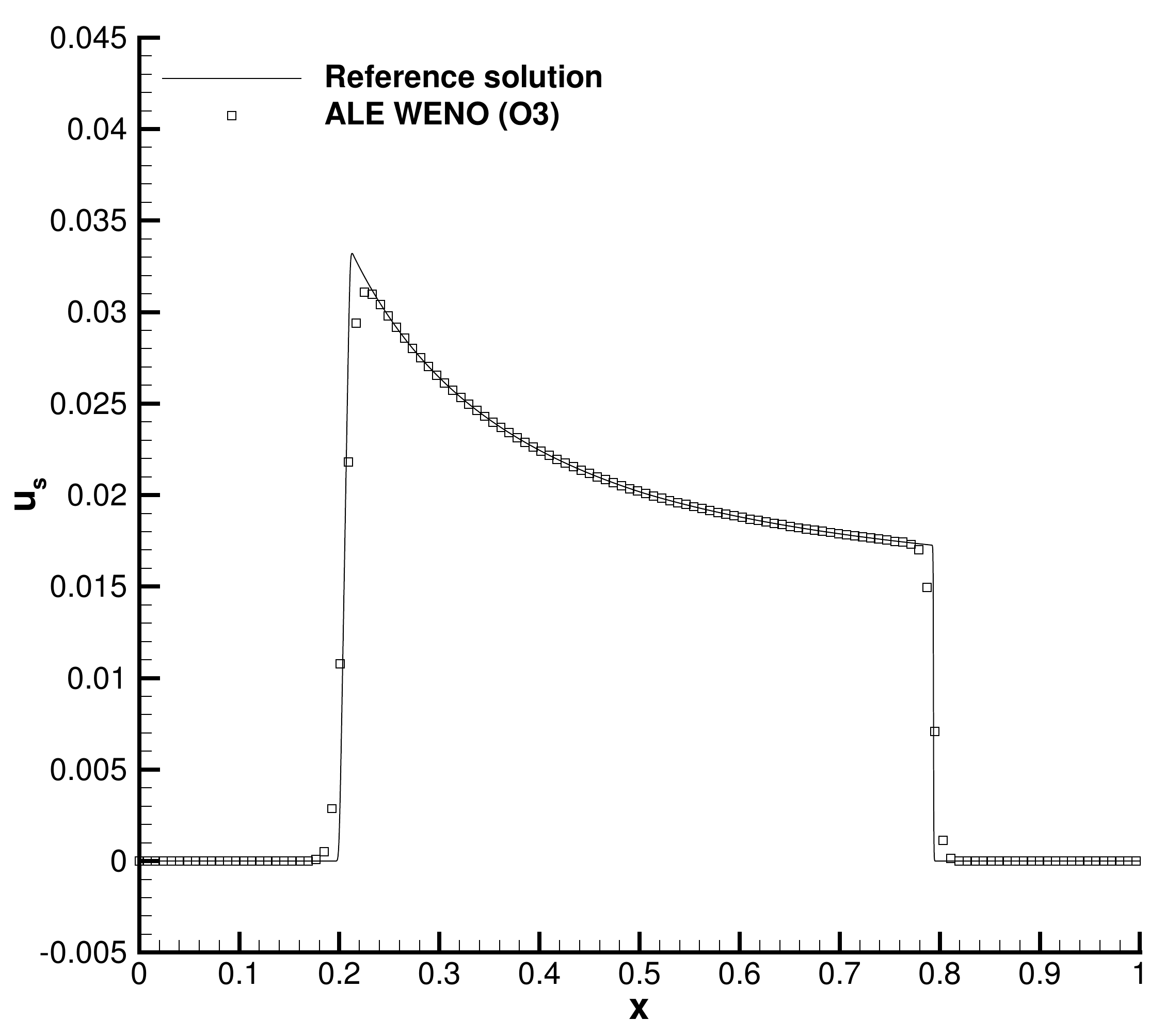}       & 
\includegraphics[width=0.4\textwidth]{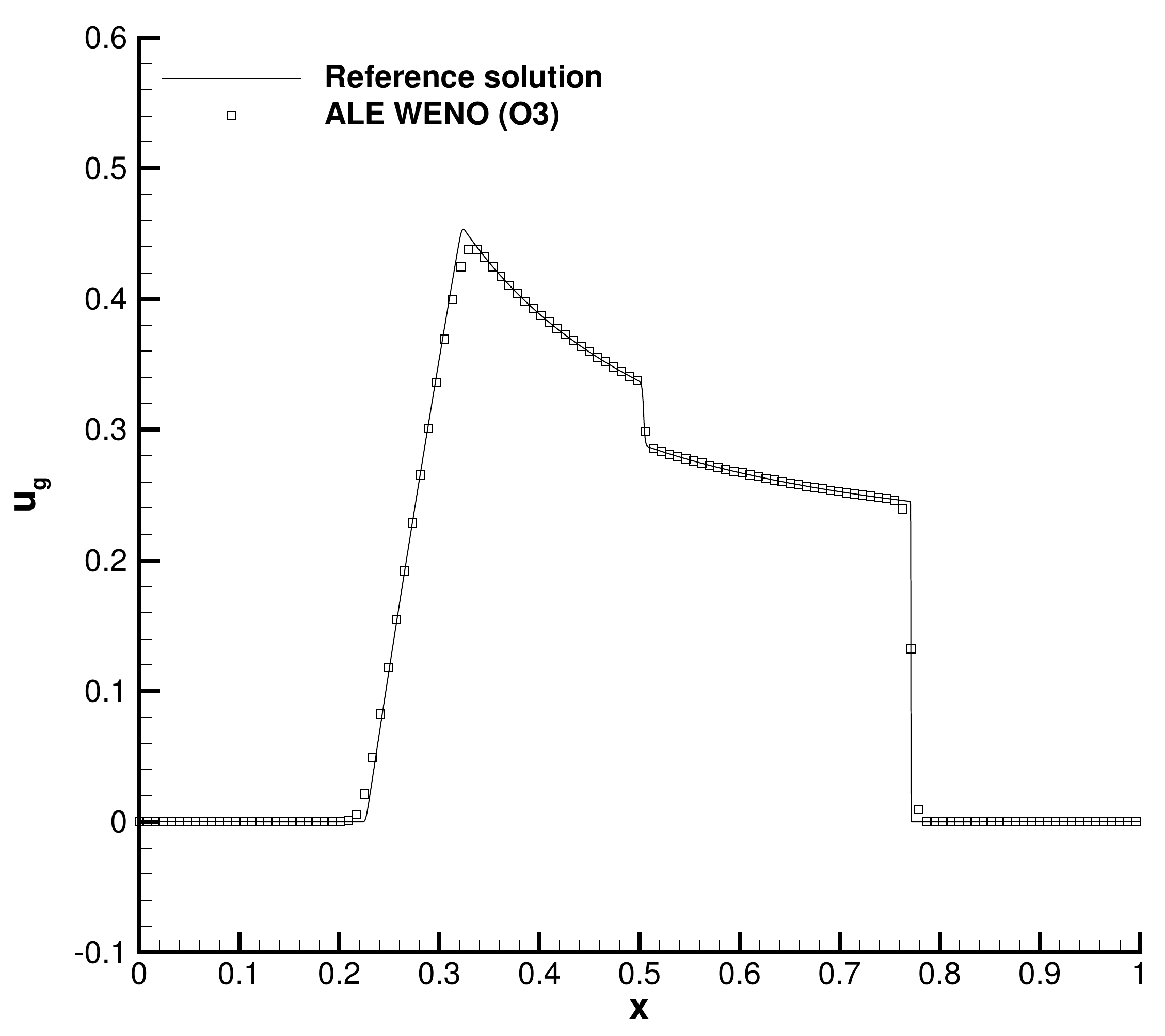}      \\ 
\includegraphics[width=0.4\textwidth]{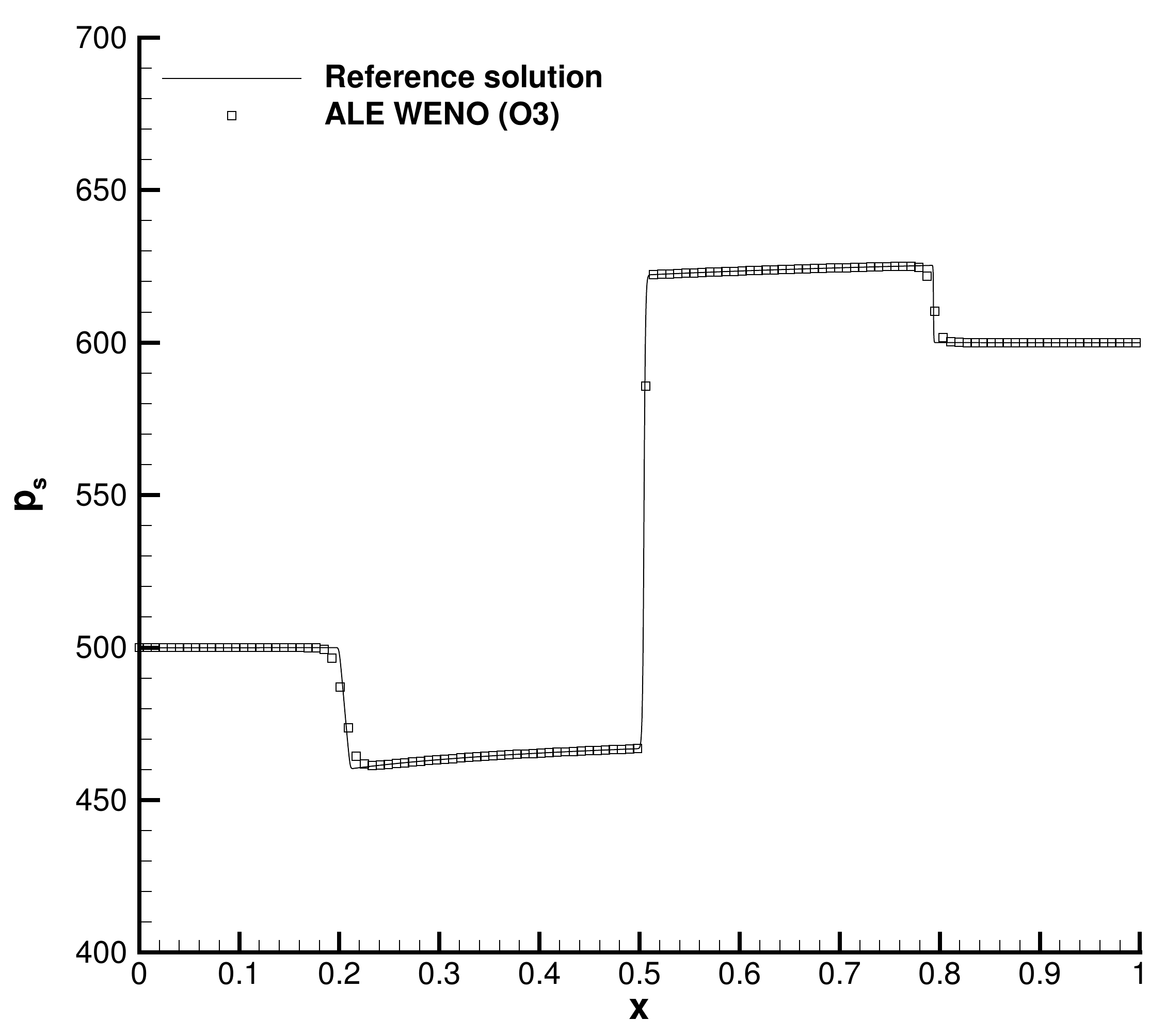}       & 
\includegraphics[width=0.4\textwidth]{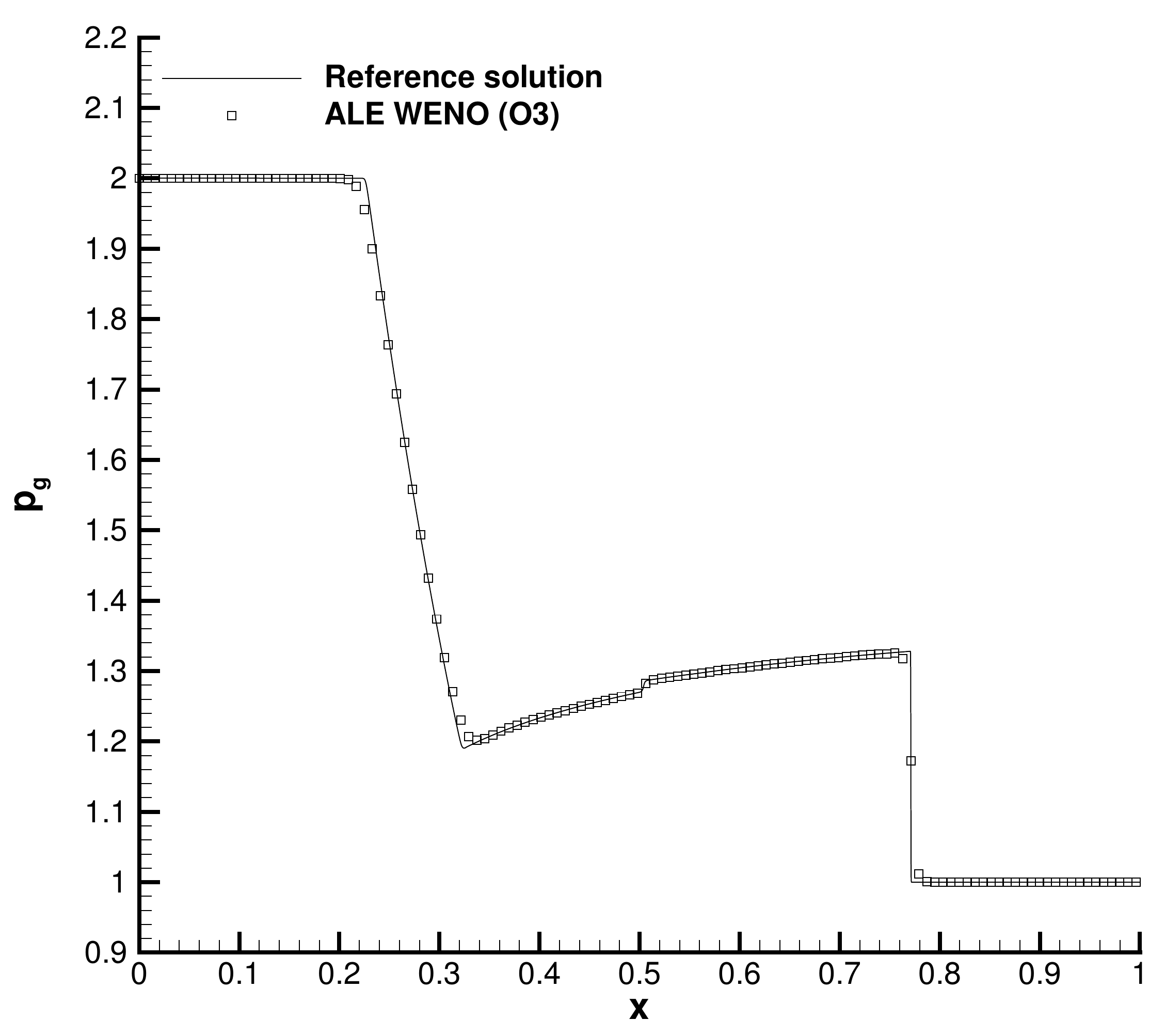}     
\end{tabular}
\caption{Results obtained for the first cylindrical explosion problem EP2 at $t=0.15$ and comparison with the reference solution.} 
\label{fig.bn.ep2}
\end{center}
\end{figure}

\begin{figure}[!htbp]
\begin{center}
\begin{tabular}{cc} 
\includegraphics[width=0.4\textwidth]{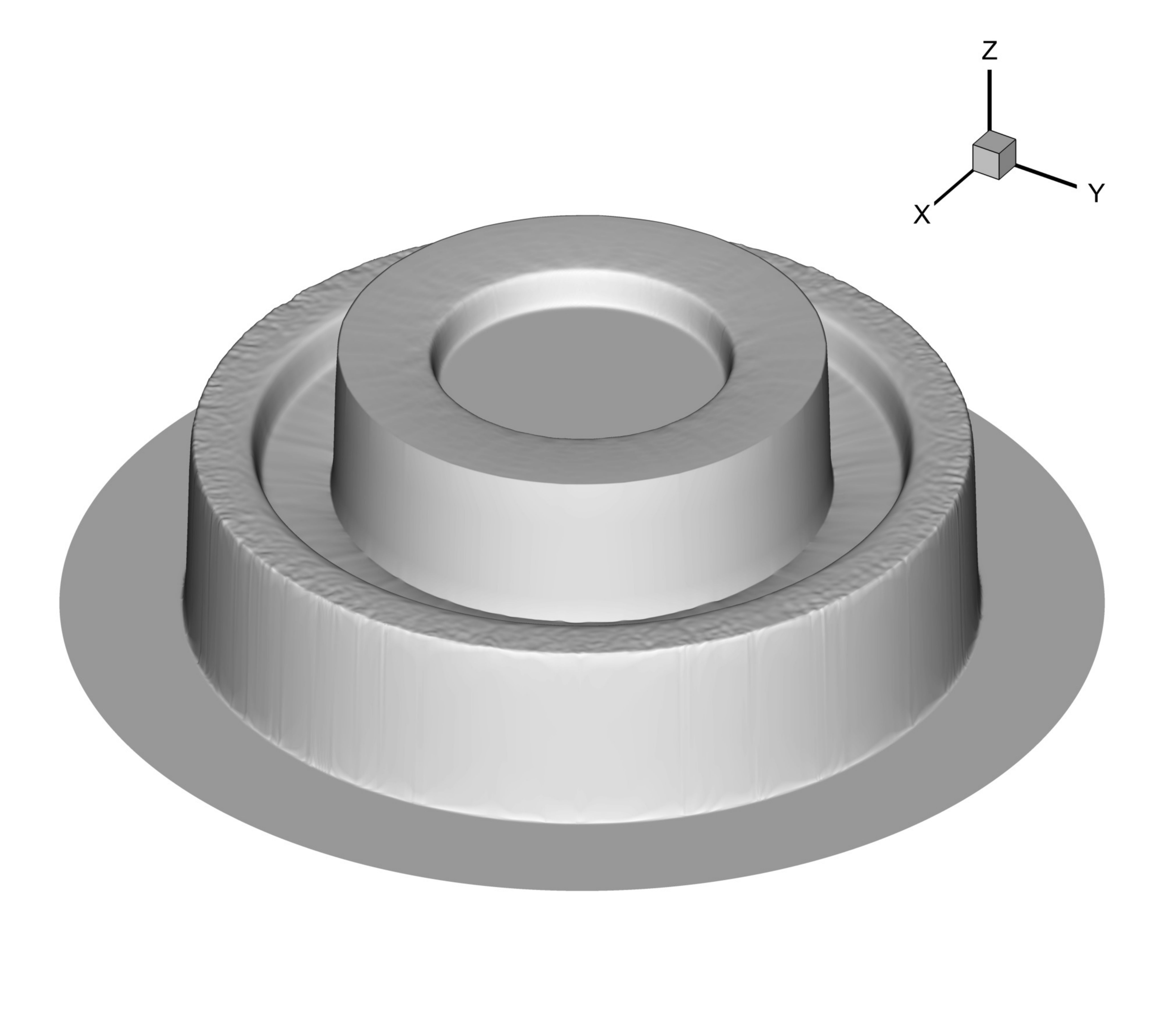}       & 
\includegraphics[width=0.4\textwidth]{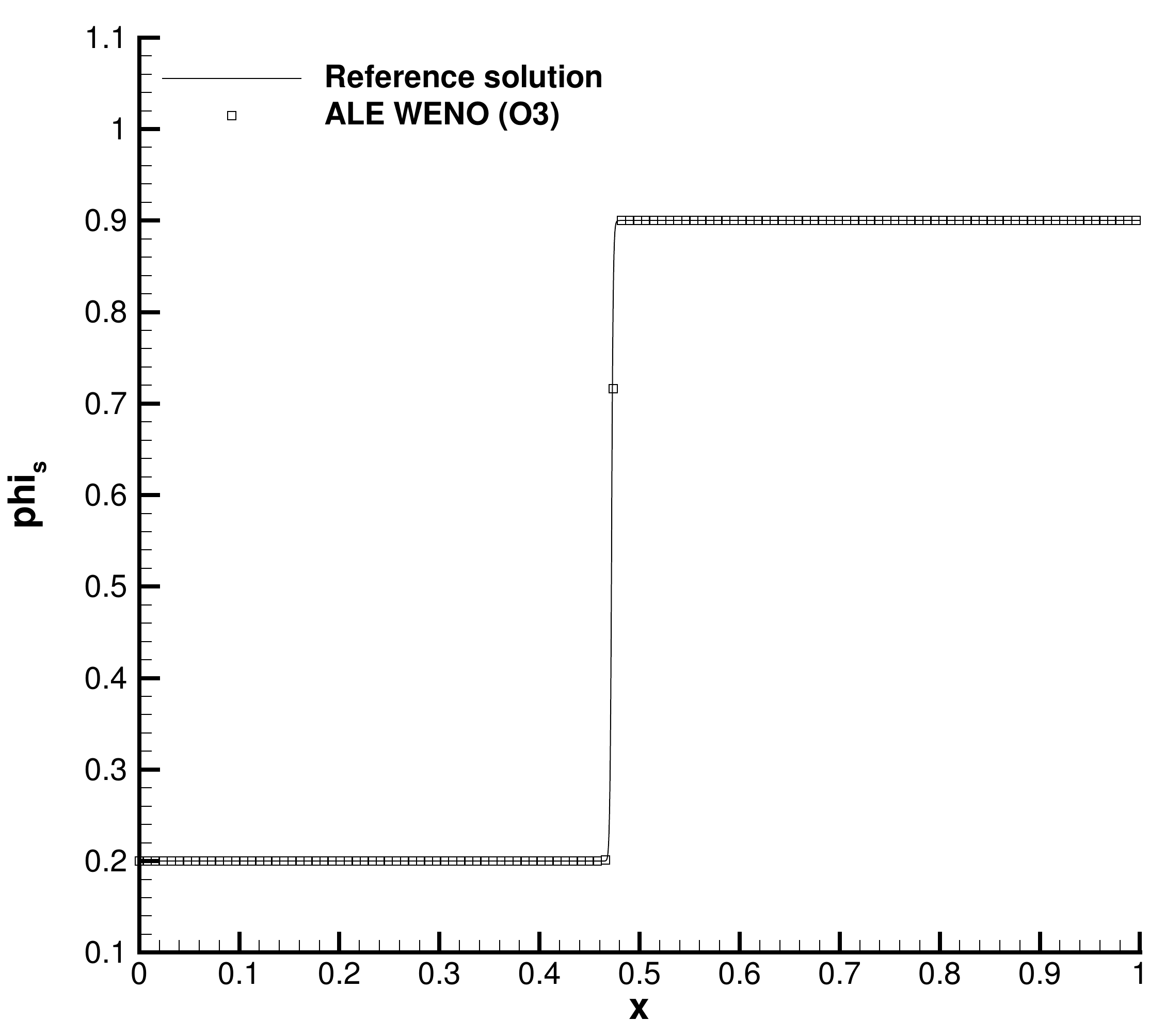}    \\ 
\includegraphics[width=0.4\textwidth]{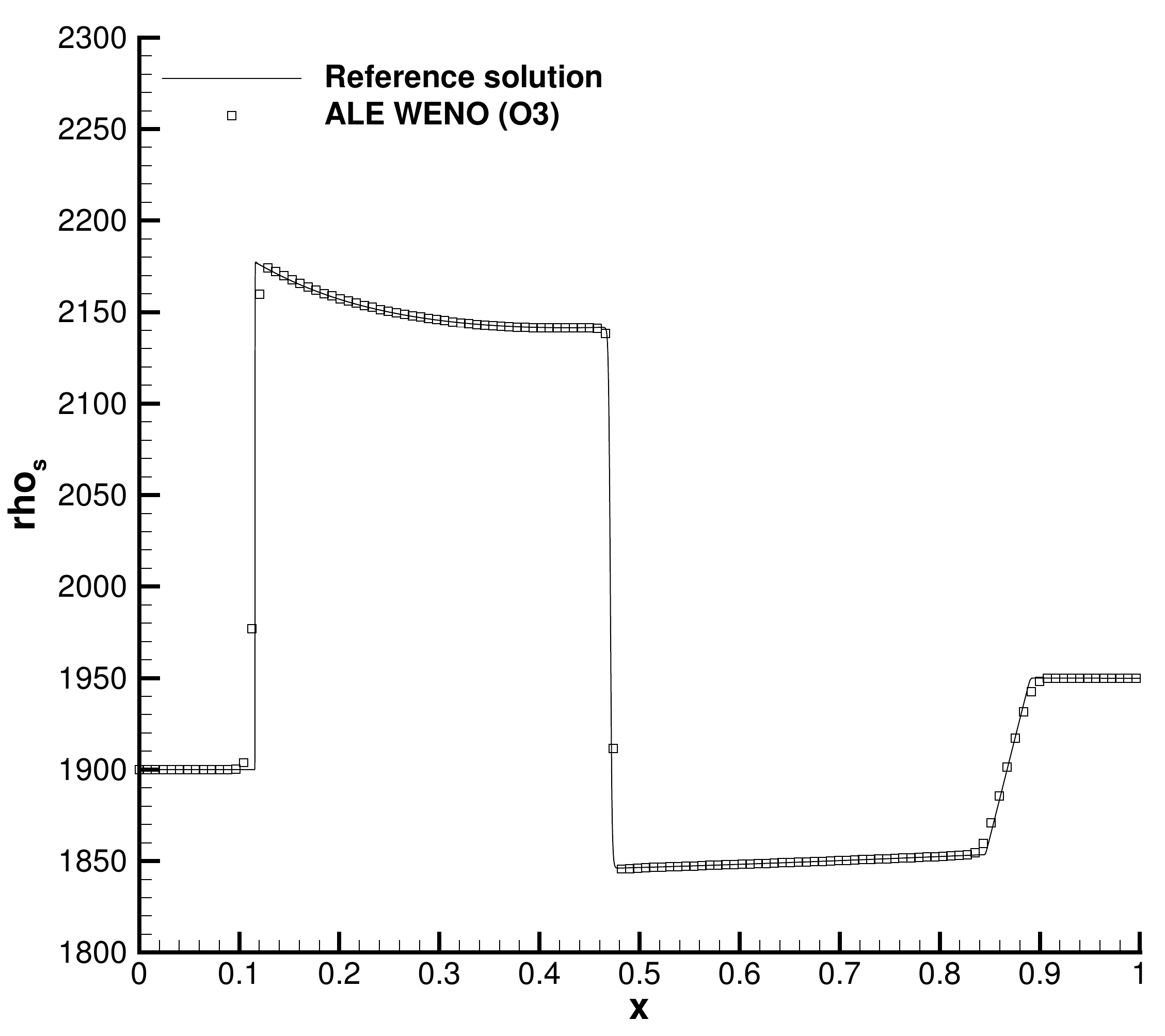}     & 
\includegraphics[width=0.4\textwidth]{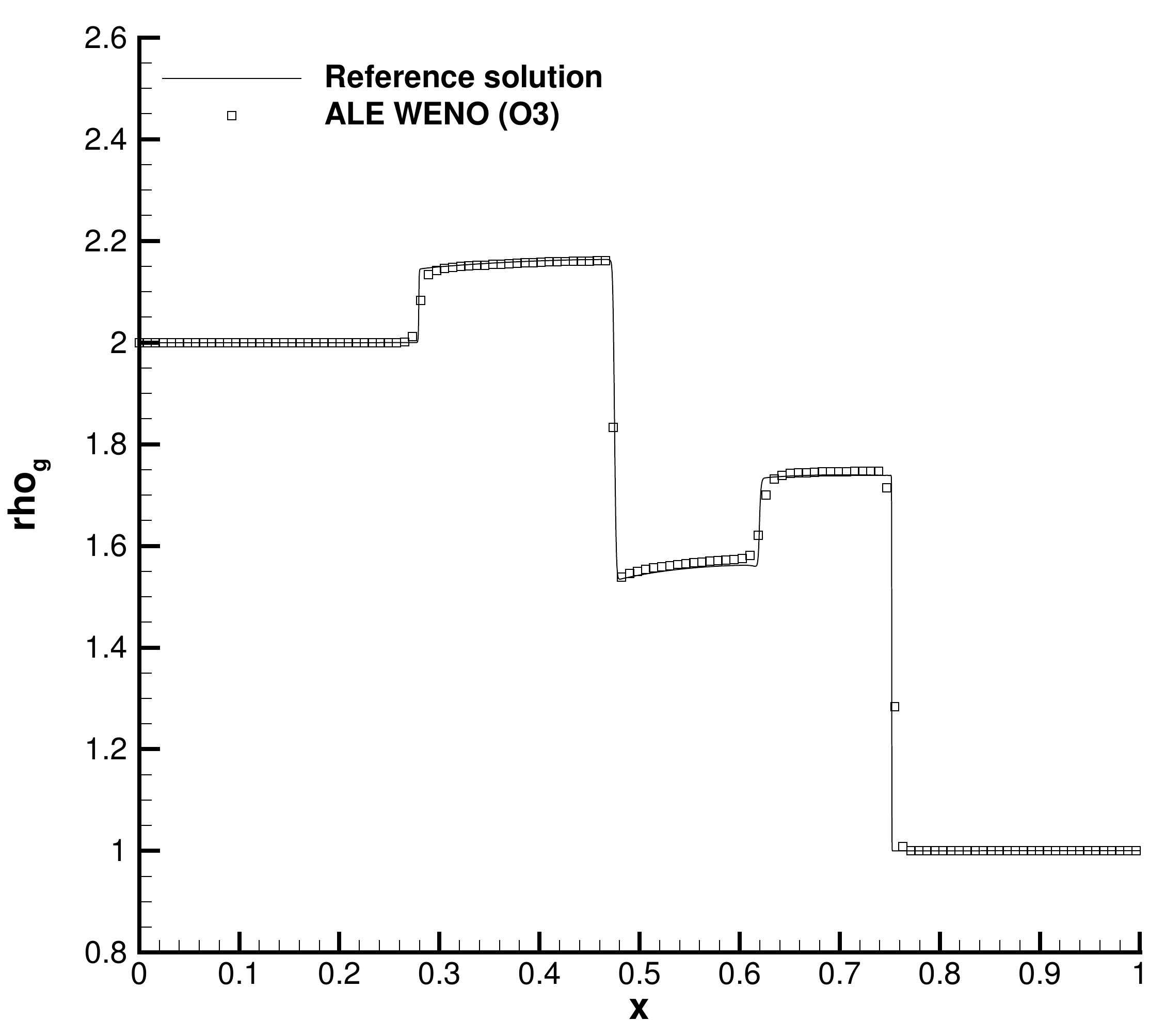}    \\ 
\includegraphics[width=0.4\textwidth]{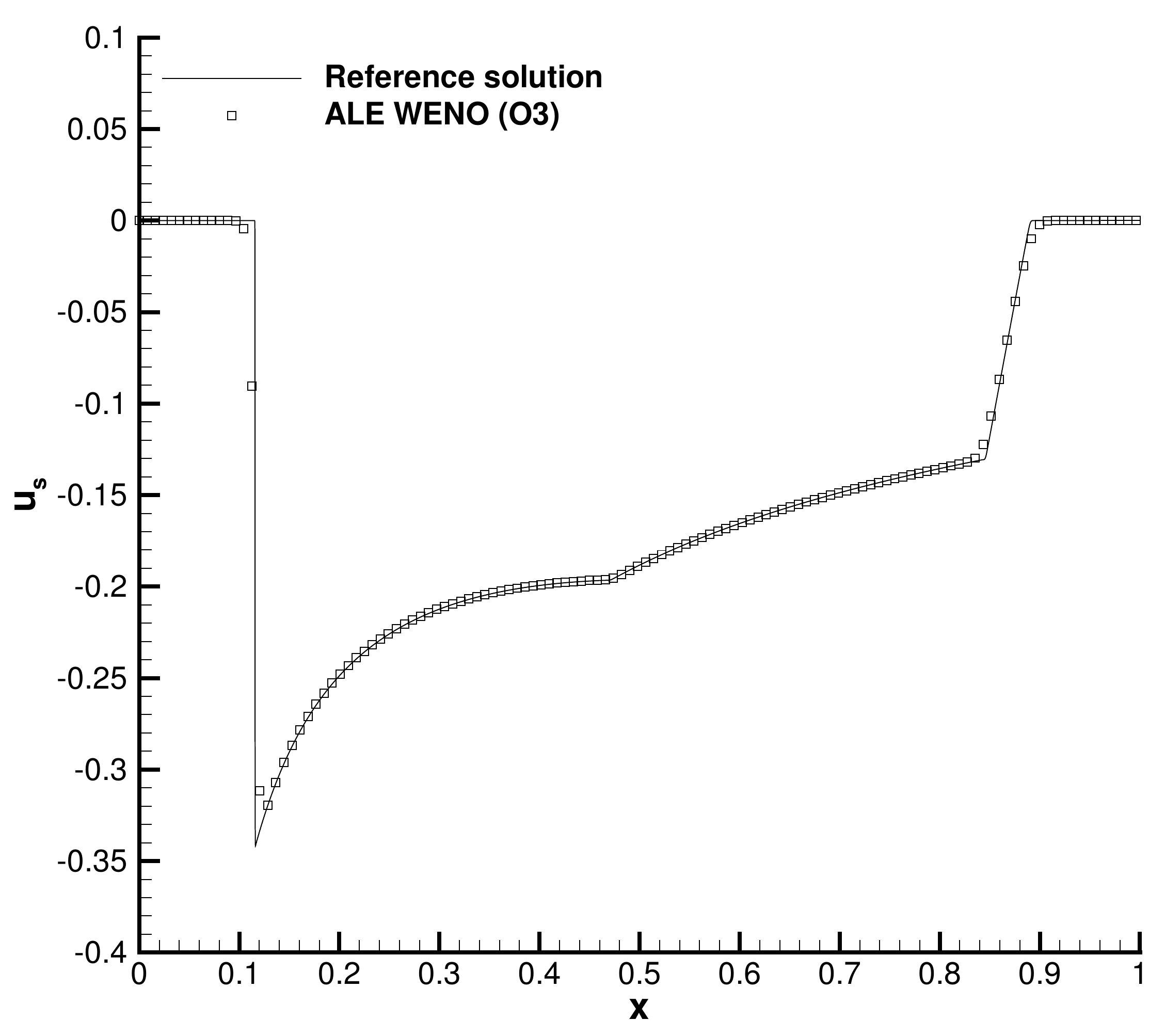}       & 
\includegraphics[width=0.4\textwidth]{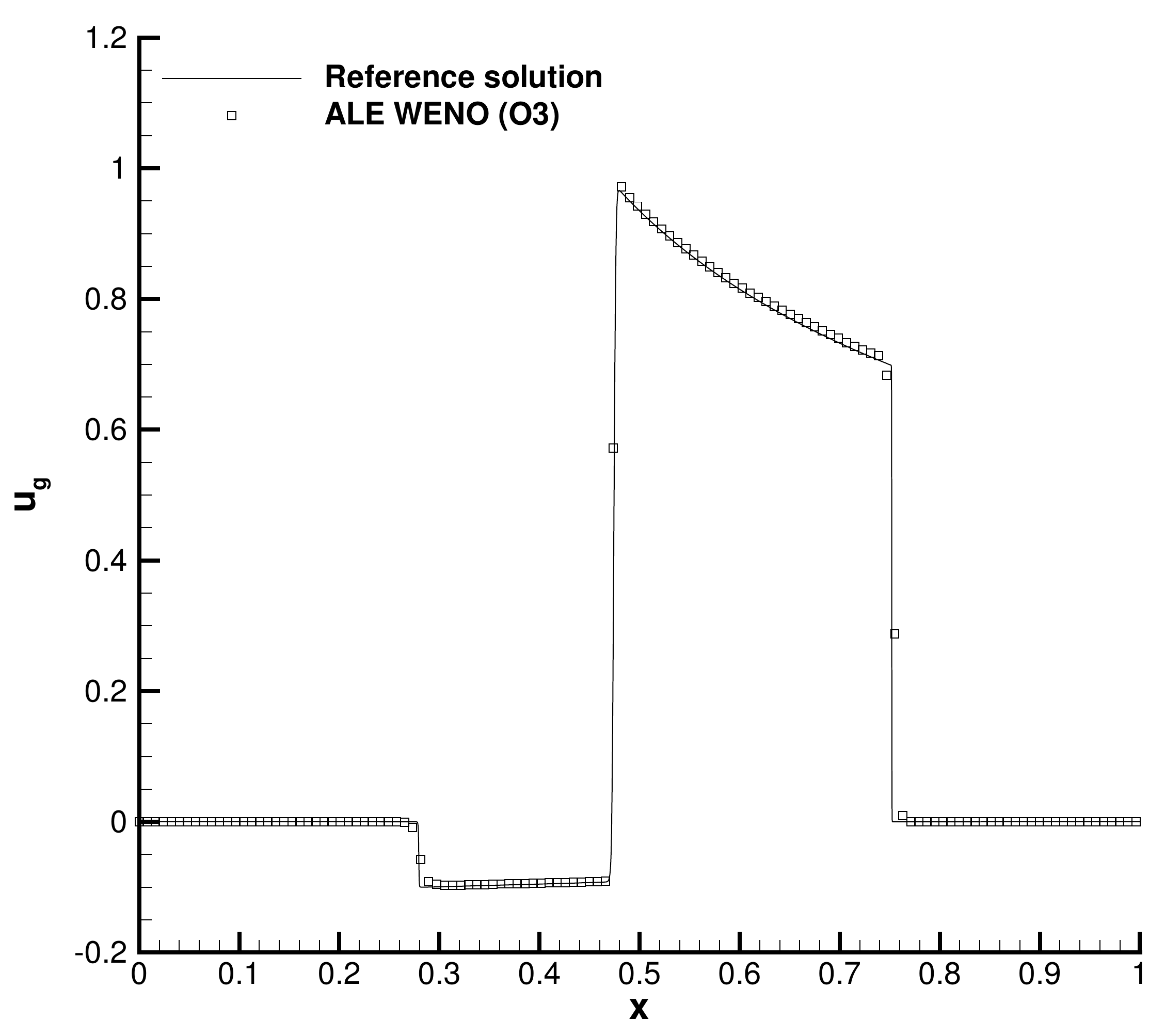}      \\ 
\includegraphics[width=0.4\textwidth]{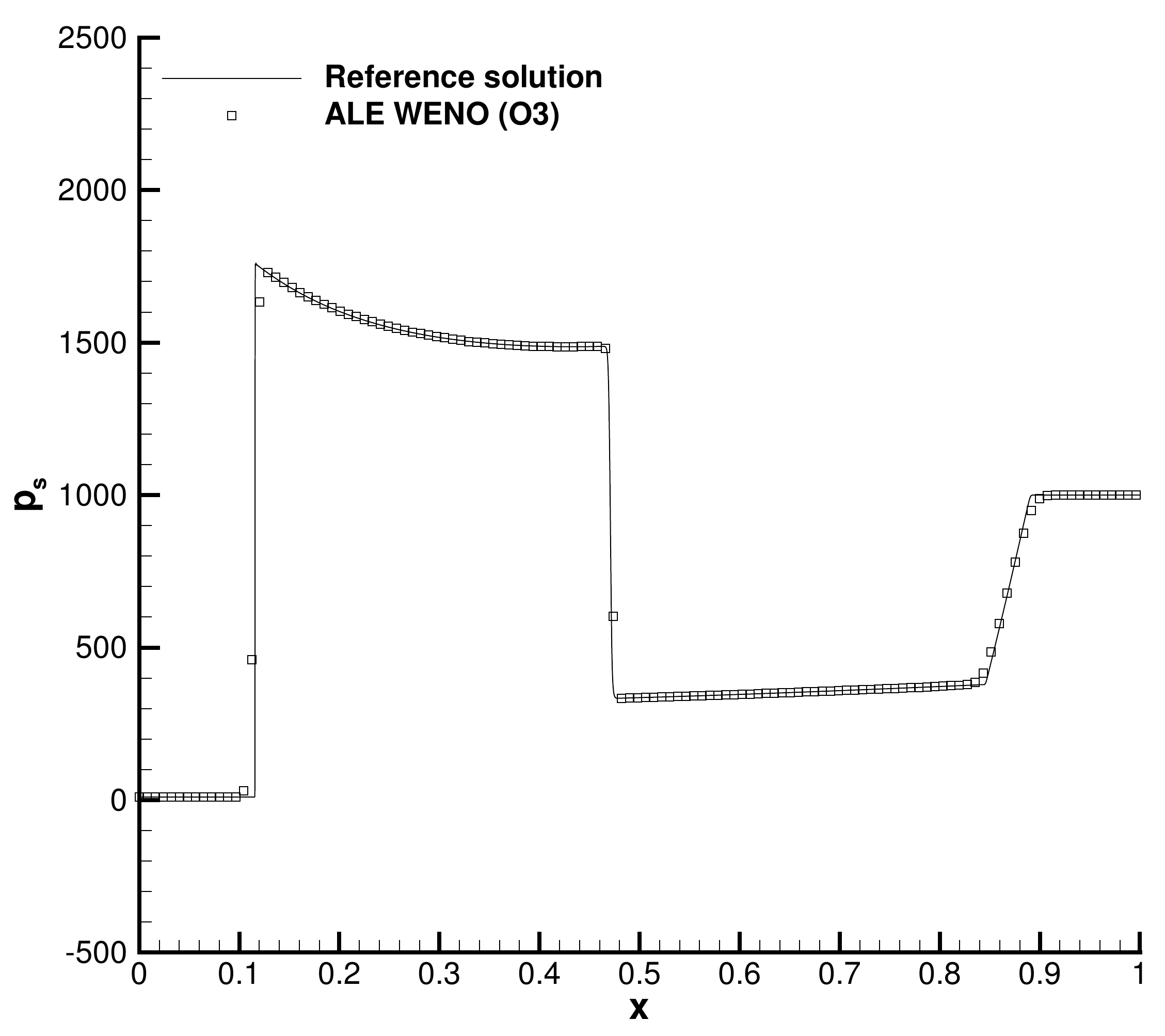}       & 
\includegraphics[width=0.4\textwidth]{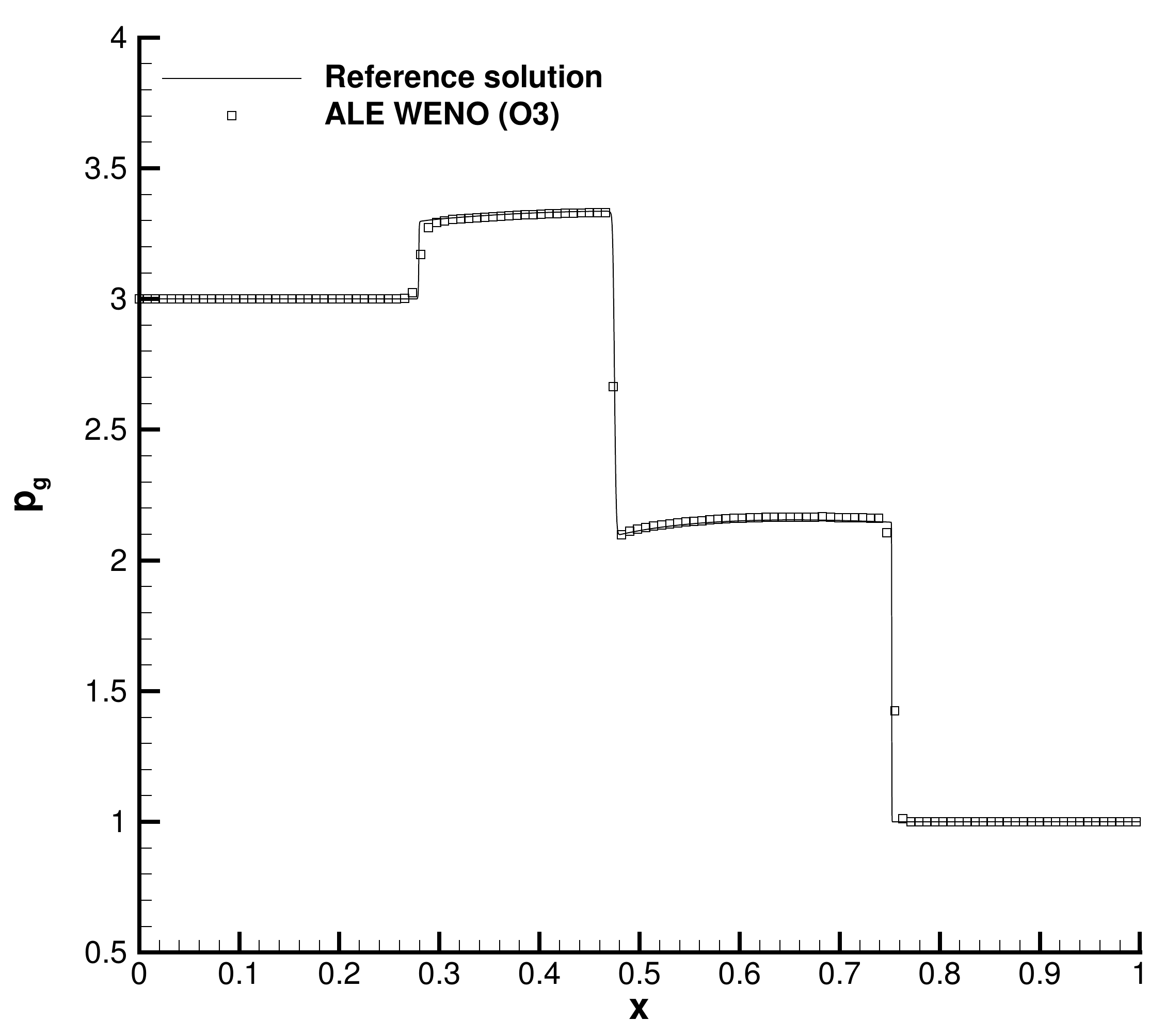}     
\end{tabular}
\caption{Results obtained for the first cylindrical explosion problem EP3 at $t=0.15$ and comparison with the reference solution.} 
\label{fig.bn.ep3}
\end{center}
\end{figure}

Here we solve the compressible Baer--Nunziato equations in two space
dimensions on a circular domain $\Omega(t)$ with initial radius 
$R=1.0$. The initial condition is in all cases given by 
\begin{equation}
  \Q(\x,0) = \left\{ \begin{array}{ll} 
                            \Q_i, & \qquad \textnormal{ if } \quad | \x |  < 0.5, \\
                            \Q_o, & \qquad \textnormal{ else.}          
                      \end{array} \right.  
\label{eqn.ep.ic}                       
\end{equation}
The reference solution is obtained by solving 
an equivalent non-conservative one-dimensional PDE in radial direction with geometric reaction source 
terms, see \cite{toro-book} for the Euler equations and \cite{USFORCE} for the Baer--Nunziato model for details. 
In our case here the reference solution has been obtained by using an Eulerian second order TVD scheme on a very fine 1D 
mesh consisting of 10,000 cells. The initial conditions for the three explosion problems solved here are taken from 
the Riemann problems solved previously, where the left state is used as inner state and the right state is used as
the outer state in \eqref{eqn.ep.ic}, respectively. In particular, the first explosion problem EP1 uses the initial 
condition of RP1, EP2 corresponds to RP2 and EP3 to RP4, respectively. Also the parameters for the equation of state are
chosen according to Table \ref{tab.rpbn.ic}. The initial mesh spacing is of characteristic size $h=1/250$, leading to a total
number of 431,224 triangular elements used to discretize $\Omega(t)$. 
The numerical results are compared with the 1D reference solution in Figures \ref{fig.bn.ep1} - \ref{fig.bn.ep3}. On the top left 
of each figure a 3D visualization of either the solid or the gas density is shown, in order to verify that the cylindrical symmetry 
is reasonably maintained on the unstructured triangular meshes used here. The other subfigures show a one--dimensional cut through the
reconstructed numerical solution $\w_h$ on 250 \textit{equidistant} sample points along the $x$-axis. In all cases an excellent agreement 
between numerical solution and reference solution is obtained. Note, in particular, the very sharp resolution of the solid contact due 
to the use of a Lagrangian framework. 

\subsection{Two--Dimensional Riemann Problems} 
\label{sec.rp2d}

\begin{figure}[!ht]
\begin{center}
\begin{tabular}{cc} 
\includegraphics[width=0.47\textwidth]{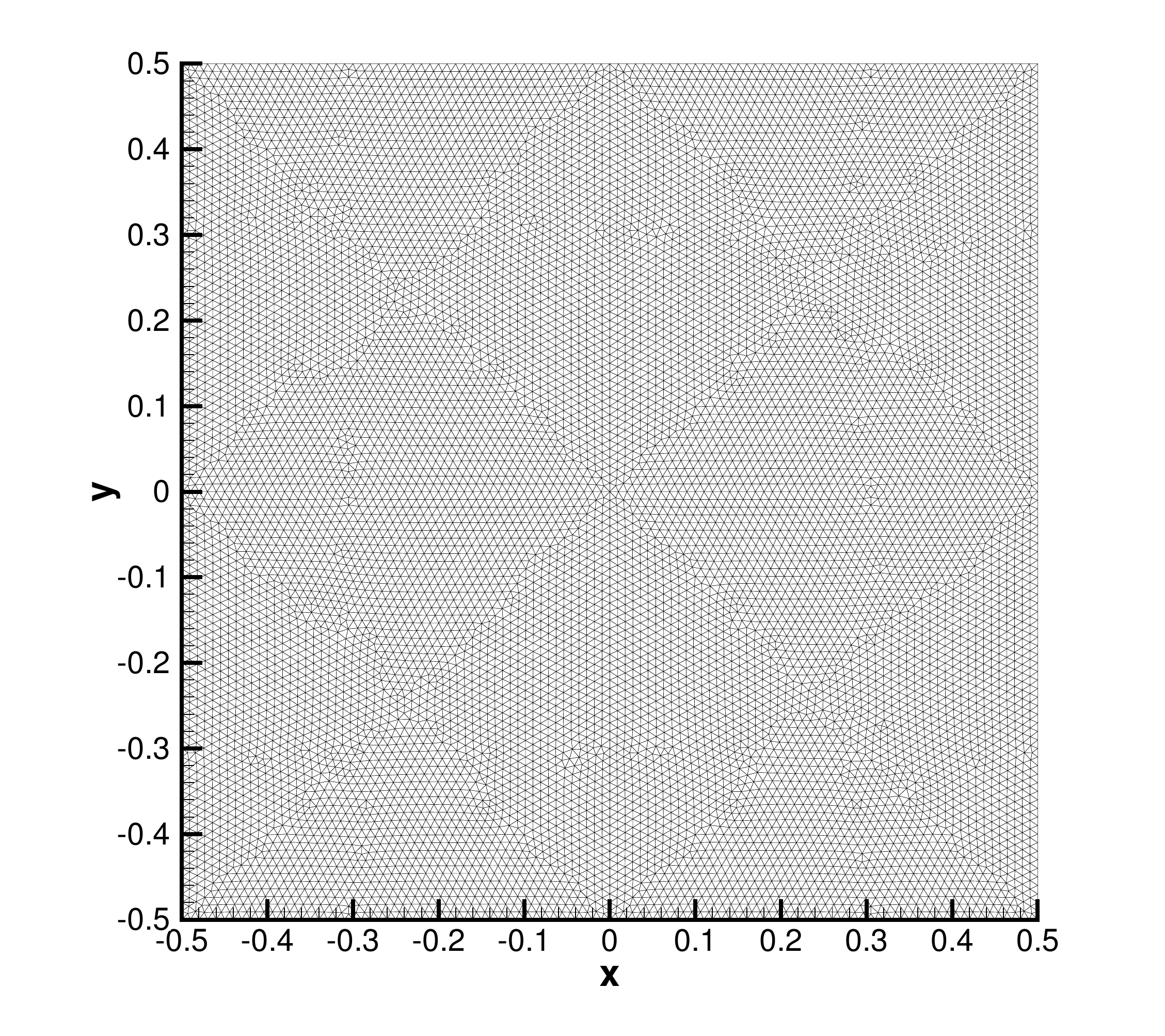}    & 
\includegraphics[width=0.47\textwidth]{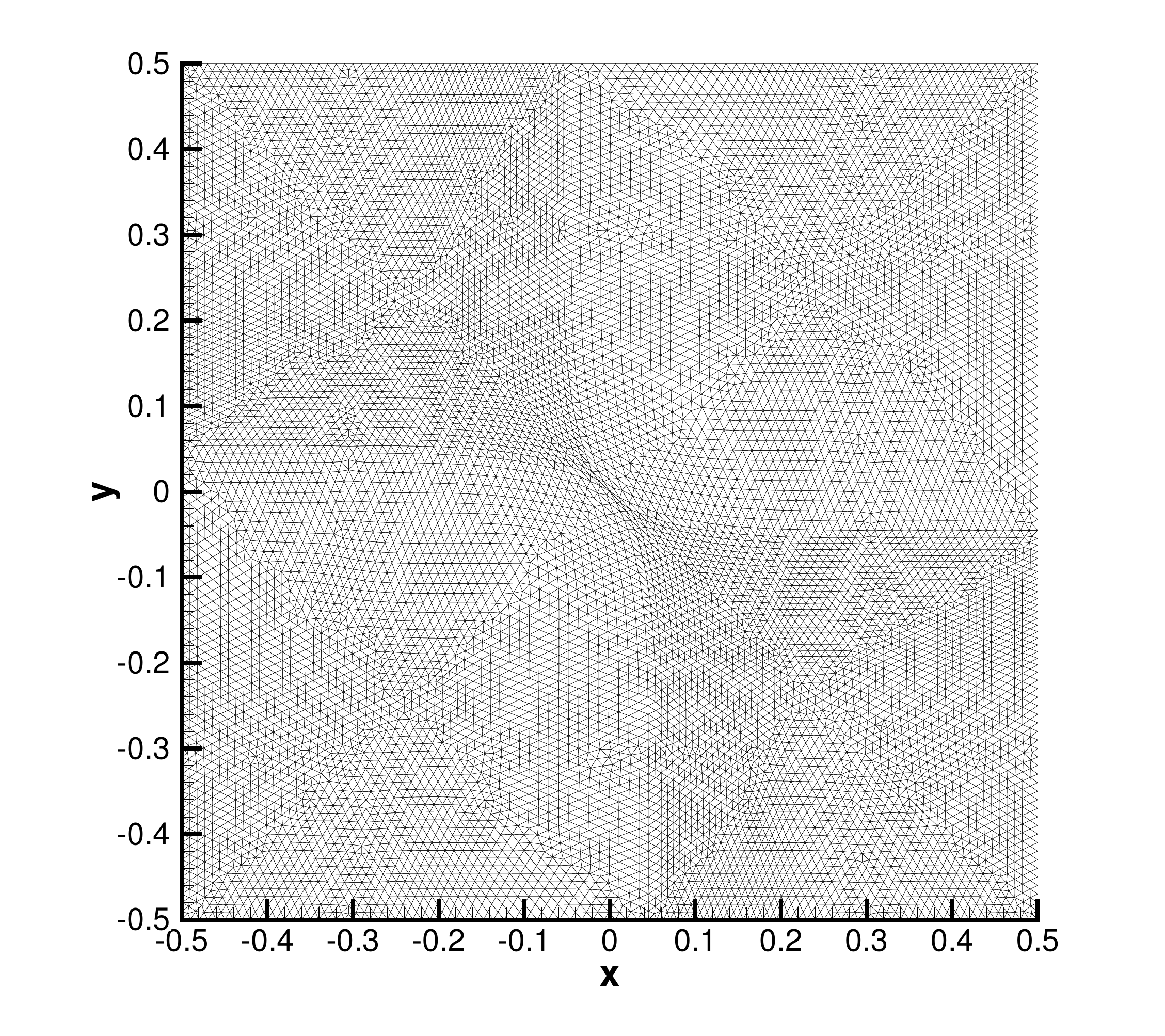}    \\ 
\vspace{-9mm}  
\end{tabular}
\caption{Mesh for configuration C1 at times $t=0$ (left) and $t=0.15$ (right).} 
\label{fig.bn.rp2d1.mesh} 
\end{center}
\end{figure}

\begin{figure}[!ht]
\begin{center}
\begin{tabular}{cc} 
\includegraphics[width=0.45\textwidth]{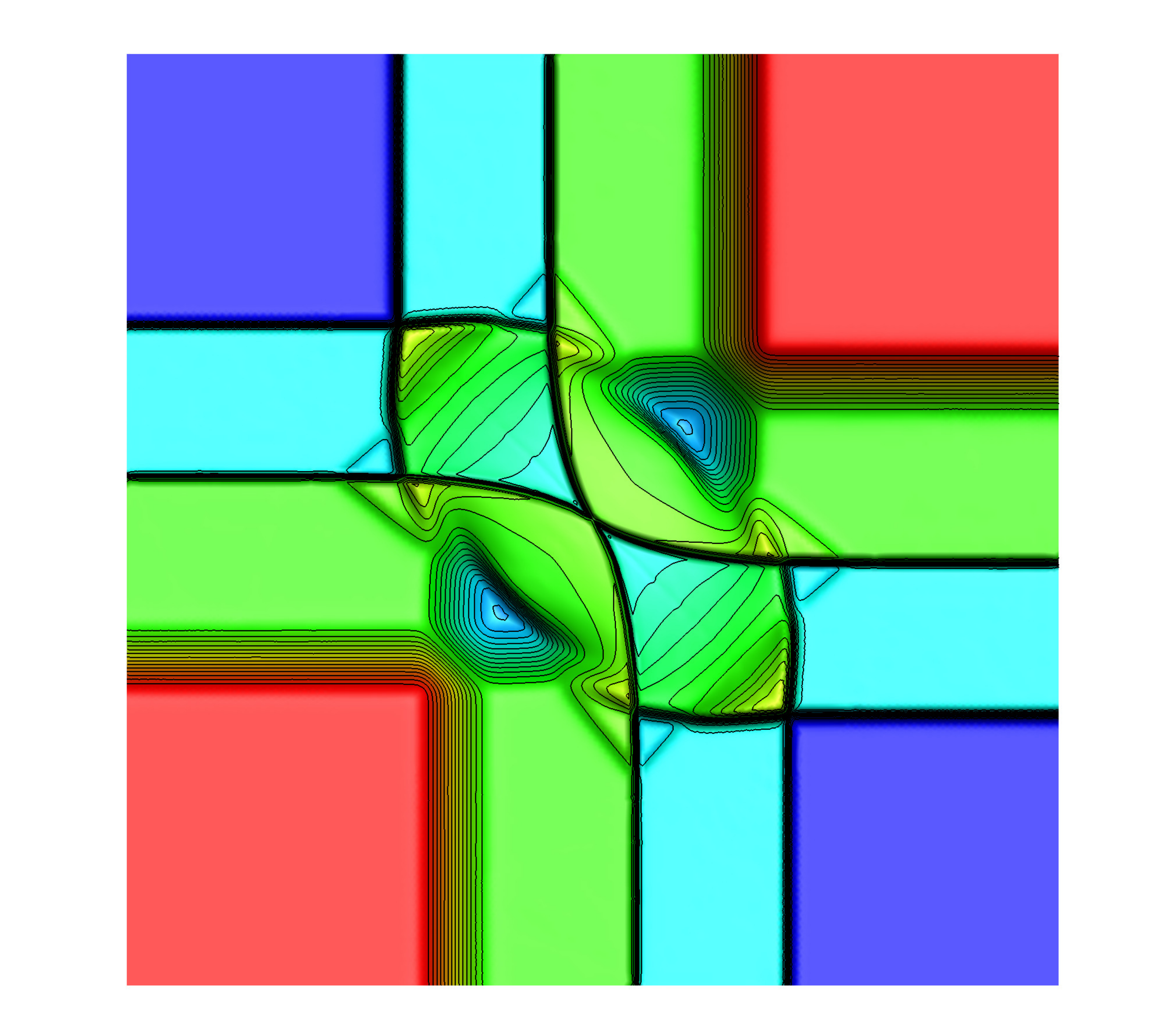}       & 
\includegraphics[width=0.45\textwidth]{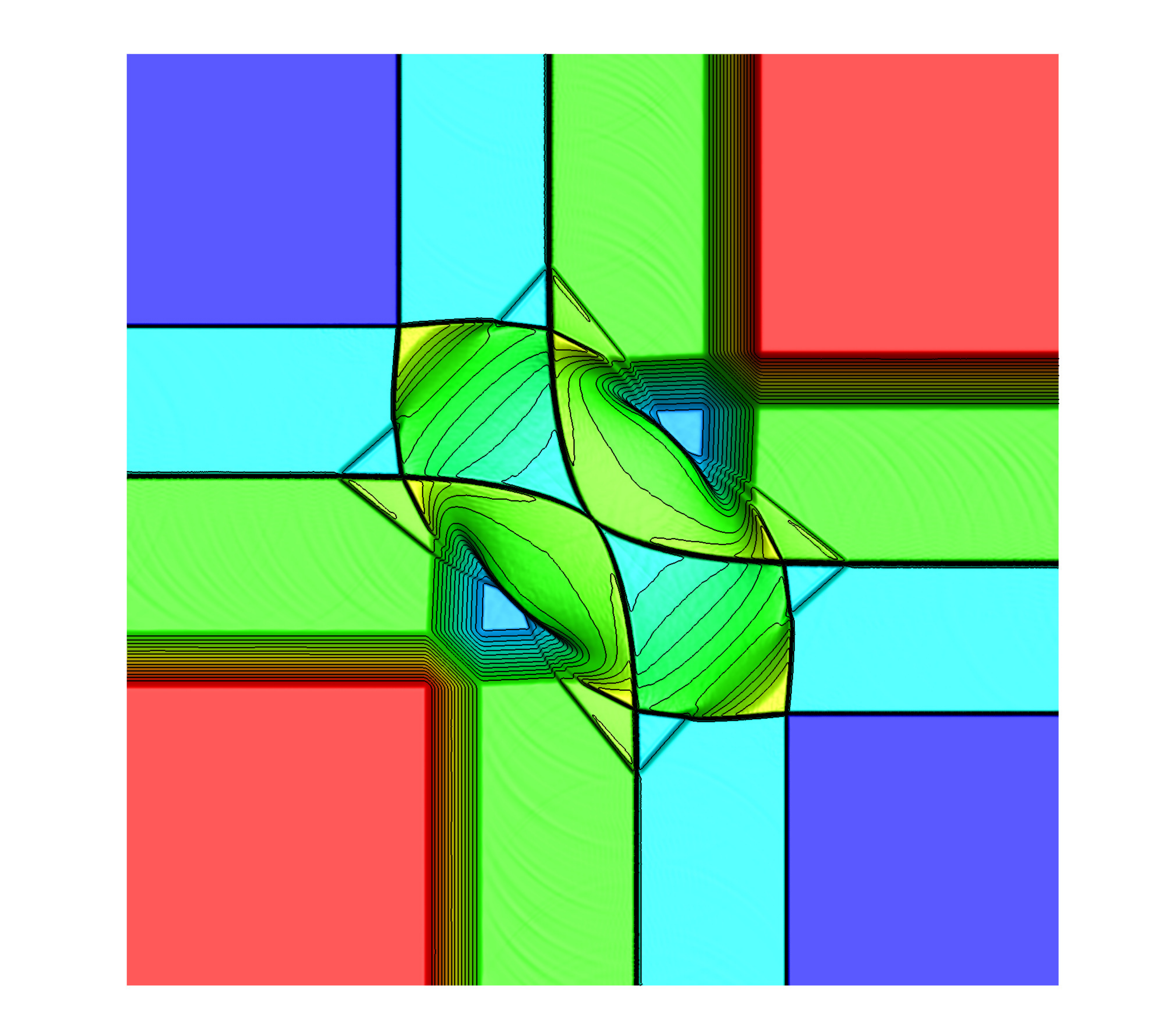}    \\ 
\includegraphics[width=0.45\textwidth]{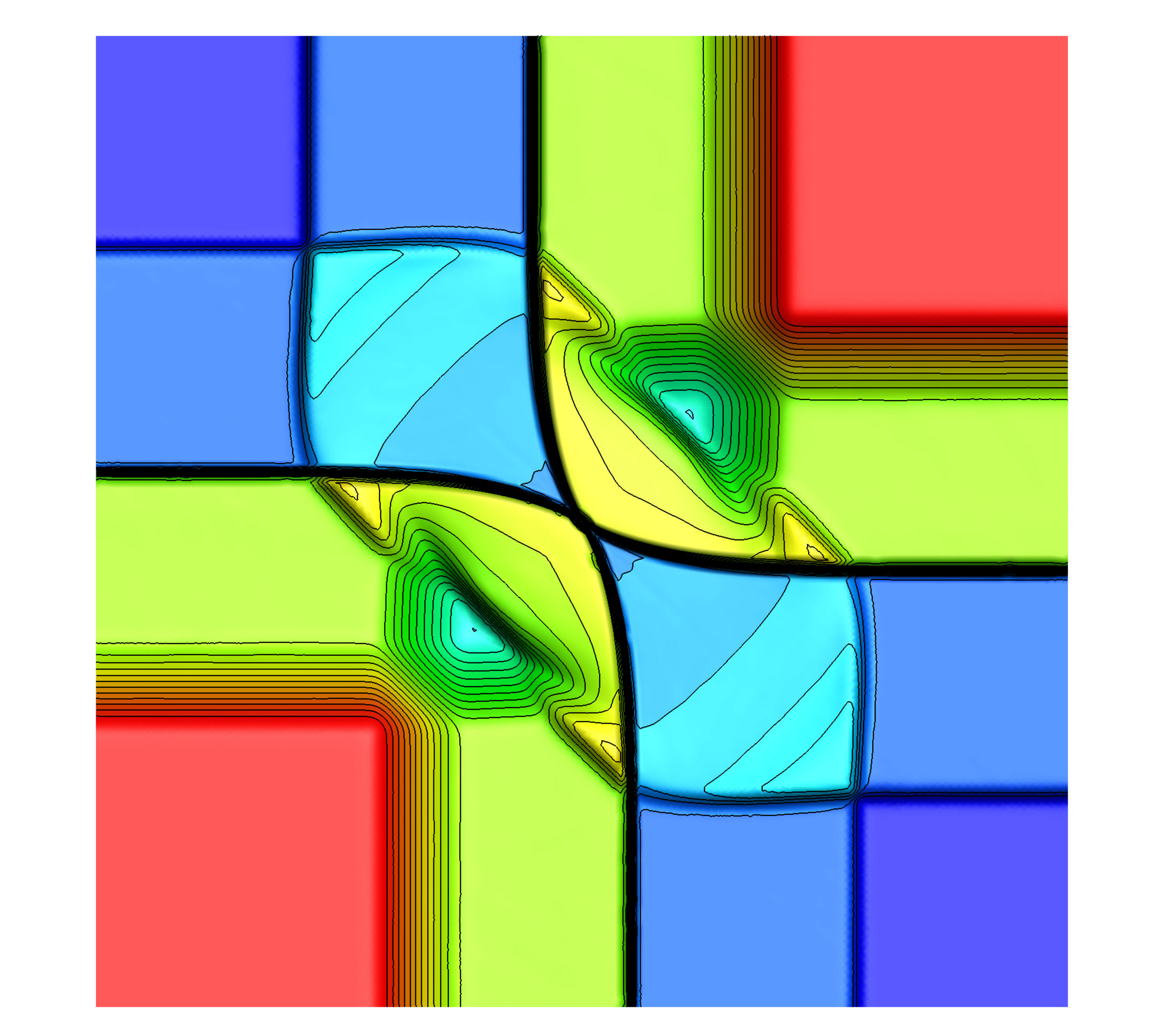}       & 
\includegraphics[width=0.45\textwidth]{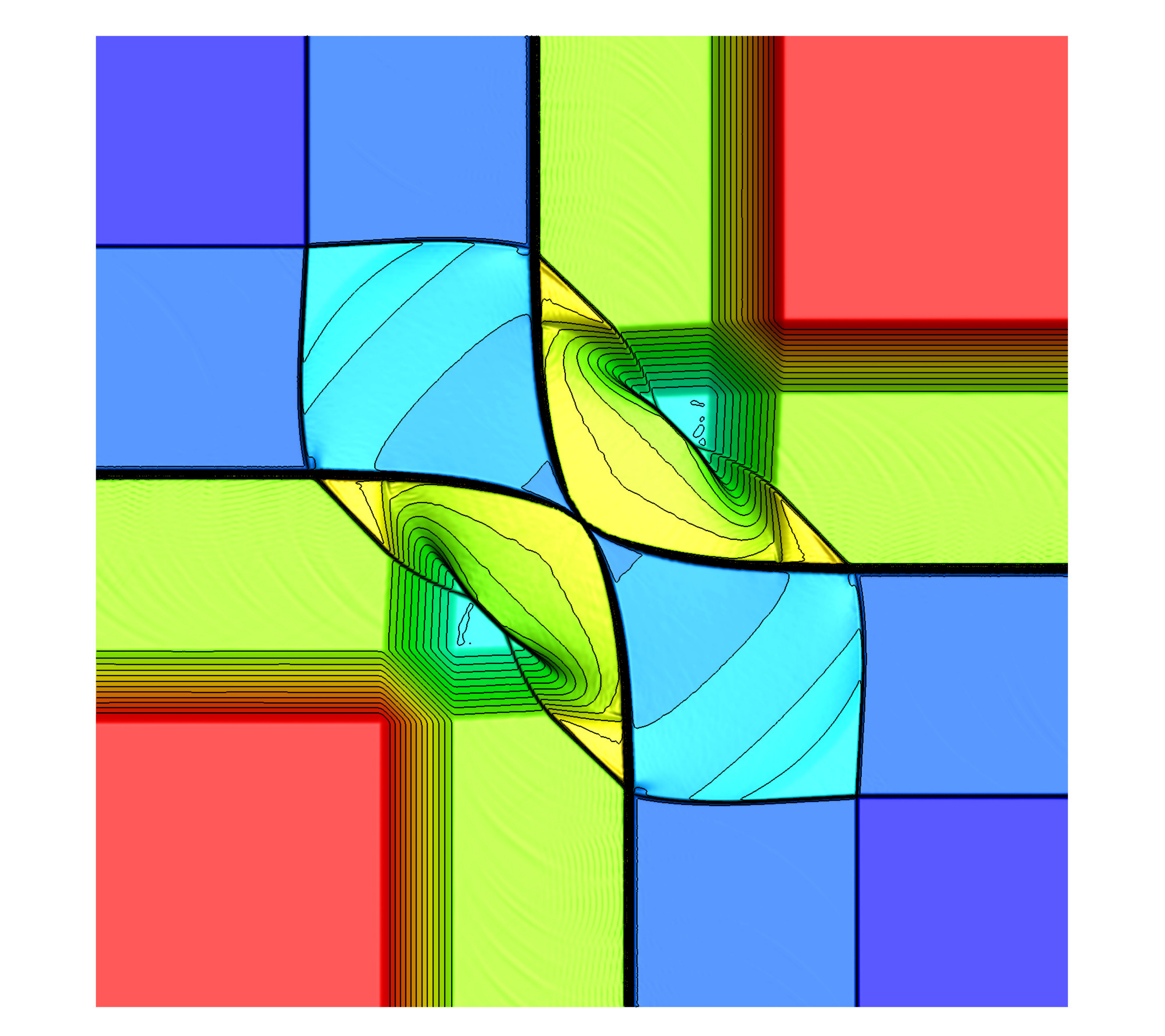}    \\ 
\includegraphics[width=0.45\textwidth]{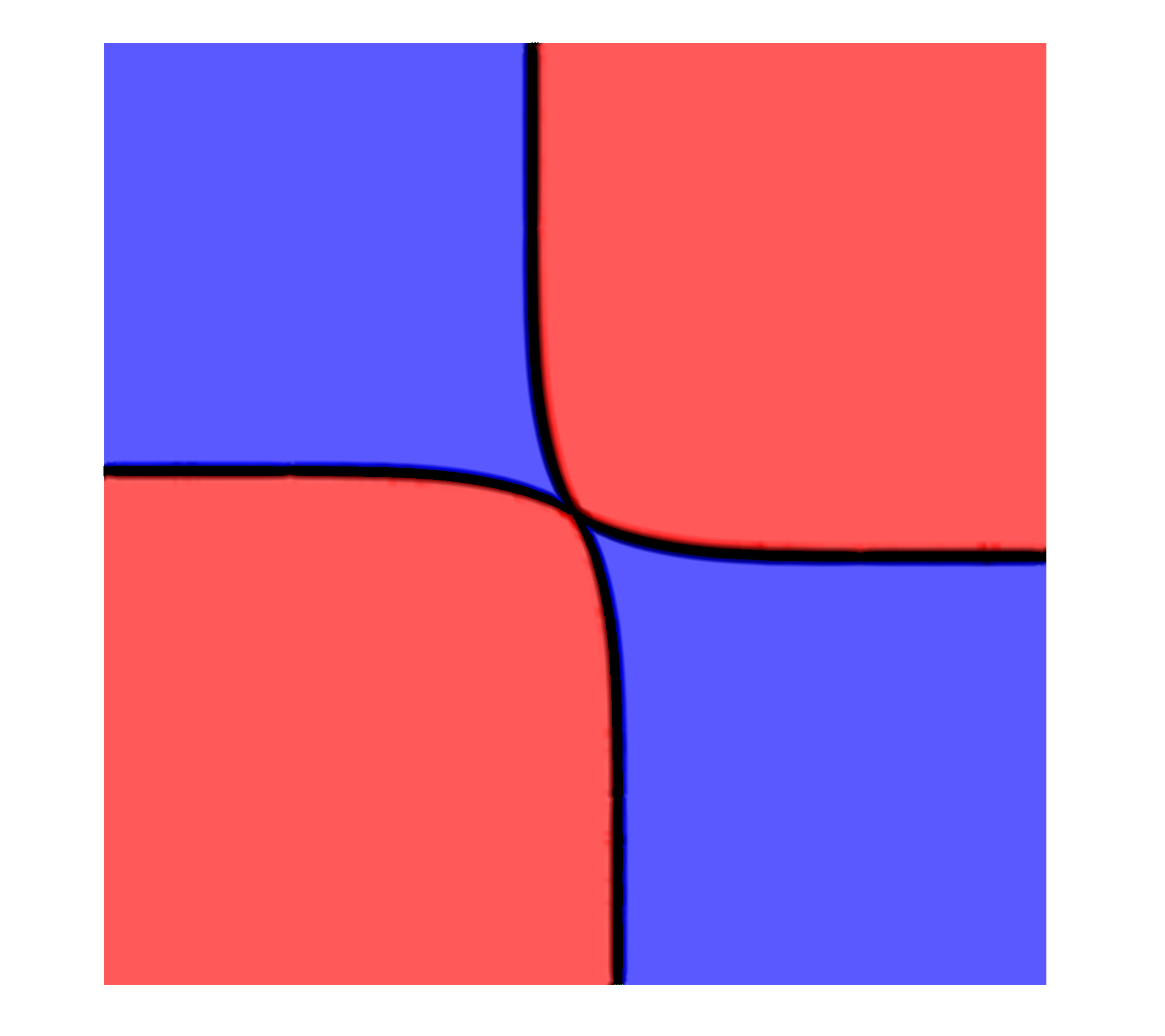}       & 
\includegraphics[width=0.45\textwidth]{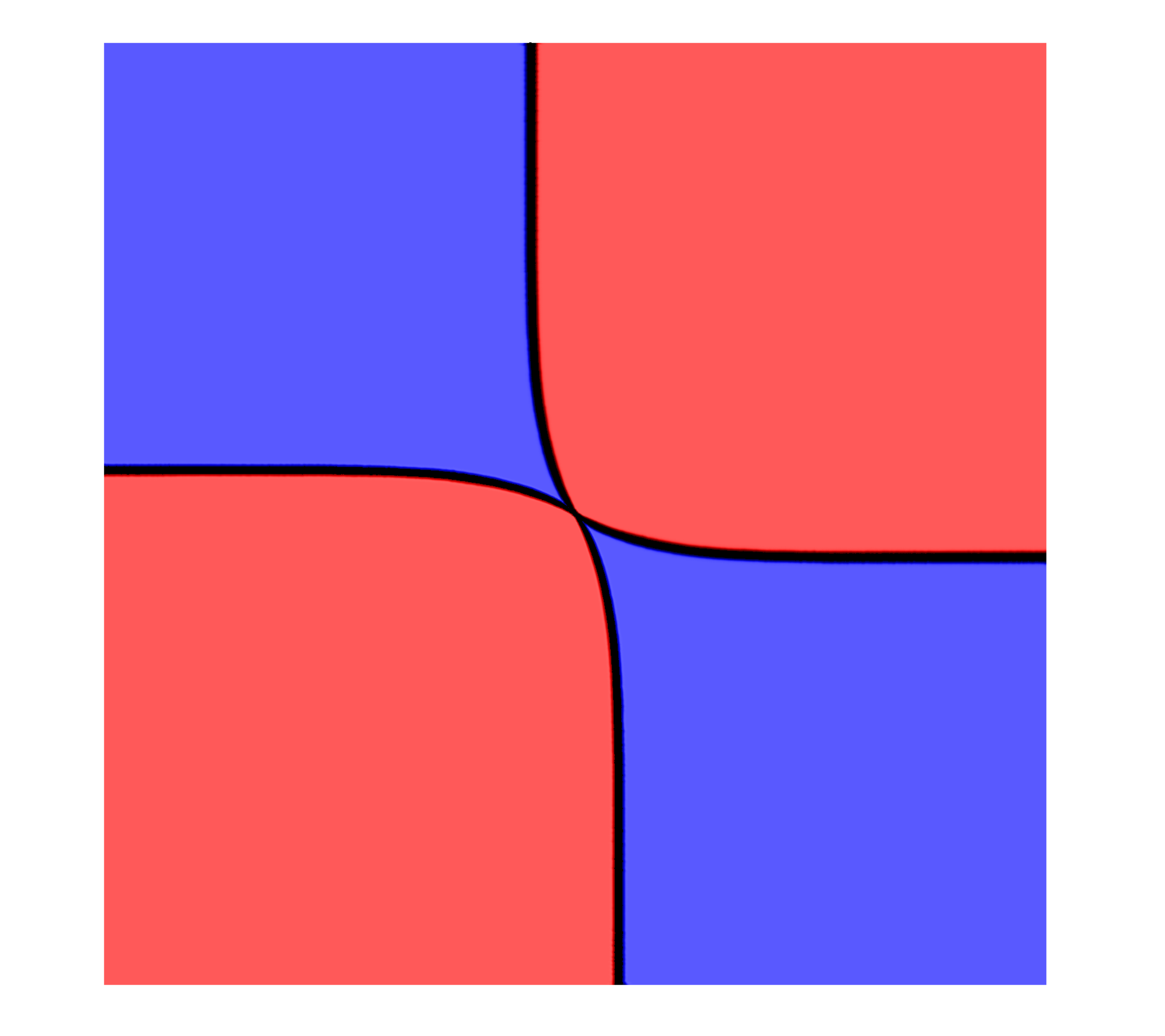}    \\ 
\end{tabular}
\caption{Results obtained with the third order Lagrangian WENO scheme for the 2D Riemann problem C1 at time $t=0.15$ (left column). 
The reference solution computed with an Eulerian method on a very fine mesh is also shown (right column). 
30 equidistant contour lines are shown for the solid density $\rho_s$ (top row), the gas density $\rho_g$ (middle row) and the 
solid volume fraction $\phi_s$ (bottom row).} 
\label{fig.bn.rp2d1}
\end{center}
\end{figure}

\begin{figure}[!htbp]
\begin{center}
\begin{tabular}{cc} 
\includegraphics[width=0.45\textwidth]{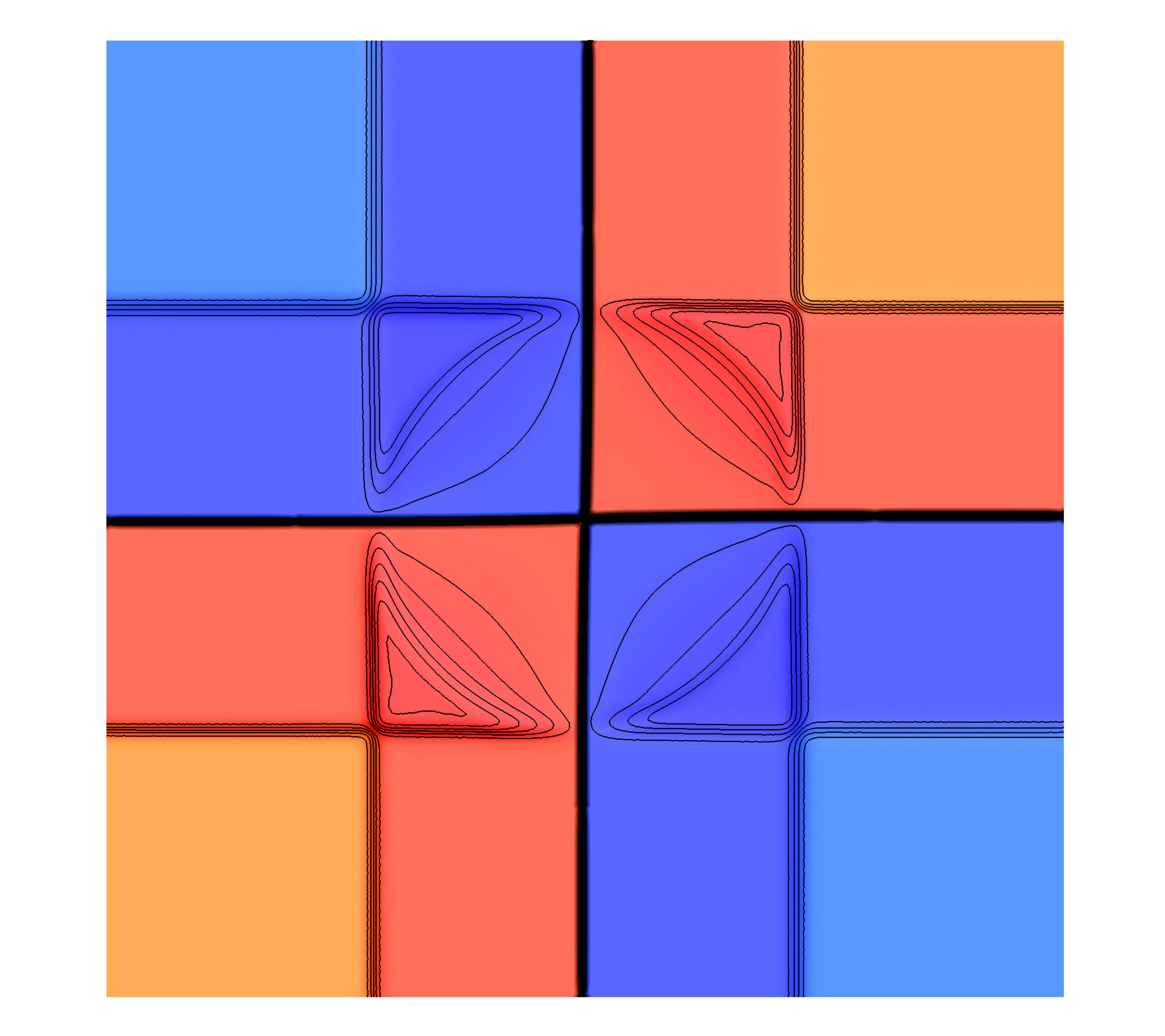}       & 
\includegraphics[width=0.45\textwidth]{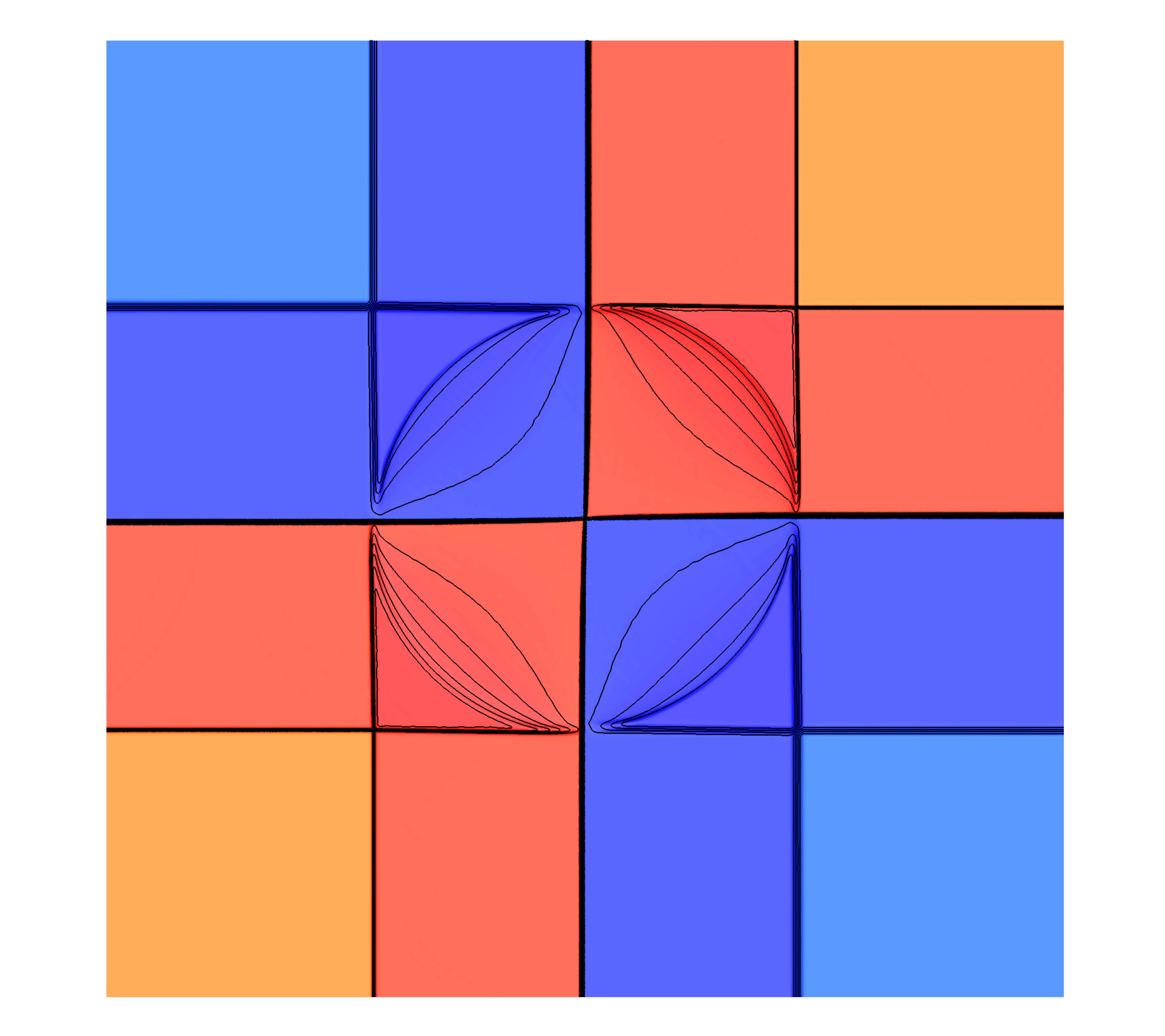}    \\ 
\includegraphics[width=0.45\textwidth]{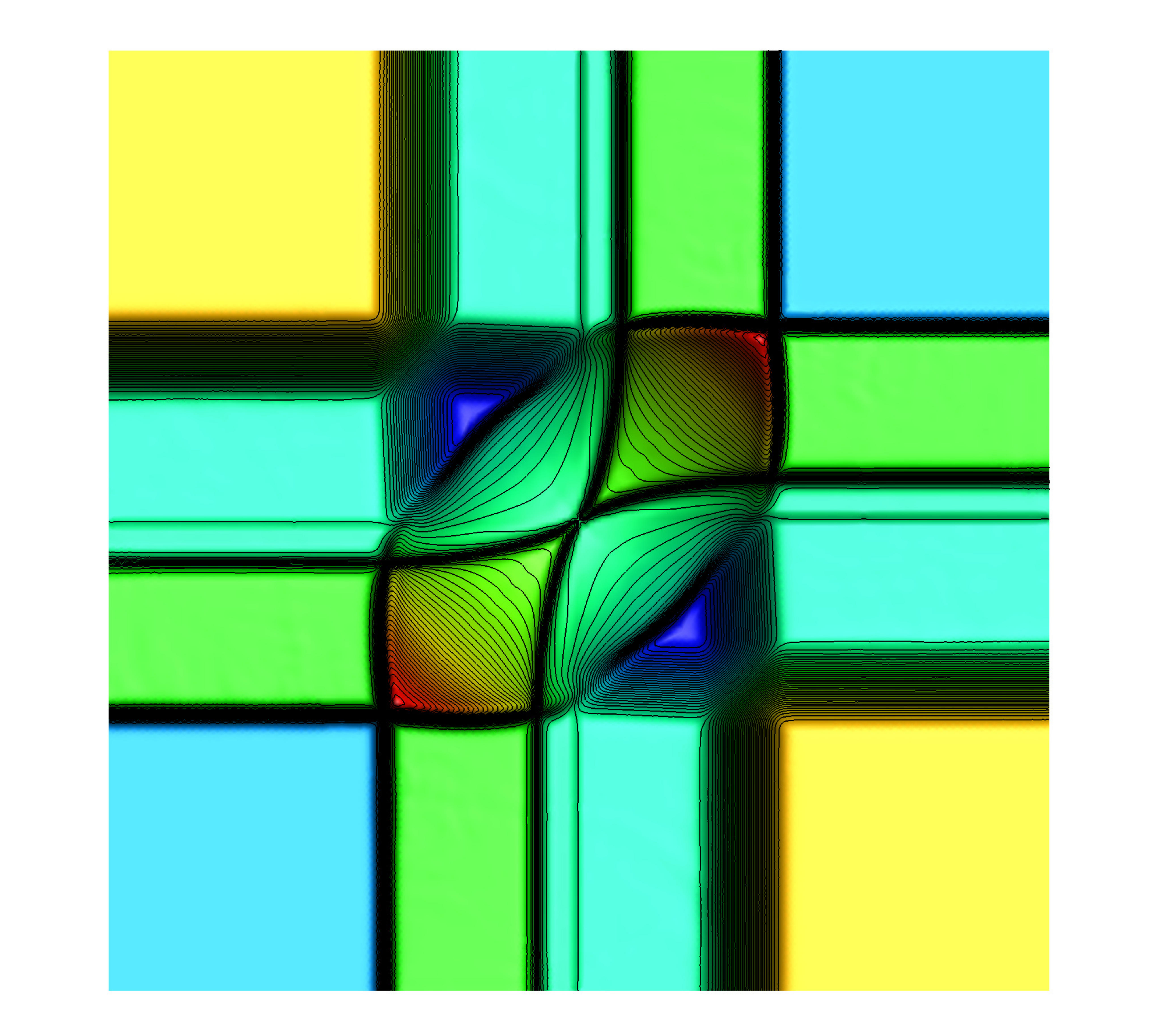}       & 
\includegraphics[width=0.45\textwidth]{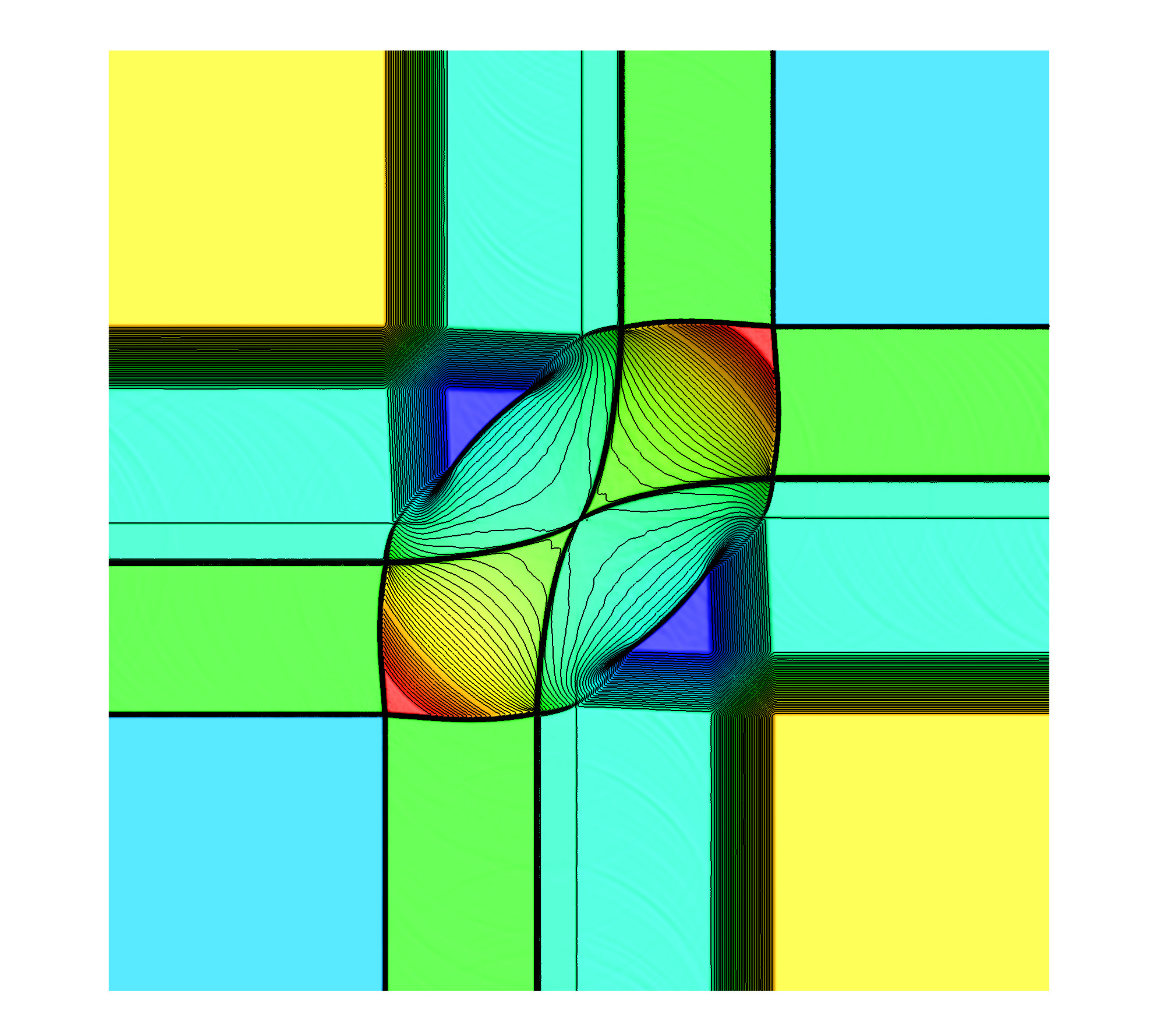}    \\ 
\includegraphics[width=0.45\textwidth]{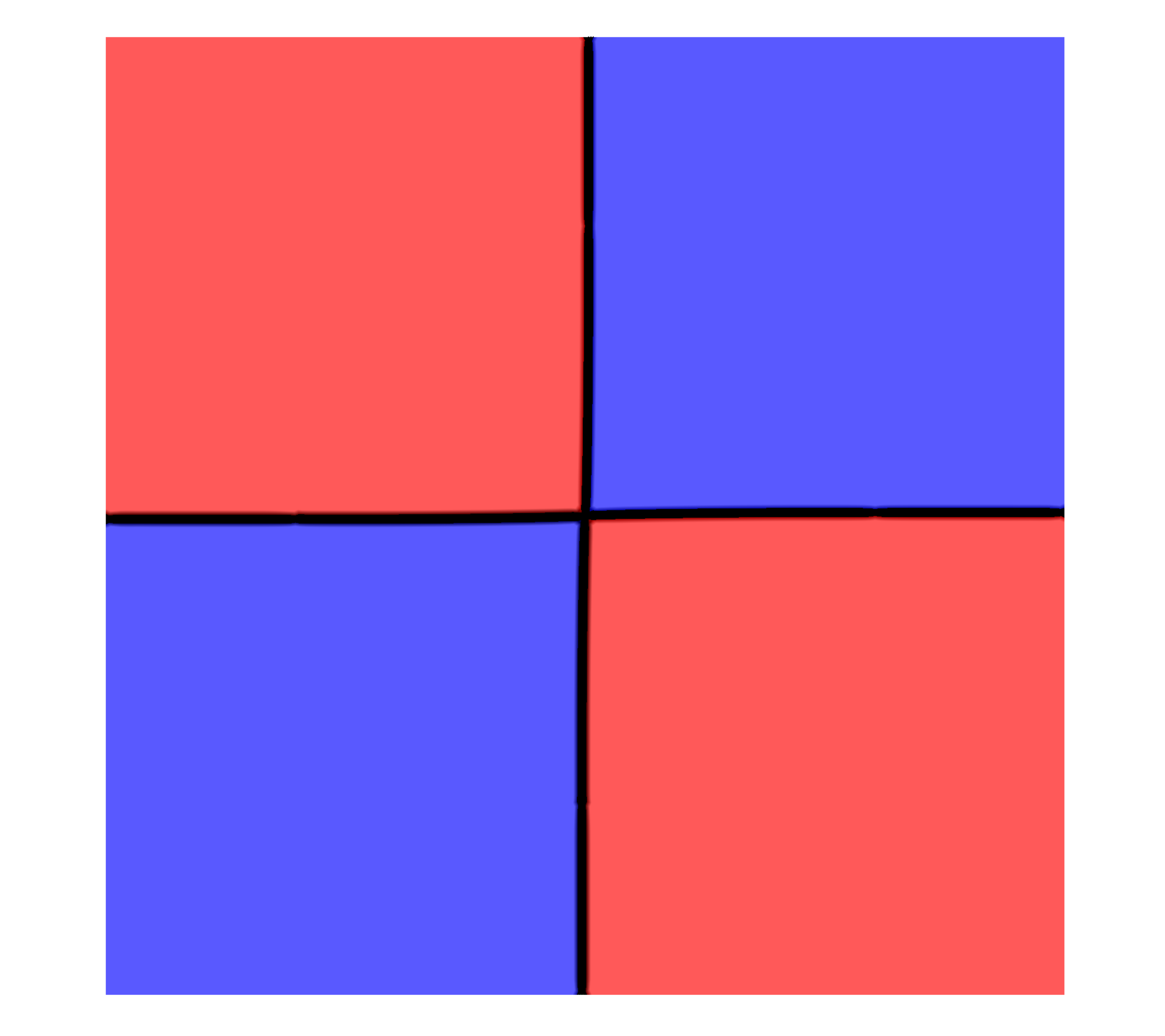}       & 
\includegraphics[width=0.45\textwidth]{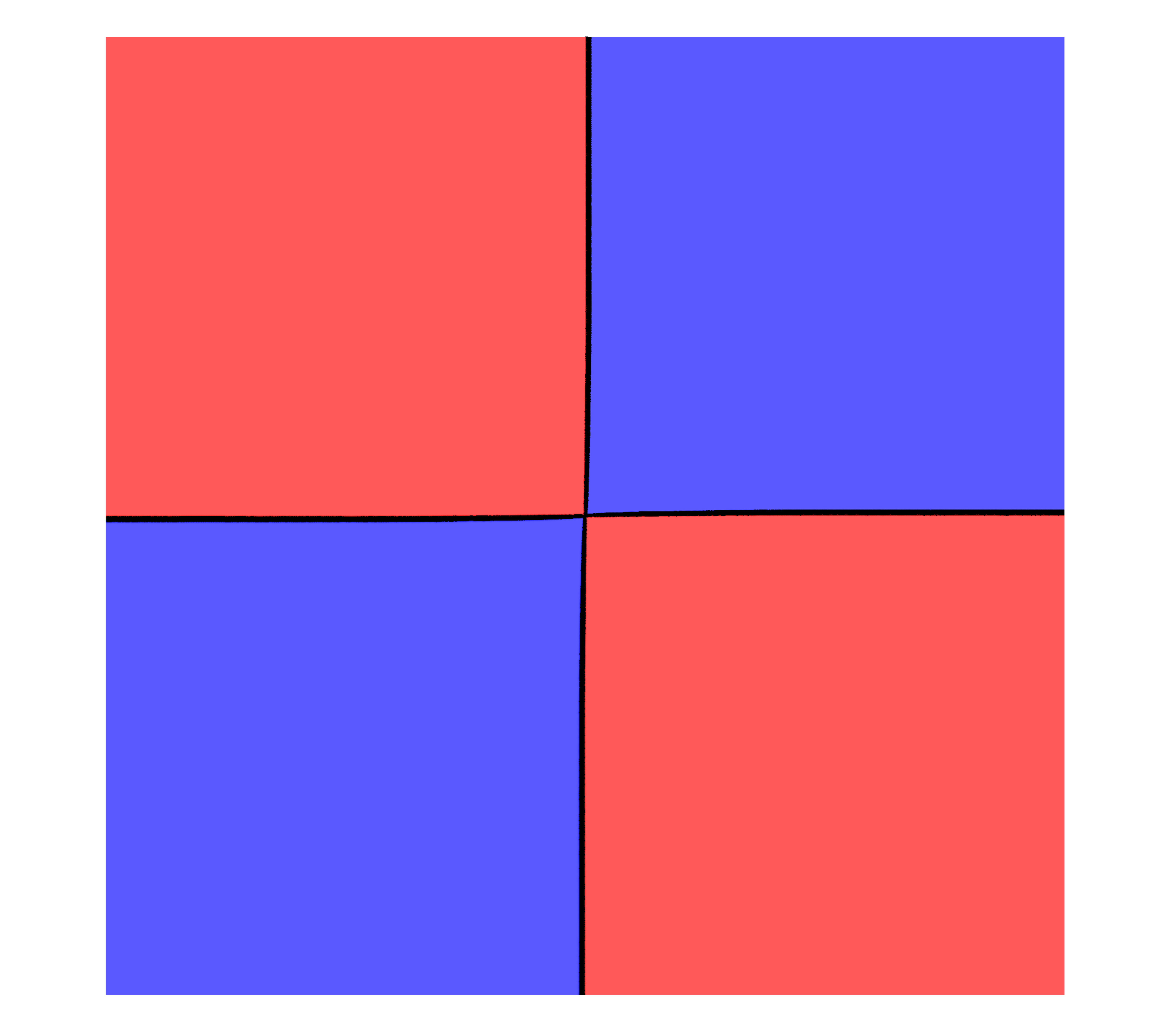}    \\ 
\end{tabular}
\caption{Results obtained with the third order Lagrangian WENO scheme for the 2D Riemann problem C2 at time $t=0.15$ (left column). 
The reference solution computed with an Eulerian method on a very fine mesh is also shown (right column). 30 equidistant contour 
lines are shown for the solid density $\rho_s$ (top row), the gas density $\rho_g$ (middle row) and the solid volume fraction 
$\phi_s$ (bottom row).} 
\label{fig.bn.rp2d2}
\end{center}
\end{figure}

In \cite{kurganovtadmor} Kurganov and Tadmor have collected a very nice set of numerical solutions for two--dimensional Riemann problems
of the compressible Euler equations. Here, we propose an extension of these 2D Riemann problems to the compressible Baer--Nunziato model. 
The initial computational domain is $\Omega(0)=[-0.5;0.5] \times [-0.5;0.5]$ and the initial condition is given by four piecewise constant 
states defined in each quadrant of the two--dimensional coordinate system:  
\begin{equation}
 \renewcommand{\arraystretch}{1.1} 
 \Q(x,0) = \left\{ \begin{array}{ccc} 
 \Q_1 & \textnormal{ if } & x > 0 \wedge y > 0,    \\ 
 \Q_2 & \textnormal{ if } & x \leq 0 \wedge y > 0, \\ 
 \Q_3 & \textnormal{ if } & x \leq 0 \wedge y \leq 0, \\ 
 \Q_4 & \textnormal{ if } & x > 0 \wedge y \leq 0.    
  \end{array} \right. 
\end{equation} 
The initial conditions for the two configurations presented in this article are listed in Table \ref{tab.rp2d.ic}. 
The simulations are carried out with a third order one--step ALE WENO finite volume scheme using an unstructured triangular mesh 
composed of 90,080 elements with an initial characteristic mesh spacing of $h=1/200$. The reference solution is computed with 
a high order Eulerian one--step scheme as presented in \cite{USFORCE,OsherNC}, using a very fine mesh composed of 2,277,668 triangles 
with characteristic mesh spacing $h=1/1000$. Reflective wall boundaries are applied on the four boundaries of the domain. 
The obtained results together with the Eulerian reference solution are depicted in Figures \ref{fig.bn.rp2d1} - \ref{fig.bn.rp2d2}, 
where we can note a very good qualitative agreement of the Lagrangian solution with the Eulerian fine--grid reference solution. 
For the first test problem, the initial and the final mesh are depicted in Fig. \ref{fig.bn.rp2d1.mesh}. 

\begin{table}[!t]   
\caption{Initial conditions for the two--dimensional Riemann problems.} 
\begin{center} 
\renewcommand{\arraystretch}{1.0}
\begin{tabular}{c|ccccccccc} 
\hline
\multicolumn{10}{c}{Configuration C1 ($\gamma_s=1.4, \gamma_g=1.4, \pi_s=\pi_g=0$)} \\ 
\hline
     & $\rho_s$ & $u_s$ & $v_s$  & $p_s$ & $\rho_g$ & $u_g$ & $v_g$  & $p_g$ & $\phi_s$  \\ 
\hline
$\Q_1: \,\, (x > 0, y > 0)$  & 2.0 & 0.0 & 0.0 & 2.0 & 1.5 & 0.0 & 0.0 & 2.0 & 0.8 \\ 
$\Q_2: \,\, (x < 0, y > 0)$  & 1.0 & 0.0 & 0.0 & 1.0 & 0.5 & 0.0 & 0.0 & 1.0 & 0.4 \\
$\Q_3: \,\, (x < 0, y < 0)$  & 2.0 & 0.0 & 0.0 & 2.0 & 1.5 & 0.0 & 0.0 & 2.0 & 0.8 \\ 
$\Q_4: \,\, (x > 0, y < 0)$  & 1.0 & 0.0 & 0.0 & 1.0 & 0.5 & 0.0 & 0.0 & 1.0 & 0.4 \\
\hline
\multicolumn{10}{c}{Configuration C2 ($\gamma_s=3.0, \gamma_g=1.4, \pi_s=100, \pi_g=0$)} \\
\hline
     & $\rho_s$ & $u_s$ & $v_s$  & $p_s$ & $\rho_g$ & $u_g$ & $v_g$  & $p_g$ & $\phi_s$ \\ 
\hline
$\Q_1: \,\, (x > 0, y > 0)$  & 1000. & 0.0 & 0.0 & 600.0 & 1.0 & 0.0 & 0.0 & 1.0 & 0.3 \\ 
$\Q_2: \,\, (x < 0, y > 0)$  & 800.  & 0.0 & 0.0 & 500.0 & 1.5 & 0.0 & 0.0 & 2.0 & 0.4 \\ 
$\Q_3: \,\, (x < 0, y < 0)$  & 1000. & 0.0 & 0.0 & 600.0 & 1.0 & 0.0 & 0.0 & 1.0 & 0.3 \\ 
$\Q_4: \,\, (x > 0, y < 0)$  & 800.  & 0.0 & 0.0 & 500.0 & 1.5 & 0.0 & 0.0 & 2.0 & 0.4 \\ 
\hline
\end{tabular} 
\end{center}
\label{tab.rp2d.ic}
\end{table}

\clearpage

\section{Application to Free Surface Flows with Moving Boundaries} 
\label{sec.freesurface} 

\subsection{Reduced Three--Equation Baer--Nunziato Model} 
\label{sec.red.bn} 

The Baer--Nunziato model \eqref{ec.BN} considered in this paper can also be applied to complex non--hydrostatic free surface flows simulations, as proposed 
in \cite{DIM2D,DIM3D}. There, a reduced three--equation Baer--Nunziato model has been considered, which can be derived from the complete seven--equation model  \eqref{ec.BN} by introducing three assumptions: firstly, all pressures are assumed to be relative pressures with respect to the atmospheric reference pressure, which is set to zero, i.e. $p_0:=0$; secondly, the gas phase that surrounds the liquid is supposed to remain always at atmospheric reference conditions, hence  obtaining $p_2 = p_0 = 0$, which allows to neglect all evolution equations related to the gas phase; finally, the pressure of the liquid phase is computed by  the Tait equation of state \cite{Batchelor1974}, which reads 
	\begin{equation}
  p_1 = k_0 \left[\left(\frac{\rho_1}{\rho_0}\right)^\gamma-1\right],
  \label{eq:TaitEOS}
  \end{equation}
where $\rho_1,\rho_0$ are the liquid density and the reference liquid density at standard conditions, respectively, $k_0$ is a constant that governs the compressibility of the fluid and $\gamma$ is a parameter commonly used to fit the equation of state with experimental data. Regarding the Tait equation \eqref{eq:TaitEOS}, we should set the proper constant values for water, i.e. $k_0=3.2 \cdot 10^8$ $Pa$, $\rho_0=1000$ $kg/m^3$ and $\gamma=7$, so that the typical speed of sound in water ($c=1500$ $m/s$) is preserved. However, since for most environmental free surface flows the maximum velocity is around $10$ $m/s$, or even less, we would have to deal with extremely low Mach numbers of $M=|\mathbf{u}|/c \ll 1$ that would cause very low time stepping and excessive numerical dissipation. In order to avoid this problem, we set artificially the Mach number to $M=0.3$, as done in \cite{DIM2D,DIM3D}, and therefore we admit density fluctuations of the order of $10\%$. The advantage of using a weakly compressible model for free surface flows is the possibility to simulate environmental--type free surface flows as well as high speed industrial free surface flows, such as they appear in water jet cutting machines or fuel 
injection systems. There, speeds may reach up to $|\mathbf{u}| \approx 1000$ m/s and thus compressibility effects can no longer 
be neglected in the liquid.  

By using the above mentioned simplifications and introducing them into system \eqref{ec.BN}, as outlined in \cite{DIM2D,DIM3D}, one can write the reduced Baer-Nunziato model as 
\begin{equation}
\begin{aligned}
&\frac{\partial}{\partial t}\left(\phi \rho \right)+\nabla \cdot \left(\phi \rho \mathbf{u}\right) = 0, \\
&\frac{\partial}{\partial t}\left(\phi \rho \mathbf{u}\right)+\nabla \cdot \left(\phi\left(\rho\mathbf{u}\mathbf{u}+\boldsymbol{\sigma} \right)\right) = \phi \rho \mathbf{g}, \\
&\frac{\partial}{\partial t}\phi + \mathbf{u}\cdot\nabla\phi = 0, 
\label{eq:ViscreducedBN}
\end{aligned}
\end{equation}
with the stress tensor $\boldsymbol{\sigma} = p \mathbf{I}$ of the inviscid liquid phase, for which we have dropped the subscript $1$ to ease 
notation. Mass and momentum equations are fully conservative in the system above, while the advection equation for the volume fraction is 
non--conservative. An evolution equation for total energy is not required, since we have used the assumption of a polytropic equation of state, 
which prescribes pressure as a function of density only. 
In \eqref{eq:ViscreducedBN} the state vector is $\mathbf{Q}=\left(\phi\rho,\phi\rho\mathbf{u},\phi\right)$ and $\mathbf{g}$ is the gravity vector 
acting along the vertical direction $y$, i.e. $\mathbf{g} = (0,g)$ with $g=9.81$ $m/s^2$. Further details on the derivation and on the assumptions  
of the reduced Baer--Nunziato model as well as a thorough validation against analytical solutions and experimental measurements can be found in  \cite{DIM2D,DIM3D}. The system \eqref{eq:ViscreducedBN} above involves a compressible inviscid phase, without viscous effects. However, viscosity plays 
an important role in turbulent flows with violent free surface motion and therefore the algorithm can be improved and extended by adding the viscous  
terms to the stress tensor $\boldsymbol{\sigma}$ of the momentum equation in \eqref{eq:ViscreducedBN}. 
The stress tensor of a Newtonian fluid using the hypothesis of Stokes reads 
\begin{equation}
\boldsymbol{\sigma} = \left(p+\frac{2}{3}\mu \nabla \cdot \mathbf{u} \right)\mathbf{I} - \mu \left(\nabla \mathbf{u} + \nabla \mathbf{u}^T \right).
\label{eq:stressTensor}
\end{equation} 
Here, $\mu$ denotes the dynamic viscosity, which for water is normally set to $10^{-3}$ $N\cdot s / m^2$. The discretization of the 
above model \eqref{eq:ViscreducedBN}, which now also contains the viscous effects, is done according to \cite{ADERNSE,ADERVRMHD,HidalgoDumbser}
and can be directly inserted into the high order one--step approach used in this paper.  

\subsection{Sloshing in a Moving Tank}
\label{sec.sloshing} 

We apply the reduced Baer--Nunziato model with viscosity \eqref{eq:ViscreducedBN} to a well known free surface flow problem in moving geometries, namely to 
the sloshing motion of a liquid in a partially filled closed tank. Such a problem can not be described by the commonly used shallow water equations, since  non--hydrostatic effects cannot be neglected. Furthermore, a numerical method is required that is able to deal with moving geometries, such as the present 
high order ALE finite volume scheme. The flow is rather complex, characterized by the presence of high amplitude oscillations and, eventually, also wave  breaking may occur. Since our algorithm is designed to be an Arbitrary Lagrangian--Eulerian (ALE) scheme, for the present problem we decide to move the 
mesh according to the motion of the sloshing tank, which is prescribed on the boundary of the domain $\partial \Omega(t)$. Inside the domain $\Omega(t)$, 
the vertices are displaced smoothly according to the \textit{Laplace equation}. In this case the evolution of the mesh and in particular the computation of 
the mesh velocity vector $\mathbf{V}(\x,t)$ is governed by the following system of elliptic PDE:  
\begin{equation}
\nabla^2 \mathbf{V}(t) = 0,  \qquad \mathbf{V}(t)=\mathbf{V}_{\!D}(t) \quad \textnormal{for} \quad \x \in \partial \Omega(t), 
\label{eq:LaplaceMesh}
\end{equation}
which is solved in each time step and where the prescribed domain motion $\mathbf{V}_{\!D}(t)$ is imposed as Dirichlet boundary condition. The solution is 
obtained using a classical P1 finite element method (FEM), where the unknowns are located at the grid nodes and the solution of the discretized system  \eqref{eq:LaplaceMesh} is computed by the conjugate gradient method. In this section the fluid motion and the mesh motion are governed by two different 
and independent laws, namely by the governing PDE system \eqref{eq:ViscreducedBN} and the Laplace equation \eqref{eq:LaplaceMesh}, respectively. By 
solving \eqref{eq:LaplaceMesh} one obtains the velocity vector $\mathbf{V}$ for each grid vertex, which allows us to move the mesh nodes and to evolve 
the geometry, i.e. to compute all the other geometric quantities needed for the computation (normal vectors, element volumes, side lengths, barycenter  positions, \textit{etc.} ). From the solution of \eqref{eq:LaplaceMesh} the space--time control volumes needed for the local space--time Galerkin 
predictor and the final one--step finite volume scheme are known a priori. 
 
The motion of liquid sloshing is a classical and widely investigated problem, since it occurs in many real world conditions whenever a sway tank is present, as in cargo ships, in liquid tank carriages on rail roads or in propellant tanks of rockets and airplanes engines. Sloshing phenomena have gained recent attention in coastal and offshore engineering with the proliferation of liquefied natural gas (LNG) and oil carriers transporting liquids in partially filled tanks \cite{LNGsloshing2012}, in order to assure and guarantee the safety of the sea transport of those fuels. The liquid sloshing motion can become very complex, even including wave breaking phenomena that can create highly localized impact pressure on tank walls, hence causing structural damages and affecting the  stability of the vehicle which carries the container. Sloshing phenomena have been widely and intensively investigated in the last decades, focusing on the development of an analytical study based on the potential flow theory, as done by Faltinsen \cite{Faltinsen1978}, who derived a linear analytical solution for liquid sloshing in a horizontally excited two--dimensional rectangular tank. Faltinsen and Timokha \cite{FaltinsenTimokha2001} extend the previous work to the nonlinear sloshing case and Hill \cite{Hill2003} analyzed more in detail the waves' behavior by relaxing some of the simplifications introduced in the previous studies. Also laboratory measurements of wave height and hydrodynamic pressure have been collected and reported \cite{Verhagen1965,Okamoto1990,Okamoto1997,Akyildiz2005}, which are very useful to validate both, theoretical solutions and numerical results. The theoretical analyses are however not suitable to describe real fluid sloshing, where viscosity and turbulence occur, so that a lot of research has been carried out in order to develop numerical models. In \cite{Faltinsen1978} Faltinsen developed a boundary element method (BEM), while Nakayama and Washizu \cite{Nakayama1980} analyzed the non-linear liquid sloshing in a two--dimensional rectangular tank under pitch excitation by using the finite element method (FEM). Wu et al. \cite{Wu1998} were the first who conducted a series of three--dimensional demonstrations on liquid sloshing based on FEM, but they did not compare the results to any experimental data. Another numerical technique is given by the finite difference method (FDM) with the use of coordinate transformations, as proposed by Chen et al. \cite{Chen1996}, who adopted a curvilinear coordinate system to map the sloshing from the non-rectangular physical domain into a rectangular computational domain. One can also solve the Navier--Stokes equations for viscous liquid sloshing \cite{ChenChiang1999,Chen2005,ChenNokes2005} or the Reynolds Averaged Navier--Stokes equations (RANS) as proposed by Armenio et al. \cite{Armenio1996}, who observed that the RANS model provides much more accurate results than the classical shallow water equations (SWE). For liquid sloshing the vertical acceleration is indeed very important and can not be neglected. In   literature one can also find Lagrangian algorithms for sloshing phenomena, see the work of Okamoto et al. \cite{Okamoto1990,Okamoto1997}, who were the first 
to apply an Arbitrary Lagrangian--Eulerian finite element method to two-- and three--dimensional liquid sloshing, and also meshless Lagrangian particle
methods, such as smoothed particle hydrodynamics (SPH) algorithms are used, as proposed in \cite{Shao2012}. Multi--fluid sloshing has been studied with the
Particle Finite Element Method (PFEM) in \cite{PFEM5}.     

From the above mentioned references we know that the intensity and the impact of sloshing do generally depend on the amplitude and the frequency of the tank motion, on the physical properties of the fluid and on the geometry of the tank, i.e. its shape.  Furthermore, liquid sloshing occurs when a tank is partially  filled with fluid, so that a free surface is present. Very frequently high Reynolds numbers are achieved in such phenomena, therefore turbulence is usually present in this kind of flows. We improve indeed the algorithm by adding the simple zero-equation Smagorinsky turbulence model \cite{Smagorinsky1963}, 
which gives an expression for the eddy viscosity $\mu_T$ that is assumed to be proportional to the velocity gradients. For the two--dimensional case the  Smagorinsky model reads
\begin{equation}
\mu_T = \rho \cdot \left(C_S\Delta\right)^2|\overline{\mathbf{S}}|, 
\label{eq:SmagorinskyMod}
\end{equation} 
where $C_S$ is a coefficient which has to be set properly, $\Delta$ is the turbulence length scale and the term $|\overline{\mathbf{S}}|$ is given by
\begin{equation}
|\overline{\mathbf{S}}| = \sqrt{2\overline{S_{ij}}\overline{S_{ij}}}, \qquad \overline{S_{ij}} = \frac{1}{2}\left(\frac{\partial u_i}{\partial x_j}+\frac{\partial u_j}{\partial x_i}\right),
\label{eq:strainTens}
\end{equation}
with the strain rate tensor $\overline{S_{ij}}$. The Smagorninsky constant is usually in the range $\left[0.1;0.24\right]$ and here we set 
$C_S=0.17$. In this paper show numerical results for an idealized two--dimensional case, where the tank is moving with a purely horizontal sinusoidal 
velocity and the vertical velocity component is zero, i.e. 
\begin{equation}
\mathbf{V}_{\!D}(t) = \left( V_A \sin\left(\omega t\right), 0 \right)^T
\label{eq:sloshingVel}
\end{equation}
where $\omega=\frac{2\pi}{T}$ is the frequency of the oscillation and $T$ is the period, while $V_{A}$ quantifies the amplitude of the sinusoidal 
velocity, which is related to the horizontal tank displacement $\delta x_A$ by $V_{A} = -\omega \delta x_A$. We always assume the fluid and the 
tank to be initially at rest, i.e. $\mathbf{V}(0)=0$. 

\begin{figure}
	\centering
		\includegraphics[width=0.6\textwidth]{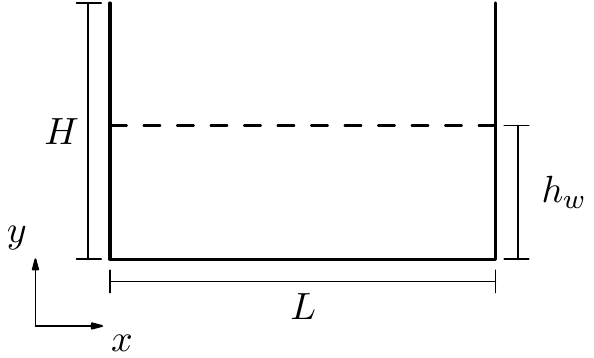}
	\caption{Sloshing tank.}
	\label{fig:SloshingGeom}
\end{figure}

The simulations presented in this paper refer to the very recent work of Shao et al. \cite{Shao2012}, who present an improved SPH method for modeling 
liquid sloshing dynamics. Their main novelty regards the use of the Reynolds averaged turbulence model incorporated into the SPH method in order to 
properly describe the effects of turbulence. Furthermore, in \cite{Shao2012} some new density and kernel gradient corrections have been introduced to 
achieve better accuracy and smoother pressure fields. The initial domain $\Omega(0)$ refers to the one used for the laboratory measurements carried out 
by Faltinsen et al. \cite{FaltinsenTimokha2000} and is represented by a two--dimensional rectangular tank of dimensions 
$\Omega(0)=[0;1.73] \times [0;1.15]$, as sketched in Figure \ref{fig:SloshingGeom}. We use an unstructured triangular mesh with a total number of 
$N_e=44746$ elements.
The lateral sides and the bottom of the domain are modelled by classical no--slip wall boundary conditions for viscous flow, while the above boundary is set 
as transmissive outflow. The computational domain is characterized by a periodic excitation along the horizontal direction $x$, whose displacement $\delta x(t)$ is described as $\delta x(t) = \delta x_A \cdot \cos(\omega t)$, where we fix the parameters $\delta x_A = 0.032$ and $T=1.3$ with the initial water  height $h_w=0.6$, as proposed in \cite{Shao2012}. The corresponding velocity law for the tank can be easily computed by \eqref{eq:sloshingVel}, with  
$V_{A} = -\omega \delta x_A$. In this case we are dealing with a turbulent flow, since the Reynolds number $Re = \frac{L\cdot u}{\nu}$ is of the order of $10^{6}$. Regarding the parameters of the equation of state \eqref{eq:TaitEOS}, we set $k_0=8.5\cdot 10^4$ $Pa$, $\rho_0=1000$ $kg/m^3$ and $\gamma=1$, as suggested in \cite{DIM2D,DIM3D}. 

\begin{figure}
	\centering
		\includegraphics[width=0.9\textwidth]{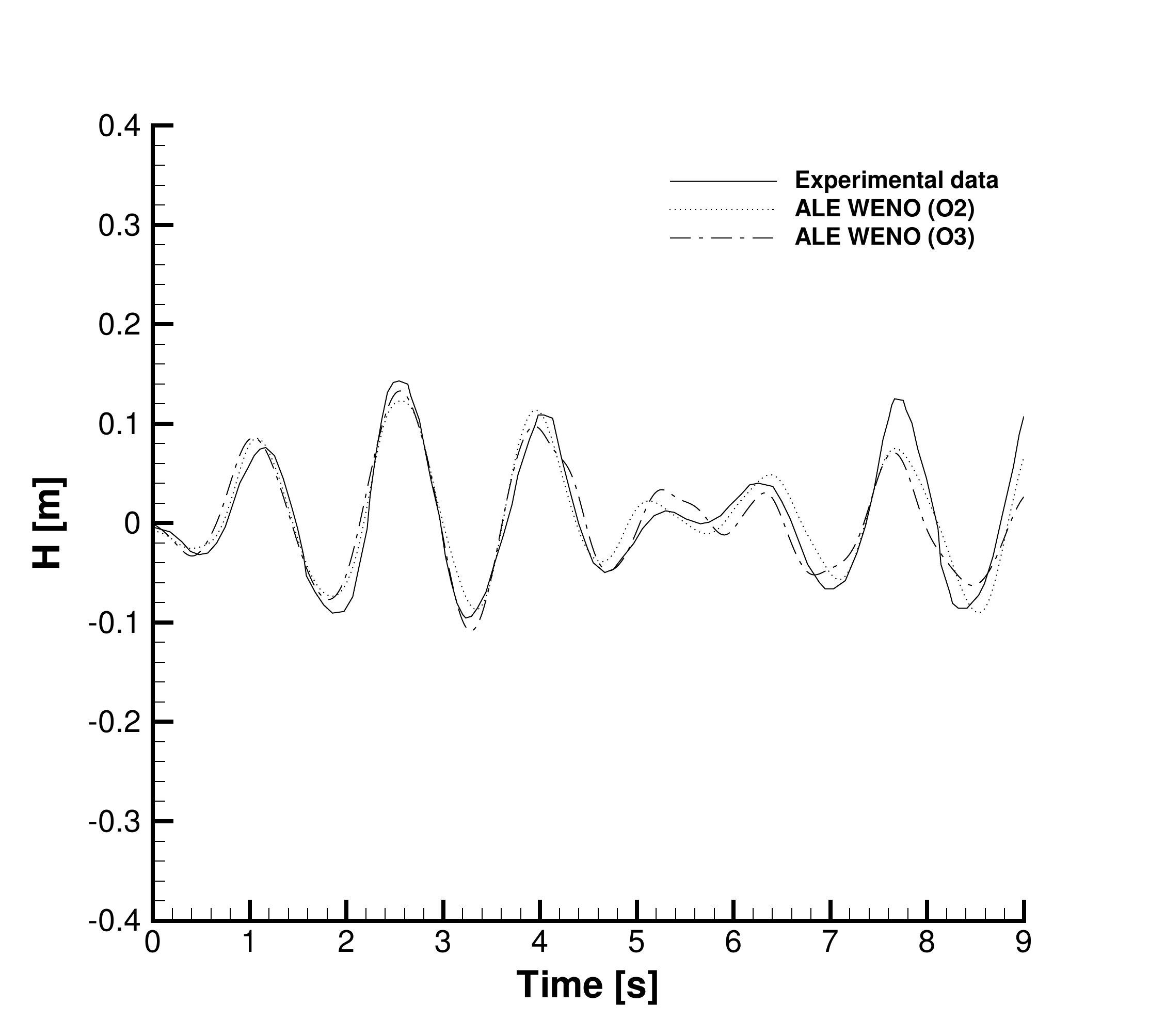}
	\caption{Sloshing motion: comparison between experimental data (solid line), second order (dotted line) and third order ALE WENO scheme (dash-dot line). $T=1.3$ $s$ and $h_w=0.6$ $m$.}
	\label{fig:slosh_T13}
\end{figure}

Figure \ref{fig:slosh_T13} shows a comparison between experimental data and numerical results: the periodic motion of the tank can be clearly identified looking at the perturbation $H$ of the free surface location with respect to the initial configuration. The numerical results have been collected at the same probe  point $\x_p$ used for the experimental data, which is placed on the initial free surface level and is $0.05$ $m$ away from the left wall, hence its initial location is $\x_p(0)=(0.05,0.6)$. Since the probe point is attached to the tank, it moves together with the domain. Furthermore the value of $H$ is evaluated at each time step as 
\begin{equation}
H = \int\limits_{h_\Omega} \phi(s) ds - h_w, \qquad \textnormal{ with } \qquad h_\Omega\in [0,1.15]
\label{eq:H_int}
\end{equation}
where the volume fraction integral along the whole height $h_\Omega$ of the domain represents indeed the water column at the probe point which is shifting horizontally according to Eqn.\eqref{eq:sloshingVel} and $h_w$ is the initial free surface elevation. As the tank motion begins, water moves  towards the left, hence decreasing the pressure ($t=0.5$), then the tank changes completely direction shifting towards the right side and therefore the free surface elevation increases its level up to $0.7$ $m$ at time $t=1.2$ $s$. The water keeps moving following the described periodic cycle of the tank boundaries, progressively increasing the wave amplitude in time. The numerical results are in good agreement with the experimental data, also with the  particular wiggles observed in the time interval $t \in [5;7]s$, see Fig. \ref{fig:slosh_T13}. 

\begin{figure}
	\centering
		\includegraphics[width=1.0\textwidth]{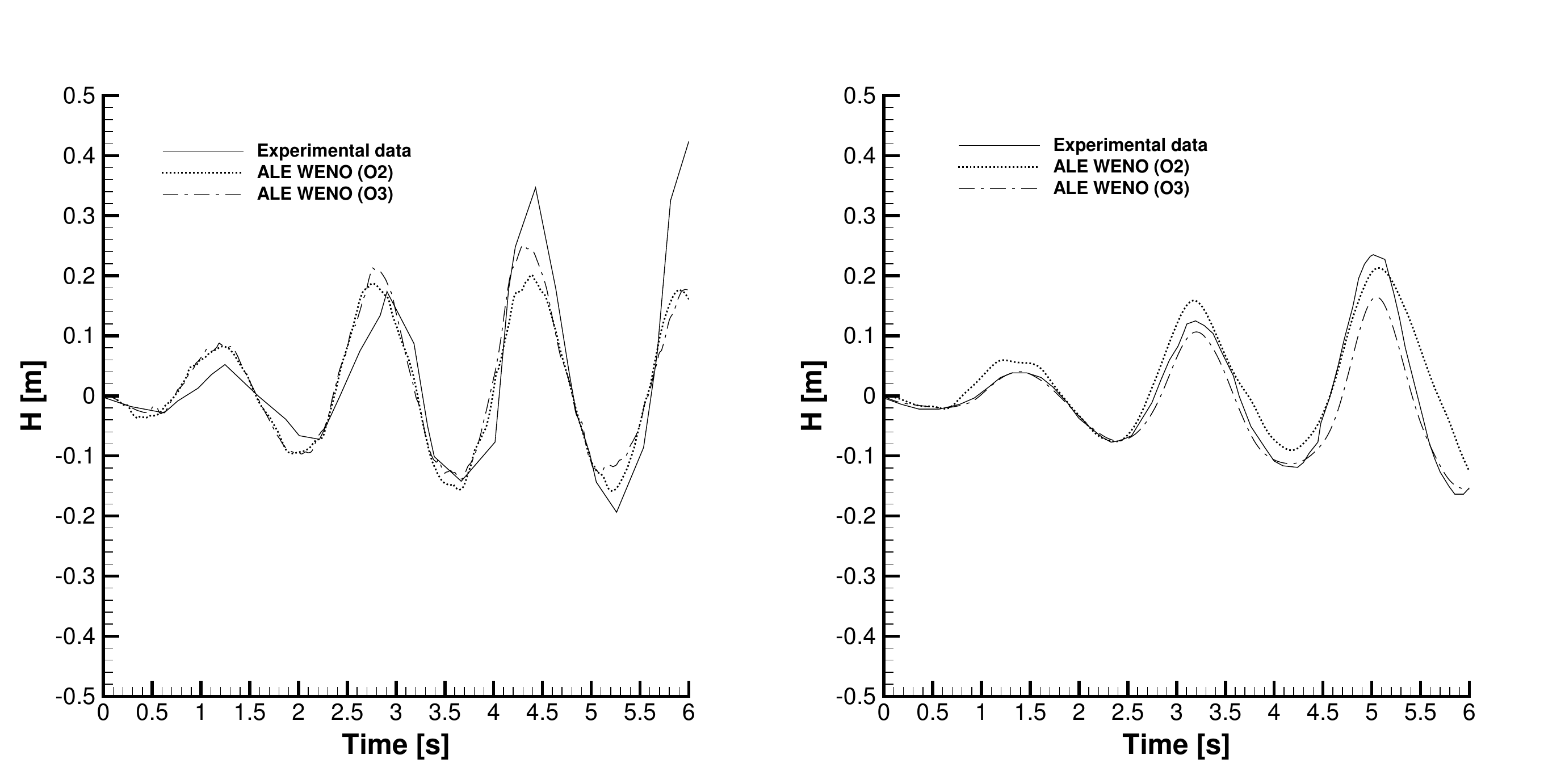}
	\caption{Comparison between experimental data (solid line) and numerical results with second (dotted line) and third (dash-dot line) order of accuracy. Left: $T=1.5$ $h_w=0.5$. Right: $T=1.875$ $h_w=0.5$.}
	\label{fig:comparison}
\end{figure}

Another comparison of numerical results with experimental data is depicted in Figure \ref{fig:comparison} for two different sets of parameters: on the left we show the case $T=1.5 s$ and $h_w=0.6 m$, while on the right we use $T=1.875 s$ and $h_w=0.5 m$. 
\begin{figure}
	\centering
		\includegraphics[width=0.95\textwidth]{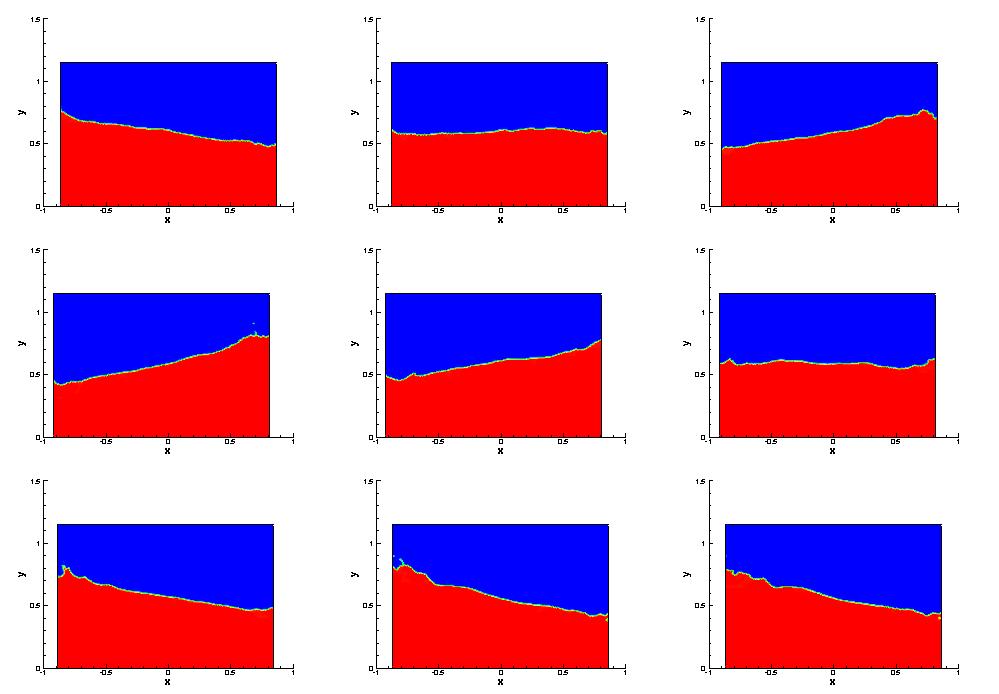}
	\caption{Water motion for the case $T=1.5 s$ and $h_w=0.6 m$ at times $t=3.0$, $3.2$, $3.4$, $3.6$, $3.8$, $4.0$, $4.2$, $4.4$ and $4.5$ $s$.}
	\label{fig:SloshPattern}
\end{figure}
The flow pattern of the sloshing for $T = 1.5 s$ and $h_w = 0.6 m$ is plotted in Figure \ref{fig:SloshPattern} at nine representative time instants within one  period. One can notice the strong oscillations of the free surface occurring in the flow while the tank is swinging, as well as the motion of the domain $\Omega(t)$.  

\section{Conclusions}
\label{sec.concl}

In this article the first better than second order Arbitrary--Lagrangian-Eulerian one--step WENO finite volume scheme on unstructured triangular meshes 
has been proposed for the solution of hyperbolic systems with non--conservative products. The method has been applied to the full seven equation 
Baer--Nunziato model of compressible multiphase flows as well as to a reduced three--equation model for the simulation of weakly compressible free 
surface flows in moving geometries. 
High order of accuracy in space and time have been verified by a numerical convergence study on a smooth unsteady test problem where an exact solution 
of the Baer--Nunziato model is available. The scheme has also been successfully applied to problems with shock waves and material interfaces. For all
test problems the numerical results have been compared with exact or numerical reference solutions in order to validate the approach. 

In future research we plan to extend the present scheme to three space dimensions in the more general framework of the new $P_{N}P_{M}$ method proposed 
in \cite{Dumbser2008}, which unifies high order finite volume and discontinuous Galerkin finite element methods in one general approach. Further work will 
also consist in a generalization to moving \textit{curved} meshes as well as the introduction of multi--dimensional Riemann solvers, such as the ones
used in \cite{StagLag,balsarahlle2d,balsarahllc2d}. Finally, more research is necessary concerning a suitable conservative and high order accurate 
remapping / remeshing strategy to deal with significant mesh distorsions arising in the presence of strong shear waves.

\section*{Acknowledgments}
The presented research has been financed by the European Research Council (ERC) under the European Union's Seventh Framework 
Programme (FP7/2007-2013) with the research project \textit{STiMulUs}, ERC Grant agreement no. 278267. Special thanks to 
Marco Nesler and Paul Maistrelli for the installation and technical support of the AMD Opteron cluster used for the simulations 
shown in this paper.


\bibliography{LagrangeNC}

\begin{thebibliography}{100}

\bibitem{abgrall_eno}
R.~Abgrall.
\newblock On essentially non-oscillatory schemes on unstructured meshes:
  analysis and implementation.
\newblock {\em Journal of Computational Physics}, 144:45--58, 1994.

\bibitem{abgrallkarni}
R.~Abgrall and S.~Karni.
\newblock A comment on the computation of non-conservative products.
\newblock {\em Journal of Computational Physics}, 229:2759--2763, 2010.

\bibitem{AboiyarIske}
T.~Aboiyar, E.H. Georgoulis, and A.~Iske.
\newblock {Adaptive ADER Methods Using Kernel-Based Polyharmonic Spline WENO
  Reconstruction}.
\newblock {\em SIAM Journal on Scientific Computing}, 32:3251--3277, 2010.

\bibitem{Akyildiz2005}
H.~Akyildiz and E.~Unal.
\newblock Experimental investigation of pressure distribution on a rectangular
  tank due to the liquid sloshing.
\newblock {\em Ocean Eng.}, 32:1503 -- 1516, 2005.

\bibitem{AndrianovWarnecke}
N.~Andrianov and G.~Warnecke.
\newblock The {Riemann} problem for the {Baer--Nunziato} two-phase flow model.
\newblock {\em Journal of Computational Physics}, 212:434--464, 2004.

\bibitem{Armenio1996}
V.~Armenio and M.~La Rocca.
\newblock On the analysis of sloshing of water in rectangular containers:
  numerical and experimental investigation.
\newblock {\em Ocean Eng.}, 23:705 -- 739, 1996.

\bibitem{BaerNunziato1986}
M.R. Baer and J.W. Nunziato.
\newblock A two-phase mixture theory for the deflagration-to-detonation
  transition {(DDT)} in reactive granular materials.
\newblock {\em J. Multiphase Flow}, 12:861–--889, 1986.

\bibitem{Balsara2004}
D.~Balsara.
\newblock Second-order accurate schemes for magnetohydrodynamics with
  divergence-free reconstruction.
\newblock {\em The Astrophysical Journal Supplement Series}, 151:149--184,
  2004.

\bibitem{balsarashu}
D.~Balsara and {C.W.} Shu.
\newblock Monotonicity preserving weighted essentially non-oscillatory schemes
  with increasingly high order of accuracy.
\newblock {\em Journal of Computational Physics}, 160:405--452, 2000.

\bibitem{balsarahlle2d}
D.S. Balsara.
\newblock {Multidimensional HLLE Riemann solver: Application to Euler and
  magnetohydrodynamic flows}.
\newblock {\em Journal of Computational Physics}, 229:1970--1993, 2010.

\bibitem{balsarahllc2d}
D.S. Balsara.
\newblock {A two-dimensional HLLC Riemann solver for conservation laws:
  Application to Euler and magnetohydrodynamic flows}.
\newblock {\em Journal of Computational Physics}, 231:7476--7503, 2012.

\bibitem{barthlsq}
T.J. Barth and P.O. Frederickson.
\newblock Higher order solution of the {Euler} equations on unstructured grids
  using quadratic reconstruction.
\newblock {\em AIAA paper no. 90-0013}, 28th Aerospace Sciences Meeting January
  1990.

\bibitem{Batchelor1974}
{G. K.} Batchelor.
\newblock {\em An Introduction to Fluid Mechanics}.
\newblock Cambridge University Press, 1974.

\bibitem{Benson1992}
D.J. Benson.
\newblock Computational methods in lagrangian and eulerian hydrocodes.
\newblock {\em Computer Methods in Applied Mechanics and Engineering},
  99:235--394, 1992.

\bibitem{MaireMM2}
M.~Berndt, J.~Breil, S.~Galera, M.~Kucharik, P.H. Maire, and M.~Shashkov.
\newblock {Two--step hybrid conservative remapping for multimaterial arbitrary
  Lagrangian–-Eulerian methods}.
\newblock {\em Journal of Computational Physics}, 230:6664--6687, 2011.

\bibitem{BoscheriDumbserLag}
W.~Boscheri and M.~Dumbser.
\newblock {Arbitrary--Lagrangian--Eulerian One--Step WENO Finite Volume Schemes
  on Unstructured Triangular Meshes}.
\newblock {\em Communications in Computational Physics}.
\newblock submitted to.

\bibitem{BoscheriDumbser}
W.~Boscheri, M.~Dumbser, and M.~Righetti.
\newblock A semi-implicit scheme for 3d free surface flows with high order
  velocity reconstruction on unstructured voronoi meshes.
\newblock {\em International Journal for Numerical Methods in Fluids}.
\newblock DOI: 10.1002/fld.3753, in press.

\bibitem{MaireMM1}
J.~Breil, S.~Galera, and P.H. Maire.
\newblock {Multi-material ALE computation in inertial confinement fusion code
  CHIC}.
\newblock {\em Computers and Fluids}, 46:161--167, 2011.

\bibitem{MaireMM3}
J.~Breil, T.~Harribey, P.H. Maire, and M.~Shashkov.
\newblock {A multi-material ReALE method with MOF interface reconstruction}.
\newblock {\em Computers and Fluids}.
\newblock DOI: 10.1016/j.compfluid.2012.08.015, in press.

\bibitem{Caramana1998}
E.J. Caramana, D.E. Burton, M.J. Shashkov, and P.P. Whalen.
\newblock The construction of compatible hydrodynamics algorithms utilizing
  conservation of total energy.
\newblock {\em Journal of Computational Physics}, 146:227--262, 1998.

\bibitem{Després2009}
G.~Carr\'e, S.~Del Pino, B.~Despr\'es, and E.~Labourasse.
\newblock A cell-centered lagrangian hydrodynamics scheme on general
  unstructured meshes in arbitrary dimension.
\newblock {\em Journal of Computational Physics}, 228:5160 -- 5183, 2009.

\bibitem{Castro2008}
M.J. Castro, J.M. Gallardo, J.A. L\'opez, and C.~Par\'es.
\newblock Well-balanced high order extensions of godunov's method for
  semilinear balance laws.
\newblock {\em SIAM Journal of Numerical Analysis}, 46:1012--1039, 2008.

\bibitem{Castro2006}
M.J. Castro, J.M. Gallardo, and C.~Par\'es.
\newblock High-order finite volume schemes based on reconstruction of states
  for solving hyperbolic systems with nonconservative products. applications to
  shallow-water systems.
\newblock {\em Mathematics of Computation}, 75:1103--–1134, 2006.

\bibitem{NCproblems}
M.J. Castro, P.G. LeFloch, M.L. {Mu\~noz-Ruiz}, and C.~Par\'es.
\newblock Why many theories of shock waves are necessary: Convergence error in
  formally path-consistent schemes.
\newblock {\em Journal of Computational Physics}, 227:8107--8129, 2008.

\bibitem{Casulli1990}
V.~Casulli.
\newblock Semi-implicit finite difference methods for the two-dimensional
  shallow water equations.
\newblock {\em Journal of Computational Physics}, 86:56--74, 1990.

\bibitem{CasulliCheng1992}
V.~Casulli and R.T. Cheng.
\newblock Semi-implicit finite difference methods for three-dimensional shallow
  water flow.
\newblock {\em International Journal of Numerical Methods in Fluids},
  15:629--648, 1992.

\bibitem{Feistauer4}
J.~Cesenek, M.~Feistauer, J.~Horacek, V.~Kucera, and J.~Prokopova.
\newblock Simulation of compressible viscous flow in time-dependent domains.
\newblock {\em Applied Mathematics and Computation}.
\newblock {DOI: 10.1016/j.amc.2011.08.077, in press.}

\bibitem{Chen2005}
{B. F.} Chen.
\newblock Viscous fluid in a tank under coupled surge, heave and pitch motions.
\newblock {\em J. Waterw. Port Coast. Ocean Eng.--ASCE}, 131:239 -- 256, 2005.

\bibitem{ChenChiang1999}
{B. F.} Chen and {H. W.} Chiang.
\newblock Complete 2d and fully nonlinear analysis of ideal fluid in tanks.
\newblock {\em J. Eng. Mech.--ASCE}, 125:70 -- 78, 1999.

\bibitem{ChenNokes2005}
{B. F.} Chen and R.~Nokes.
\newblock Time--independent finite difference analysis of 2d and nonlinear
  viscous liquid sloshing in a rectangular tank.
\newblock {\em J. comput. Phys.}, 209:47 -- 81, 2005.

\bibitem{Chen1996}
W.~Chen, {M. A.} Haroun, and F.~Liu.
\newblock Large amplitude liquid sloshing in seismically excited tanks.
\newblock {\em Earthquake Eng. Struct. Dyn.}, 25:653 -- 669, 1996.

\bibitem{chengshu1}
J.~Cheng and C.W. Shu.
\newblock {A high order ENO conservative Lagrangian type scheme for the
  compressible Euler equations}.
\newblock {\em Journal of Computational Physics}, 227:1567--1596, 2007.

\bibitem{chengshu3}
J.~Cheng and C.W. Shu.
\newblock {A cell-centered Lagrangian scheme with the preservation of symmetry
  and conservation properties for compressible fluid flows in two-dimensional
  cylindrical geometry}.
\newblock {\em Journal of Computational Physics}, 229:7191--7206, 2010.

\bibitem{chengshu4}
J.~Cheng and C.W. Shu.
\newblock {Improvement on spherical symmetry in two-dimensional cylindrical
  coordinates for a class of control volume Lagrangian schemes}.
\newblock {\em Communications in Computational Physics}, 11:1144--1168, 2012.

\bibitem{MOOD}
S.~Clain, S.~Diot, and R.~Loub\`ere.
\newblock {A high--order finite volume method for systems of conservation laws
  -- Multi--dimensional Optimal Order Detection (MOOD)}.
\newblock {\em Journal of Computational Physics}, 230:4028--4050, 2011.

\bibitem{Depres2012}
A.~Claisse, B.~Despr\'es, E.Labourasse, and F.~Ledoux.
\newblock A new exceptional points method with application to cell-centered
  lagrangian schemes and curved meshes.
\newblock {\em Journal of Computational Physics}, 231:4324--4354, 2012.

\bibitem{CBS-book}
B.~Cockburn, G.~E. Karniadakis, and {C.W.} Shu.
\newblock {\em Discontinuous {Galerkin} Methods}.
\newblock Lecture Notes in Computational Science and Engineering. Springer,
  2000.

\bibitem{CIR}
R.~Courant, E.~Isaacson, and M.~Rees.
\newblock On the solution of nonlinear hyperbolic differential equations by
  finite differences.
\newblock {\em Comm. Pure Appl. Math.}, 5:243--255, 1952.

\bibitem{DeledicquePapalexandris}
V.~Deledicque and M.V. Papalexandris.
\newblock An exact {Riemann} solver for compressible two-phase flow models
  containing non-conservative products.
\newblock {\em Journal of Computational Physics}, 222:217--245, 2007.

\bibitem{DepresMazeran2003}
B.~Despr\'es and C.~Mazeran.
\newblock Symmetrization of lagrangian gas dynamic in dimension two and
  multimdimensional solvers.
\newblock {\em C.R. Mecanique}, 331:475--480, 2003.

\bibitem{Despres2005}
B.~Despr\'es and C.~Mazeran.
\newblock Lagrangian gas dynamics in two-dimensions and lagrangian systems.
\newblock {\em Archive for Rational Mechanics and Analysis}, 178:327--372,
  2005.

\bibitem{Feistauer1}
L.~Dubcova, M.~Feistauer, J.~Horacek, and P.~Svacek.
\newblock {Numerical simulation of interaction between turbulent flow and a
  vibrating airfoil}.
\newblock {\em Computing and Visualization in Science}, 12:207--225, 2009.

\bibitem{Dubiner}
M.~Dubiner.
\newblock Spectral methods on triangles and other domains.
\newblock {\em Journal of Scientific Computing}, 6:345--390, 1991.

\bibitem{DIM3D}
M.~Dumbser.
\newblock {A Diffuse Interface Method for Complex Three-Dimensional Free
  Surface Flows}.
\newblock {\em Computer Methods in Applied Mechanics and Engineering}.
\newblock DOI: , in press.

\bibitem{ADERNSE}
M.~Dumbser.
\newblock Arbitrary high order {PNPM} schemes on unstructured meshes for the
  compressible {Navier--Stokes} equations.
\newblock {\em Computers \& Fluids}, 39:60--76, 2010.

\bibitem{DIM2D}
M.~Dumbser.
\newblock A simple two-phase method for the simulation of complex free surface
  flows.
\newblock {\em Computer Methods in Applied Mechanics and Engineering},
  200:1204--–1219, 2011.

\bibitem{Dumbser2008}
M.~Dumbser, D.~Balsara, E.F. Toro, and C.D. Munz.
\newblock A unified framework for the construction of one--step finite--volume
  and discontinuous {Galerkin} schemes.
\newblock {\em Journal of Computational Physics}, 227:8209--8253, 2008.

\bibitem{ADERVRMHD}
M.~Dumbser and D.S. Balsara.
\newblock High--order unstructured one-step {PNPM} schemes for the viscous and
  resistive {MHD} equations.
\newblock {\em CMES - Computer Modeling in Engineering \& Sciences},
  54:301--333, 2009.

\bibitem{ADERNC}
M.~Dumbser, M.~Castro, C.~{Par\'es}, and E.F. Toro.
\newblock {ADER} schemes on unstructured meshes for non-conservative hyperbolic
  systems: Applications to geophysical flows.
\newblock {\em Computers and Fluids}, 38:1731--–1748, 2009.

\bibitem{DumbserEnauxToro}
M.~Dumbser, C.~Enaux, and E.F. Toro.
\newblock Finite volume schemes of very high order of accuracy for stiff
  hyperbolic balance laws.
\newblock {\em Journal of Computational Physics}, 227:3971--4001, 2008.

\bibitem{USFORCE2}
M.~Dumbser, A.~Hidalgo, M.~Castro, C.~Par\'es, and E.F. Toro.
\newblock {FORCE} schemes on unstructured meshes {II}: Non--conservative
  hyperbolic systems.
\newblock {\em Computer Methods in Applied Mechanics and Engineering},
  199:625--647, 2010.

\bibitem{DumbserKaeser06b}
M.~Dumbser and M.~K\"aser.
\newblock Arbitrary high order non-oscillatory finite volume schemes on
  unstructured meshes for linear hyperbolic systems.
\newblock {\em Journal of Computational Physics}, 221:693--723, 2007.

\bibitem{DumbserKaeser07}
M.~Dumbser, M.~K\"aser, V.A Titarev, and E.F. Toro.
\newblock Quadrature-free non-oscillatory finite volume schemes on unstructured
  meshes for nonlinear hyperbolic systems.
\newblock {\em Journal of Computational Physics}, 226:204--243, 2007.

\bibitem{OsherUniversal}
M.~Dumbser and E.~F. Toro.
\newblock On universal {Osher}--type schemes for general nonlinear hyperbolic
  conservation laws.
\newblock {\em Communications in Computational Physics}, 10:635--671, 2011.

\bibitem{OsherNC}
M.~Dumbser and E.~F. Toro.
\newblock A simple extension of the {Osher} {Riemann} solver to
  non-conservative hyperbolic systems.
\newblock {\em Journal of Scientific Computing}, 48:70--88, 2011.

\bibitem{Dumbser2012}
M.~Dumbser, A.~Uuriintsetseg, and O.~Zanotti.
\newblock {On Arbitrary--Lagrangian--Eulerian One--Step WENO Schemes for Stiff
  Hyperbolic Balance Laws}.
\newblock {\em Communications in Computational Physics}, 14:301--327, 2013.

\bibitem{Faltinsen1978}
{O. M.} Faltinsen.
\newblock A numerical nonlinear method of sloshing in tanks with
  two--dimensional flow.
\newblock {\em J. Ship Res.}, 22:193 -- 202, 1978.

\bibitem{FaltinsenTimokha2000}
{O. M.} Faltinsen, {O. F.} Rognebakke, and {I. A.}~Lukovsky nd~{A. N.}~Timokha.
\newblock Adaptive multimodal approach to nonlinear sloshing in a rectangular
  tank.
\newblock {\em J. Fluid Mech.}, 407:201 -- 234, 2000.

\bibitem{FaltinsenTimokha2001}
{O. M.} Faltinsen and {A. N.} Timokha.
\newblock Adaptive multimodal approach to nonlinear sloshing in a rectangular
  tank.
\newblock {\em J. Fluid Mech.}, 432:167 -- 200, 2001.

\bibitem{FedkiwEtAl1}
R.~Fedkiw, T.~Aslam, B.~Merriman, and S.~Osher.
\newblock A non-oscillatory {Eulerian} approach to interfaces in multimaterial
  flows (the ghost fluid method).
\newblock {\em Journal of Computational Physics}, 152:457--492, 1999.

\bibitem{FedkiwEtAl2}
R.P. Fedkiw, T.~Aslam, and S.~Xu.
\newblock The {Ghost Fluid} method for deflagration and detonation
  discontinuities.
\newblock {\em Journal of Computational Physics}, 154:393--427, 1999.

\bibitem{Feistauer3}
M.~Feistauer, J.~Horacek, M.~Ruzicka, and P.~Svacek.
\newblock Numerical analysis of flow-induced nonlinear vibrations of an airfoil
  with three degrees of freedom.
\newblock {\em Computers and Fluids}, 49:110--127, 2011.

\bibitem{Feistauer2}
M.~Feistauer, V.~Kucera, J.~Prokopova, and J.~Horacek.
\newblock {The ALE discontinuous Galerkin method for the simulatio of air flow
  through pulsating human vocal folds}.
\newblock {\em AIP Conference Proceedings}, 1281:83--86, 2010.

\bibitem{SPHLagrange}
A.~Ferrari, M.~Dumbser, E.F. Toro, and A.~Armanini.
\newblock {A New Stable Version of the SPH Method in Lagrangian Coordinates}.
\newblock {\em Communications in Computational Physics}, 4:378--404, 2008.

\bibitem{SPH3D}
A.~Ferrari, M.~Dumbser, E.F. Toro, and A.~Armanini.
\newblock {A new 3D parallel SPH scheme for free surface flows}.
\newblock {\em Computers \& Fluids}, 38:1203--1217, 2009.

\bibitem{Dambreak3D}
A.~Ferrari, L.~Fraccarollo, M.~Dumbser, E.F. Toro, and A.~Armanini.
\newblock Three--dimensional flow evolution after a dambreak.
\newblock {\em Journal of Fluid Mechanics}, 663:456--477, 2010.

\bibitem{FerrariLevelSet}
A.~Ferrari, C.D. Munz, and B.~Weigand.
\newblock A high order sharp interface method with local timestepping for
  compressible multiphase flows.
\newblock {\em Communications in Computational Physics}, 9:205--230, 2011.

\bibitem{friedrich}
O.~Friedrich.
\newblock Weighted essentially non-oscillatory schemes for the interpolation of
  mean values on unstructured grids.
\newblock {\em Journal of Computational Physics}, 144:194--212, 1998.

\bibitem{GassnerDumbserMunz}
G.~Gassner, M.~Dumbser, F.~Hindenlang, and C.D. Munz.
\newblock Explicit one--step time discretizations for discontinuous {Galerkin}
  and finite volume schemes based on local predictors.
\newblock {\em Journal of Computational Physics}, 230:4232--4247, 2011.

\bibitem{eno}
A.~Harten, B.~Engquist, S.~Osher, and S.~Chakravarthy.
\newblock Uniformly high order essentially non-oscillatory schemes, {III}.
\newblock {\em Journal of Computational Physics}, 71:231--303, 1987.

\bibitem{ALE1996FV}
{R.W.} Healy and {T.F.} Russel.
\newblock Solution of the advection-dispersion equation in two dimensions by a
  finite-volume eulerian-lagrangian localized adjoint method.
\newblock {\em Advances in Water Resources}, 21:11--26, 1998.

\bibitem{HidalgoDumbser}
A.~Hidalgo and M.~Dumbser.
\newblock {ADER} schemes for nonlinear systems of stiff
  advection–diffusion–reaction equations.
\newblock {\em Journal of Scientific Computing}, 48:173--189, 2011.

\bibitem{Hill2003}
{D. F.} Hill.
\newblock Transient and steady--state amplitudes of forced waves in rectangular
  basins.
\newblock {\em Phys. Fluid}, 15:1576 -- 1587, 2003.

\bibitem{Hirt1974}
C.~Hirt, A.~Amsden, and J.~Cook.
\newblock An arbitrary lagrangian–eulerian computing method for all flow
  speeds.
\newblock {\em Journal of Computational Physics}, 14:227–253, 1974.

\bibitem{HirtNichols}
C.~W. Hirt and B.~D. Nichols.
\newblock Volume of fluid ({VOF}) method for dynamics of free boundaries.
\newblock {\em Journal of Computational Physics}, 39:201--225, 1981.

\bibitem{HuShuTri}
C.~Hu and {C.W.} Shu.
\newblock Weighted essentially non-oscillatory schemes on triangular meshes.
\newblock {\em Journal of Computational Physics}, 150:97--127, 1999.

\bibitem{HuangQiul2011}
{C.S.} Huang, T.~Arbogast, and J.~Qiu.
\newblock An eulerian-lagrangian weno finite volume scheme for advection
  problems.
\newblock {\em Journal of Computational Physics}, 231:4028--4052, 2012.

\bibitem{PFEM1}
S.~R. Idelsohn, E.~O{\~{n}}ate, and F.~Del Pin.
\newblock {The Particle Finite Element Method: a powerful tool to solve
  incompressible flows with free-surfaces and breaking waves}.
\newblock {\em International Journal for Numerical Methods in Engineering},
  61:964--984, 2004.

\bibitem{PFEM5}
S.R. Idelsohn, M.~Mier-Torrecilla, and E.~O{\~{n}}ate.
\newblock {Multi--fluid flows with the Particle Finite Element Method}.
\newblock {\em Comput. Methods Appl. Mech. Engrg.}, 198:2750--2767, 2009.

\bibitem{shu_efficient_weno}
{G.S.} Jiang and {C.W.} Shu.
\newblock Efficient implementation of weighted {ENO} schemes.
\newblock {\em Journal of Computational Physics}, 126:202--228, 1996.

\bibitem{orth-basis}
G.~E. Karniadakis and S.~J. Sherwin.
\newblock {\em Spectral/hp Element Methods in CFD}.
\newblock Oxford University Press, 1999.

\bibitem{kaeserjcp}
M.~K\"aser and A.~Iske.
\newblock {ADER} schemes on adaptive triangular meshes for scalar conservation
  laws.
\newblock {\em Journal of Computational Physics}, 205:486--508, 2005.

\bibitem{kurganovtadmor}
A.~Kurganov and E.~Tadmor.
\newblock Solution of two-dimensional {Riemann} problems for gas dynamics
  without {Riemann} problem solvers.
\newblock {\em Numer. Methods Partial Differential Equations}, 18:584--608,
  2002.

\bibitem{PFEM4}
A.~Larese, R.~Rossi, E.~O{\~{n}}ate, and S.R. Idelsohn.
\newblock {Validation of the Particle Finite Element Method (PFEM) for
  Simulation of the Free-Surface Flows}.
\newblock {\em Engineering Computations}, 25:385--425, 2008.

\bibitem{laxwendroff}
{P.D.} Lax and B.~Wendroff.
\newblock Systems of conservation laws.
\newblock {\em Communications in Pure and Applied Mathematics}, 13:217--237,
  1960.

\bibitem{LentineEtAl2011}
M.~Lentine, J\'on~T\'omas Gr\'etarsson, and R.~Fedkiw.
\newblock An unconditionally stable fully conservative semi-lagrangian method.
\newblock {\em Journal of Computational Physics}, 230:2857--2879, 2011.

\bibitem{chengshu2}
W.~Liu, J.~Cheng, and C.W. Shu.
\newblock {High order conservative Lagrangian schemes with Lax–Wendroff type
  time discretization for the compressible Euler equations}.
\newblock {\em Journal of Computational Physics}, 228:8872--8891, 2009.

\bibitem{Loehner}
R.~L\"ohner, C.~Yang, and E.~Onate.
\newblock On the simulation of flows with violent free surface motion.
\newblock {\em Computer Methods in Applied Mechanics and Engineering},
  195:5597--5620, 2006.

\bibitem{MaireCyl2}
R.~Loub\`ere, P.H. Maire, and M.~Shashkov.
\newblock {ReALE: A Reconnection Arbitrary--Lagrangian–-Eulerian method in
  cylindrical geometry}.
\newblock {\em Computers and Fluids}, 46:59--69, 2011.

\bibitem{StagLag}
R.~Loub\`ere, P.H. Maire, and P.~V\'achal.
\newblock {A second--order compatible staggered Lagrangian hydrodynamics scheme
  using a cell--centered multidimensional approximate Riemann solver}.
\newblock {\em Procedia Computer Science}, 1:1931--1939, 2010.

\bibitem{Maire2011}
P.-H. Maire.
\newblock A high-order one-step sub-cell force-based discretization for
  cell-centered lagrangian hydrodynamics on polygonal grids.
\newblock {\em Computers and Fluids}, 46(1):341--347, 2011.

\bibitem{Maire2010}
P.-H. Maire.
\newblock A unified sub-cell force-based discretization for cell-centered
  lagrangian hydrodynamics on polygonal grids.
\newblock {\em International Journal for Numerical Methods in Fluids},
  65:1281–1294, 2011.

\bibitem{MaireCyl1}
P.H. Maire.
\newblock {A high-order cell--centered Lagrangian scheme for compressible fluid
  flows in two--dimensional cylindrical geometry }.
\newblock {\em Journal of Computational Physics}, 228:6882--6915, 2009.

\bibitem{Maire2009}
{P.H.} Maire.
\newblock A high-order cell-centered lagrangian scheme for two-dimensional
  compressible fluid flows on unstructured meshes.
\newblock {\em Journal of Computational Physics}, 228:2391 -- 2425, 2009.

\bibitem{Maire2007}
P.H. Maire, R.~Abgrall, J.~Breil, and J.~Ovadia.
\newblock A cell-centered lagrangian scheme for two-dimensional compressible
  flow problems.
\newblock {\em SIAM Journal on Scientific Computing}, 29:1781–1824, 2007.

\bibitem{Maire2009b}
{P.H.} Maire and B.~Nkonga.
\newblock {Multi-scale Godunov-type method for cell-centered discrete
  Lagrangian hydrodynamics}.
\newblock {\em Journal of Computational Physics}, 228:799--821, 2009.

\bibitem{DLMtheory}
G.~Dal Maso, P.G. LeFloch, and F.~Murat.
\newblock Definition and weak stability of nonconservative products.
\newblock {\em J. Math. Pures Appl.}, 74:483–--548, 1995.

\bibitem{Monaghan1994}
J.J. Monaghan.
\newblock Simulating free surface flows with {SPH}.
\newblock {\em Journal of Computational Physics}, 110:399--406, 1994.

\bibitem{Munoz2007}
M.L. {Mu\~noz} and C.~Par\'es.
\newblock Godunov method for nonconservative hyperbolic systems.
\newblock {\em Mathematical Modelling and Numerical Analysis}, 41:169--185,
  2007.

\bibitem{levelset2}
W.~Mulder, S.~Osher, and J.A. Sethian.
\newblock Computing interface motion in compressible gas dynamics.
\newblock {\em Journal of Computational Physics}, 100:209--228, 1992.

\bibitem{munz94}
C.D. Munz.
\newblock {On Godunov--type schemes for Lagrangian gas dynamics}.
\newblock {\em SIAM Journal on Numerical Analysis}, 31:17--42, 1994.

\bibitem{MurroneGuillard}
A.~Murrone and H.~Guillard.
\newblock A five equation reduced model for compressible two phase flow
  problems.
\newblock {\em Journal of Computational Physics}, 202:664--698, 2005.

\bibitem{Nakayama1980}
T.~Nakayama and K.~Washizu.
\newblock Nonlinear analysis of liquid motion in a container subjected to
  forced pitching oscillation.
\newblock {\em Int. J. Numer. Method Eng.}, 15:1207 -- 1220, 1980.

\bibitem{Okamoto1990}
T.~Okamoto and M.~Kawahara.
\newblock Two--dimensional sloshing analysis by lagrangian finite element
  method.
\newblock {\em Int. Numer. Method Fluid}, 11:453 -- 477, 1990.

\bibitem{Okamoto1997}
T.~Okamoto and M.~Kawahara.
\newblock 3d sloshing analysis by an arbitrary lagrangian–-eulerian finite
  element method.
\newblock {\em Int. J. Comput. Fluid Dyn.}, 8:129 -- 146, 1997.

\bibitem{PFEM6}
E.~O{\~{n}}ate, M.~Celigueta, S.~Idelsohn, F.~Salazar, and B.~Suarez.
\newblock {Possibilities of the Particle Finite Element Method for
  fluid--soil--structure interaction problems}.
\newblock {\em Journal of Computational Mechanics}, 48:307--318, 2011.

\bibitem{PFEM3}
E.~O{\~{n}}ate, S.R. Idelsohn, M.A. Celigueta, and R.~Rossi.
\newblock {Advances in the Particle Finite Element Method for the Analysis of
  Fluid-Multibody Interaction and Bed Erosion in Free-surface Flows}.
\newblock {\em Computer Methods in Applied Mechanics and Engineering},
  197:1777--1800, 2008.

\bibitem{scovazzi1}
A.~L\'opez Ortega and G.~Scovazzi.
\newblock {A geometrically--conservative, synchronized, flux--corrected remap
  for arbitrary Lagrangian--Eulerian computations with nodal finite elements}.
\newblock {\em Journal of Computational Physics}, 230:6709--6741, 2011.

\bibitem{levelset1}
S.~Osher and J.A. Sethian.
\newblock Fronts propagating with curvature--dependent speed: Algorithms based
  on {Hamilton--Jacobi} formulations.
\newblock {\em Journal of Computational Physics}, 79:12--49, 1988.

\bibitem{osherandsolomon}
S.~Osher and F.~Solomon.
\newblock Upwind difference schemes for hyperbolic conservation laws.
\newblock {\em Math. Comput.}, 38:339--374, 1982.

\bibitem{Pares2006}
C.~Par\'es.
\newblock Numerical methods for nonconservative hyperbolic systems: a
  theoretical framework.
\newblock {\em SIAM Journal on Numerical Analysis}, 44:300--–321, 2006.

\bibitem{Pares2004}
C.~Par\'es and M.J. Castro.
\newblock On the well-balance property of roe's method for nonconservative
  hyperbolic systems. applications to shallow-water systems.
\newblock {\em Mathematical Modelling and Numerical Analysis}, 38:821--852,
  2004.

\bibitem{Peery2000}
J.S. Peery and D.E. Carroll.
\newblock Multi-material ale methods in unstructured grids,.
\newblock {\em Computer Methods in Applied Mechanics and Engineering},
  187:591--619, 2000.

\bibitem{PFEM2}
F.~Del Pin, S.~R. Idelsohn, E.~O{\~{n}}ate, and R.~Aubry.
\newblock {The ALE/Lagrangian Particle Finite Element Method: A new approach to
  computation of free-surface flows and fluid-object interactions}.
\newblock {\em Computers and Fluids}, 36:27--38, 2007.

\bibitem{QuiShu2011}
Jing-Mei Qiu and Chi-Wang Shu.
\newblock Conservative high order semi-lagrangian finite difference weno
  methods for advection in incompressible flow.
\newblock {\em Journal of Computational Physics}, 230:863--889, 2011.

\bibitem{Rhebergen2008}
S.~Rhebergen, O.~Bokhove, and J.J.W. van~der Vegt.
\newblock Discontinuous {Galerkin} finite element methods for hyperbolic
  nonconservative partial differential equations.
\newblock {\em Journal of Computational Physics}, 227:1887--–1922, 2008.

\bibitem{Rieber}
M.~Rieber and A.~Frohn.
\newblock A numerical study on the mechanism of splashing.
\newblock {\em International Journal of Heat and Fluid Flow}, 20:455--461,
  1999.

\bibitem{ALE2000Belgium}
K.~Riemslagh, J.~Vierendeels, and E.~Dick.
\newblock An arbitrary lagrangian-eulerian finite-volume method for the
  simulation of rotary displaecment pump flow.
\newblock {\em Applied Numerical Mathematics}, 32:419--433, 2000.

\bibitem{SaurelAbgrall}
R.~Saurel and R.~Abgrall.
\newblock A multiphase {Godunov} method for compressible multifluid and
  multiphase flows.
\newblock {\em Journal of Computational Physics}, 150:425--467, 1999.

\bibitem{Schwendeman}
D.W. Schwendeman, C.W. Wahle, and A.K. Kapila.
\newblock The {Riemann} problem and a high-resolution {Godunov} method for a
  model of compressible two-phase flow.
\newblock {\em Journal of Computational Physics}, 212:490--526, 2006.

\bibitem{scovazzi2}
G.~Scovazzi.
\newblock {Lagrangian shock hydrodynamics on tetrahedral meshes: A stable and
  accurate variational multiscale approach}.
\newblock {\em Journal of Computational Physics}, 231:8029--8069, 2012.

\bibitem{Shao2012}
{J. R.} Shao, {H. Q.} Li, {G. R.} Liu, and {M. B.} Liu.
\newblock An improved sph method for modeling liquid sloshing dynamics.
\newblock {\em Computers and Structures}, 101:18 -- 26, 2012.

\bibitem{Smagorinsky1963}
S.~Smagorinsky.
\newblock General circulation experiments with the primitive equations.
\newblock {\em Mon. Weather Rev.}, 91:99 -- 164, 1963.

\bibitem{Smith1999}
{R.W.} Smith.
\newblock {AUSM(ALE)}: a geometrically conservative arbitrary
  lagrangian--eulerian flux splitting scheme.
\newblock {\em Journal of Computational Physics}, 150:268–286, 1999.

\bibitem{sonar}
T.~Sonar.
\newblock On the construction of essentially non-oscillatory finite volume
  approximations to hyperbolic conservation laws on general triangulations:
  polynomial recovery, accuracy and stencil selection.
\newblock {\em Computer Methods in Applied Mechanics and Engineering},
  140:157--181, 1997.

\bibitem{stroud}
{A.H.} Stroud.
\newblock {\em Approximate Calculation of Multiple Integrals}.
\newblock Prentice-Hall Inc., Englewood Cliffs, New Jersey, 1971.

\bibitem{TianToro}
B.~Tian, E.F. Toro, and C.E. Castro.
\newblock A path-conservative method for a five-equation model of two-phase
  flow with an hllc-type riemann solver.
\newblock {\em Computers and Fluids}, 46:122--132, 2011.

\bibitem{MixedWENO2D}
V.A. Titarev, P.~Tsoutsanis, and D.~Drikakis.
\newblock {WENO schemes for mixed--element unstructured meshes}.
\newblock {\em Communications in Computational Physics}, 8:585--609, 2010.

\bibitem{TokarevaToro}
S.A. Tokareva and E.F. Toro.
\newblock Hllc--type riemann solver for the baer–nunziato equations of
  compressible two-phase flow.
\newblock {\em Journal of Computational Physics}, 229:3573--3604, 2010.

\bibitem{toro-book}
E.F. Toro.
\newblock {\em {Riemann} Solvers and Numerical Methods for Fluid Dynamics}.
\newblock Springer, second edition, 1999.

\bibitem{USFORCE}
E.F. Toro, A.~Hidalgo, and M.~Dumbser.
\newblock {FORCE} schemes on unstructured meshes {I}: Conservative hyperbolic
  systems.
\newblock {\em Journal of Computational Physics}, 228:3368--3389, 2009.

\bibitem{Toumi1992}
I.~Toumi.
\newblock A weak formulation of {Roe's} approximate {Riemann} solver.
\newblock {\em Journal of Computational Physics}, 102:360--373, 1992.

\bibitem{MixedWENO3D}
P.~Tsoutsanis, V.A. Titarev, and D.~Drikakis.
\newblock {WENO schemes on arbitrary mixed-element unstructured meshes in three
  space dimensions}.
\newblock {\em Journal of Computational Physics}, 230:1585--1601, 2011.

\bibitem{leer5}
B.~{van Leer}.
\newblock Towards the ultimate conservative difference scheme {V}: A second
  order sequel to {Godunov}'s method.
\newblock {\em Journal of Computational Physics}, 32:101--136, 1979.

\bibitem{Verhagen1965}
{H. G.} Verhagen and L.~Wijingaarden.
\newblock Nonlinear oscillation of fluid in a container.
\newblock {\em J. Fluid. Mech.}, 22:737 -- 751, 1965.

\bibitem{Neumann1950}
J.~von Neumann and R.D. Richtmyer.
\newblock A method for the calculation of hydrodynamics shocks.
\newblock {\em Journal of Applied Physics}, 21:232--237, 1950.

\bibitem{Wu1998}
{G. X.} Wu, {Q. A.} Ma, and {R. E.} Taylor.
\newblock Numerical simulation of sloshing waves in a 3d tank based on a finite
  element method.
\newblock {\em Appl. Ocean Res.}, 20:337 -- 355, 1998.

\bibitem{ZhangShu3D}
{Y.T.} Zhang and {C.W.} Shu.
\newblock Third order {WENO} scheme on three dimensional tetrahedral meshes.
\newblock {\em Communications in Computational Physics}, 5:836--848, 2009.

\bibitem{LNGsloshing2012}
W.~Zhao, J.~Yang, Z.~Hu, and L.~Xiao.
\newblock Experimental investigation of effects of inner--tank sloshing on
  hydrodynamics of an {FLNG} system.
\newblock {\em Journal of Hydrodynamics}, 24:107 -- 115, 2012.

\end{thebibliography}
\bibliographystyle{plain}

\end{document}